\documentclass[12pt]{article}
\usepackage{tikz,makeidx}
\usepackage{amssymb, xr, multicol,textcomp,amsthm,amsmath,url,enumerate,hyperref}
\usepackage[english]{babel} 
\def\hfl#1{\smash{\mathop{\hbox to 10mm{\rightarrowfill}}\limits^{\textstyle
#1}}}
\makeindex
\newtheorem{proposition}[equation]{Proposition}
\newtheorem{corollary}[equation]{Corollary}
\newtheorem{theorem}[equation]{Theorem}
\newtheorem{exa}[equation]{Example}
\newtheorem{ex}[equation]{Exercise}
\newtheorem{s-ex}[equation]{Side-exercise}
\newtheorem{exas}[equation]{Examples}
\newtheorem{lemma}[equation]{Lemma}

\newtheorem{remar}[equation]{Remark}
\newtheorem{remars}[equation]{Remarks}
\newtheorem{nota}[equation]{Notation}
\newtheorem{sremar}[equation]{Side-remark}
\newtheorem{definitio}[equation]{Definition}

\newenvironment{remark}{\begin{remar} \rm }{\end{remar}}

\newenvironment{example}{\begin{exa} \rm }{\end{exa}}

\newcommand{\tilth}{d_e}

\newcommand{\CaC}{{\cal C}}
\newcommand{\CD}{{\cal D}}

\newcommand{\CJ}{{\cal J}}

\newcommand{\Link}{L}

\newcommand{\halfbra}{{|}}

\newcommand{\ZZ}{\mathbb{Z}}
\newcommand{\RR}{\mathbb{R}}
\newcommand{\QQ}{\mathbb{Q}}
\newcommand{\CC}{\mathbb{C}}

\newcommand{\Asc}{{\cal A}}
\newcommand{\Bdesc}{{\cal B}}
\newcommand{\ballB}{B}
\newcommand{\Euhler}{E}
\newcommand{\fMorse}{f}

\newcommand{\metrigo}{\frak{g}} 
\newcommand{\twocycG}{G}
\newcommand{\handlebodA}{H_{\cal A}}
\newcommand{\handA}{H_{\cal A}}
\newcommand{\handB}{H_{\cal B}}
\newcommand{\homotop}{h}
\newcommand{\linh}{{\cal L}}
\newcommand{\matchingo}{\frak{m}}
\newcommand{\matchm}{m}
\newcommand{\nornu}{\nu}
\newcommand{\pointp}{p}
\newcommand{\pointq}{q}
\newcommand{\twist}{\frak{t}}
\newcommand{\normbun}{\frak{V}}
\newcommand{\rpt}{\mathbb{RP}^3}
\newcommand{\pointw}{w}
\newcommand{\ComX}{X}
\newcommand{\ComY}{Y}
\newcommand{\ComZ}{Z}
\newcommand{\bp}{\noindent {\sc Proof: }}
\newcommand{\eop}{\nopagebreak \hspace*{\fill}{$\Box$} \medskip}

\begin{document}
\title{A combinatorial definition of the $\Theta$-invariant from {H}eegaard diagrams}
\author{Christine Lescop \thanks{Institut Fourier, UJF Grenoble, CNRS}}
\maketitle
\begin{abstract}
The invariant $\Theta$ is the simplest $3$-manifold invariant defined by counting graph configurations. It is actually an invariant of rational homology $3$-spheres $M$ equipped with a combing $\ComX$ over the complement of a point.
The invariant $\Theta(M,X)$ is the sum of $6 \lambda(M)$ and $\frac{p_1(\ComX)}{4}$,
where $\lambda$ denotes the Casson-Walker invariant, and $p_1$ is an invariant of combings that is an extension of a first relative Pontrjagin class, and that is simply related to a Gompf invariant $\theta_G$. In \cite{lesHC}, we proved a combinatorial formula for the $\Theta$-invariant in terms of decorated Heegaard diagrams. In this article,
we study the variations of the invariants $p_1$ or $\theta_G$ when the decorations of the Heegaard diagrams that define the combings change, independently.
Then we prove that the formula of \cite{lesHC} defines an invariant of combed once punctured rational homology $3$-spheres without referring to configuration spaces. Finally, we prove that this invariant is the sum of $6 \lambda(M)$ and $\frac{p_1(\ComX)}{4}$ for integral homology spheres, by proving surgery formulae both for the combinatorial invariant and for $p_1$.
\end{abstract}

\tableofcontents

\section{Introduction}
\label{secintro}

In this article, a {\em $\QQ$--sphere\/} or {\em rational homology sphere\/} is a smooth closed oriented $3$-manifold that has the same rational homology as $S^3$.

\subsection{General introduction}
\label{subgenintro}

The work of Witten \cite{witten} pioneered the introduction of many $\QQ$--sphere invariants, among which the Le-Murakami-Ohtsuki universal finite type invariant \cite{lmo} and the Kontsevich configuration space invariant \cite{ko} that was proved to be equivalent to the LMO invariant for integral homology spheres by
G. Kuperberg and D. Thurston \cite{kt}.
The construction of the Kontsevich configuration space invariant for a $\QQ$--sphere $M$ involves a point $\infty$ in $M$,
an identification of a neighborhood of $\infty$ with a neighborhood $S^3 \setminus B(1)$ of $\infty$ in $S^3=\RR^3 \cup \{\infty\}$, and a parallelization
$\tau$ of $(\check{M}=M \setminus \{\infty\} )$ that coincides with the standard parallelization of $\RR^3$ on $\RR^3 \setminus B(1)$, where $B(r)$ denotes the ball centered at $0$ with radius $r$ in $\RR^3$.
The Kontsevich configuration space invariant is in fact an invariant of $(M,\tau)$.
Its degree one part
$\Theta(M,\tau)$ is
the sum of $6 \lambda(M)$ and $\frac{p_1(\tau)}{4}$, where $\lambda$ is the Casson-Walker invariant and $p_1$ is a Pontrjagin number associated with $\tau$, according to a theorem of G. Kuperberg and D. Thurston \cite{kt} generalized to rational homology spheres in \cite{lessumgen}.
Here, the Casson-Walker invariant $\lambda$ is normalized like in \cite{akmc,gm,mar} for integral homology spheres, and like $\frac{1}{2}\lambda_W$ for rational homology spheres where $\lambda_W$ is the Walker normalisation in \cite{wal}.

Let $B_M$ denote the complement in $M$ of the neighborhood of $\infty$ identified with $S^3 \setminus B(1)$, $B_M$ is a rational homology ball.
An \emph{$\infty$-combing} of such a rational homology sphere $M$ is a section
of the unit tangent bundle $U\check{M}$ of $\check{M}$ that is constant on $\check{M} \setminus B_M$ (via the identifications above with $\RR^3 \setminus B(1)$ near $\infty$), up to homotopies through this kind of sections. 
As it is shown in \cite{lesHC}, $\Theta(M,.)$ is actually an invariant of rational homology spheres equipped with such an $\infty$-combing.

Every closed oriented $3$--manifold $M$ can be written as the union of two handlebodies $\handA$ and $\handB$ glued along their common boundary that is a genus $g$ surface as
\index{HA@$\handA$}
$$M = \handA \cup_{\partial \handA} \handB$$
where $\partial \handA= -\partial \handB$.
Such a decomposition is called a {\em Heegaard decomposition\/} of $M$.
A {\em system of meridians\/} for $\handA$ is a system of $g$ disjoint curves $\alpha_i$ of $\partial \handA$ that bound disjoint disks $D(\alpha_i)$ properly embedded in $\handA$ such that the union of the $\alpha_i$ does not separate $\partial \handA$. For a positive integer $g$, we will denote the set $\{1,2,\dots,g\}$ by $\underline{g}$.
Let $(\alpha_i)_{i\in \underline{g}}$ be a system of meridians for $\handA$ and let $(\beta_j)_{j\in \underline{g}}$ be such a system for $\handB$. Then the surface equipped with the collections of the curves $\alpha_i$ and the curves $\beta_j=\partial D(\beta_j)$ determines $M$.
When the collections $(\alpha_i)_{i\in \underline{g}}$ and
$(\beta_j)_{j\in \underline{g}}$ are transverse, the data
$\CD=(\partial \handA, (\alpha_i)_{i\in \underline{g}}, (\beta_j)_{j\in \underline{g}})$ is called a {\em Heegaard diagram.\/}

Such a Heegaard diagram may be obtained from a Morse function $\fMorse_M$ of $M$ that has one minimum mapped to $(-3)$, one maximum mapped to $9$, $g$ index one points $a_i$ and $g$ index $2$ points $b_j$, such that $\fMorse_M$ maps index $1$ points to $1$ and index $2$ points to $5$, and $\fMorse_M$ satisfies generic Morse-Smale conditions ensuring transversality of descending and ascending manifolds of critical points, with respect to a Euclidean metric $\metrigo$ of $M$.
Thus the surface $\partial \handA$ is $\fMorse_M^{-1}(3)$, the ascending manifolds of the $a_i$ intersect $\handA$ as disks $D(\alpha_i)$ bounded by the $\alpha_i$ and the descending manifolds of the $b_j$ intersect $\handB$ as disks $D(\beta_j)$ bounded by the $\beta_i$, and the flow line closures from
$a_i$ to $b_j$ are in natural one-to-one correspondence with the crossings of $\alpha_i \cap \beta_j$.
Conversely, for any Heegaard diagram, there exists a Morse function $\fMorse_M$ with the properties above.

A {\em matching\/} in a genus $g$ Heegaard diagram \index{Dcal@$\CD$}$$\CD=(\partial \handA,\{\alpha_i\}_{i=1,\dots,g},\{\beta_j\}_{j=1,\dots,g})$$ is a set $\matchingo$ of $g$ crossings such that every curve of the diagram contains one crossing of $\matchingo$. 
An \emph{exterior point} in such a diagram $\CD$ is a point of $\partial \handA \setminus \left(\coprod_{i=1}^g \alpha_i \cup \coprod_{j=1}^g \beta_j\right)$. The choice of a matching $\matchingo$ and of an exterior point $\pointw$ in a diagram $\CD$ of $M$ equips ${M}$ with the following $\infty$-combing $\ComX(\pointw,\matchingo )=\ComX(\CD,\pointw,\matchingo )$. 

Remove an open ball around a flow line from the minimum to the maximum that goes through $\pointw$, so that we are left with a rational homology ball $$B_M(2)=B_M \cup_{\partial B(1)=\partial B_M} B(2) \setminus \mathring{B}(1)$$ where the gradient field of $\fMorse_M$ is vertical near the boundary.
Reversing the gradient field along the flow lines $\gamma(c)$ associated with the crossings $c$ of $\matchingo$ as in Subsection~\ref{subcombX} below
produces the $\infty$-combing $\ComX(\pointw,\matchingo )$ of $M$. 

Let $\theta_G$ denote the invariant of combings of rational homology spheres introduced by Gompf in \cite[Section 4]{Go}.
A choice of a standard modification described in Subsection~\ref{subclscomb} of $\ComX(\pointw,\matchingo )$ in the fixed neighborhood of $\infty$ identified with $S^3 \setminus B(2)$ transforms $\ComX(\pointw,\matchingo )$ into a combing $\ComX(M,\pointw,\matchingo)$ such that $p_1(\ComX(\pointw,\matchingo)) - \theta_G(\ComX(M,\pointw,\matchingo))$ is independent of $(M,\pointw,\matchingo)$.

In \cite[Theorem~1.5]{lesHC}, we express $\Theta(M,\ComX(\pointw,\matchingo))$ as a combination
\index{Theta@$\tilde{\Theta}$}$$\tilde{\Theta}(\CD,\pointw,\matchingo)=\ell_2(\CD)+s_{\ell}(\CD,\matchingo)-e(\CD,\pointw,\matchingo)$$
of invariants of Heegaard diagrams $\CD$ equipped with a matching $\matchingo$ and an exterior point $\pointw$. First combinatorial expressions of the ingredients $\ell_2(\CD)$, $s_{\ell}(\CD,\matchingo)$, and $e(\CD,\pointw,\matchingo)$ are given in the end of this introduction section whereas Section~\ref{seccombdefmore} provides alternative expressions and properties of these quantities.

In this article, we give several expressions of the variations of $p_1(\ComX(\pointw,\matchingo))$, or, equivalently of $\theta_G(\ComX(M,\pointw,\matchingo))$, when $\pointw$ and $\matchingo$ vary, for a fixed Heegaard diagram.
Expressions in terms of linking numbers are given in Subsection~\ref{subpone} and derived combinatorial expressions can be found
in Section~\ref{secpXWCP}.

The latter ones allow us to give combinatorial proofs that 
$\left(4\tilde{\Theta}(\CD,\pointw,\matchingo)-{p_1(\ComX(\pointw,\matchingo))}\right)$
is independent of $(\pointw,\matchingo)$ in Section~\ref{secvartilWCP}. 
We prove that $$\tilde{\lambda}(\CD)=\frac1{24}\left(4\tilde{\Theta}(\CD,\pointw,\matchingo)-{p_1(\ComX(\pointw,\matchingo))}\right)$$ only depends on the presented rational homology sphere $M$,
combinatorially, in Section~\ref{secinvcomb}. We set $\tilde{\lambda}(M)=\tilde{\lambda}(\CD)$ so that $\tilde{\lambda}$ is a topological invariant of rational homology $3$-spheres.

Then we give a direct combinatorial proof
that $\tilde{\lambda}$ satisfies the Casson surgery formula for $\frac1{n}$--Dehn surgeries along null-homologous knots in Section~\ref{secsurgcomb}.
This implies that $\tilde{\lambda}$ coincides with the Casson invariant for integral homology $3$-spheres.
Our proof also yields a surgery formula for $p_1$ that is stated in Theorem~\ref{thmponesurg}.

Thus this article contains an independent construction of the Casson invariant that includes a direct proof of the Casson surgery formula,
and an independent combinatorial proof of the formula of \cite[Theorem~3.8]{lesHC} for the $\Theta$-invariant in terms of Heegaard diagrams in the case of $\ZZ$-spheres. It also describes the behaviour of the four quantities $\ell_2(\CD)$, $s_{\ell}(\CD,\matchingo)$, $e(\CD,\pointw,\matchingo)$ and $p_1(\ComX(\CD,\pointw,\matchingo))$ (or equivalently $\theta_G(\ComX(M,\pointw,\matchingo))$) associated with Heegaard diagrams $\CD$ decorated with $(\pointw,\matchingo)$ under standard modifications of Heegaard diagrams, and Dehn surgeries. These quantities might show up in combinatorial definitions of other invariants from Heegaard diagrams, as $\theta_G$ that Gripp and Huang use to define the Heegaard Floer homology $\widehat{HF}$ grading in \cite{GH}.

The definitions introduced in \cite[Theorem~3.8]{lesHC} are also given here, causing overlaps.
I thank Jean-Mathieu Magot for useful conversations.

\subsection{Conventions and notations}
Unless otherwise mentioned, all manifolds are oriented. Boundaries are oriented by the outward normal first convention.
Products are oriented by the order of the factors. More generally, unless otherwise mentioned, the order of appearance of coordinates or parameters orients chains or manifolds.
The normal bundle $\normbun(A)$ of an oriented submanifold $A$ is oriented so that the normal bundle followed by the tangent bundle of the submanifold induce the orientation of the ambient manifold (fiberwise).
The transverse intersection of two submanifolds $A$ and $B$ of a manifold $C$ is oriented so that the normal bundle $\normbun_x(A\cap B)$ of $A\cap B$ at $x$ is oriented as $(\normbun_x(A) \oplus \normbun_x(B))$.
When the dimensions of two such submanifolds add up to the dimension of $C$, each intersection point is equipped with a sign $\pm 1$ that is $1$ if and only if $(\normbun_x(A) \oplus \normbun_x(B))$ (or equivalently $(T_x(A) \oplus T_x(B))$) induces the orientation of $C$. When $A$ is compact, the sum of the signs of the intersection points is the \emph{algebraic intersection number} $\langle A,B \rangle_{C}$.
The \emph{linking number} $lk(L_1,L_2)=lk_C(L_1,L_2)$ of two disjoint null-homologous cycles $L_1$ and $L_2$ of respective dimensions $d_1$ and $d_2$ in an oriented $(d_1+d_2+1)$-manifold $C$ is the algebraic intersection $\langle L_1 , W_2 \rangle_{C}$ of $L_1$ with a chain $W_2$ bounded by $L_2$ in $C$. 
This definition extends to rationally null-homologous cycles by bilinearity.

\subsection{Introduction to the combinatorial definition of \texorpdfstring{$\tilde{\Theta}$}{Theta} }
\setcounter{equation}{0}
\label{subseccombdef}

In the end of this section, we give explicit formulas for the ingredients $\ell_2(\CD)$, $s_{\ell}(\CD,\matchingo)$ and $e(\CD,\pointw,\matchingo)$ in the formula
$$\tilde{\Theta}(\CD,\pointw,\matchingo)=\ell_2(\CD)+s_{\ell}(\CD,\matchingo)-e(\CD,\pointw,\matchingo)$$
for a Heegaard diagram $\CD$ equipped with a matching $\matchingo$ and an exterior point $\pointw$. These ingredients will be studied in more details in Section~\ref{seccombdefmore}.

Let $\CD=(\partial \handA, (\alpha_i)_{i\in \underline{g}}, (\beta_j)_{j\in \underline{g}})$ be a Heegaard diagram of a rational homology $3$-sphere. A {\em crossing \/} $c$ of $\CD$ is an intersection point of a curve $\alpha_{i(c)}=\alpha(c)$ and a curve $\beta_{j(c)}=\beta(c)$. Its sign $\sigma(c)$ is $1$ if $\partial \handA$ is oriented by the oriented tangent vector of $\alpha(c)$ followed by the oriented tangent vector of $\beta(c)$ at $c$ as above. It is $(-1)$ otherwise.
The set of crossings of $\CD$ is denoted by $\CaC$.

Let
$$[\CJ_{ji}]_{(j,i) \in \underline{g}^2}=[\langle \alpha_i,\beta_j \rangle_{\partial \handA}]^{-1}$$
denote the inverse matrix of the intersection matrix.
$$\sum_{i=1}^g\CJ_{ji}\langle \alpha_i,\beta_k \rangle_{\partial \handA}=\delta_{jk}=\left\{\begin{array}{ll} 1 \;& \mbox{if} \; j=k\\ 0 & \mbox{otherwise.}\end{array}\right.$$

When $d$ and $e$ are two crossings of $\alpha_i$, $[d,e]_{\alpha_i}=[d,e]_{\alpha}$ denotes the set of crossings from $d$ to $e$ (including them) along $\alpha_i$, or the closed arc from $d$ to $e$ in $\alpha_i$ depending on the context. Then $[d,e[_{\alpha}=[d,e]_{\alpha} \setminus \{e\}$,
$]d,e]_{\alpha}=[d,e]_{\alpha} \setminus \{d\}$ and $]d,e[_{\alpha}=[d,e[_{\alpha} \setminus \{d\}$.

Now, for such a part $I$ of $\alpha_i$, $$\langle I,\beta_j \rangle=\langle I,\beta_j \rangle_{\partial \handA}= \sum_{c\in I \cap \beta_j} \sigma(c).$$
We shall also use the notation $\halfbra$ for ends of arcs to say that an end is half-contained in an arc, and that it must be counted with coefficient $1/2$.
(``$[d,e\halfbra_{\alpha}=[d,e]_{\alpha} \setminus \{e\}/2$''). We agree that $\halfbra d,d\halfbra_{\alpha}=\emptyset$.

We use the same notation for arcs $[d,e\halfbra_{\beta_j}=[d,e\halfbra_{\beta}$ of $\beta_j$.
For example, if $d$ is a crossing of $\alpha_i\cap\beta_j$, then
$$\langle [d,d\halfbra_{\alpha},\beta_j \rangle=\frac{\sigma(d)}{2}$$
and
$$\langle [c,d\halfbra_{\alpha},[e,d\halfbra_{\beta} \rangle= \frac{\sigma(d)}{4} + \sum_{c\in [c,d[_{\alpha} \cap [e,d[_{\beta}} \sigma(c).$$

\begin{example}\label{exaheegone}
 In the Heegaard diagrams of $\rpt$ in Figure~\ref{figheegrpt},
$\langle [c,c\halfbra_{\alpha}, [c,c\halfbra_{\beta}\rangle=\frac14, 
\langle [c,c\halfbra_{\alpha}, [c,d\halfbra_{\beta}\rangle=\langle [c,d\halfbra_{\alpha}, [c,c\halfbra_{\beta}\rangle=\frac12, \langle [c,d\halfbra_{\alpha}, [c,d\halfbra_{\beta}\rangle=\frac54,\langle [c,c\halfbra_{\alpha}, \beta_1\rangle=\frac12$ and $\langle [c,d\halfbra_{\alpha}, \beta_1\rangle=\frac32.$
 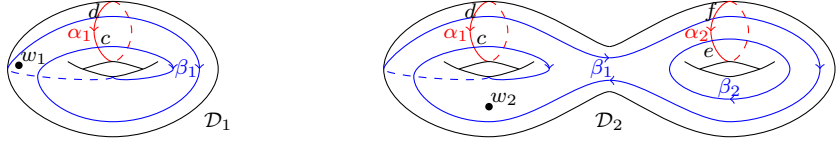
\begin{figure}\begin{center}
\begin{tikzpicture} \useasboundingbox (.2,-1) rectangle (12,1); 
\draw [red,->]  (1.88,.45) node[left]{\scriptsize $\alpha_1$}  (2,.9) .. controls (1.9,.9) and (1.75,.7) .. (1.75,.5);  
\draw  [red] (1.75,.5) .. controls (1.75,.3) and (1.9,.1) .. (2,.1);
\draw [red,dashed] (2,.9) .. controls (2.1,.9) and (2.25,.7) .. (2.25,.5) .. controls (2.25,.3) and (2.1,.1) .. (2,.1);
\draw plot[smooth] coordinates{(1.4,.1) (1.6,0) (2,-.1) (2.4,0) (2.6,.1)};
\draw plot[smooth] coordinates{(1.6,0) (2,.1) (2.4,0)};
\fill (.75,.05) circle (0.05);
\draw (.95,.15) node{\scriptsize $\pointw_1$};
\draw (.6,0) .. controls (.6,.45) and (1.3,.9)  .. (2,.9) .. controls (2.7,.9) and (3.4,.45) .. (3.4,0) .. controls  (3.4,-.45) and (2.7,-.9) .. (2,-.9) .. controls (1.3,-.9) and (.6,-.45)  ..  (.6,0);
\draw [blue,->]  (.6,0) .. controls (.6,.1) and (1.2,.7)  .. (2,.7) .. controls (2.6,.7) and (3.15,.35) .. (3.15,0);
\draw  [blue,->]    (2.98,0) node{\scriptsize $\beta_1$} (3.15,0) .. controls  (3.15,-.35) and (2.6,-.7) .. (2,-.7) .. controls (1.4,-.7) and (1,-.35)  ..  (1,-.1)  .. controls (1,.1) and (1.5,.3)  .. (2,.3) .. controls (2.3,.3) and (2.8,.15) .. (2.8,0) ;
\draw  [blue] (2.8,0) .. controls  (2.8,-.1) and (2,-.2) .. (2,-.1);
\draw [blue,dashed] (2,-.1)  .. controls (2,-.2) and (.6,-.1) ..(.6,0);
\draw  (3.4,-.7) node{\scriptsize $\CD_1$}  (1.75,.78) node{\scriptsize $d$}   (1.9,.4) node{\scriptsize $c$};
\begin{scope}[xshift=5cm] 
\draw [red,->]  (1.88,.45) node[left]{\scriptsize $\alpha_1$}  (2,.9) .. controls (1.9,.9) and (1.75,.7) .. (1.75,.5);  
\draw  [red] (1.75,.5) .. controls (1.75,.3) and (1.9,.1) .. (2,.1);
\draw [red,dashed] (2,.9) .. controls (2.1,.9) and (2.25,.7) .. (2.25,.5) .. controls (2.25,.3) and (2.1,.1) .. (2,.1);
\draw plot[smooth] coordinates{(1.4,.1) (1.6,0) (2,-.1) (2.4,0) (2.6,.1)};
\draw plot[smooth] coordinates{(1.6,0) (2,.1) (2.4,0)};
\fill (2,-.5) circle (0.05);
\draw (2.2,-.4) node{\scriptsize $\pointw_2$};
\draw (.6,0) .. controls (.6,.45) and (1.3,.9)  .. (2,.9) .. controls (2.7,.9) and (3.4,.3) .. (3.6,.3) (3.6,-.3) .. controls  (3.4,-.3) and (2.7,-.9) .. (2,-.9) .. controls (1.3,-.9) and (.6,-.45)  ..  (.6,0);
\draw [blue,->]  (.6,0) .. controls (.6,.1) and (1.2,.7)  .. (2,.7) .. controls (2.6,.7) and (3.15,.15) .. (3.6,.15) ;
\draw [blue,->]  (3.6,.15) .. controls (4.05,.15) and (4.7,.7) .. (5.2,.7) .. controls (5.7,.7) and (6.4,.4) .. (6.4,0) ;
\draw  [blue,->]  (6.4,0) .. controls (6.4,-.4) and  (5.7,-.7) .. (5.2,-.7) ..  controls (4.7,-.7) and (4.05,-.15) .. (3.6,-.15);
\draw  [blue,->]    (3.5,0) node{\scriptsize $\beta_1$} (3.6,-0.15) .. controls  (3.15,-.15) and (2.6,-.7) .. (2,-.7) .. controls (1.4,-.7) and (1,-.35)  ..  (1,-.1)  .. controls (1,.1) and (1.5,.3)  .. (2,.3) .. controls (2.3,.3) and (2.8,.15) .. (2.8,0) ;
\draw  [blue] (2.8,0) .. controls  (2.8,-.1) and (2,-.2) .. (2,-.1);
\draw [blue,dashed] (2,-.1)  .. controls (2,-.2) and (.6,-.1) ..(.6,0);
\draw    (1.75,.78) node{\scriptsize $d$}   (1.9,.4) node{\scriptsize $c$};
\draw  (3.6,-.7) node{\scriptsize $\CD_2$};
\end{scope}
\begin{scope}[xshift=8.2cm] 
\draw (.4,.3) .. controls (.6,.3) and (1.3,.9)  .. (2,.9) .. controls (2.7,.9) and (3.4,.45) .. (3.4,0) .. controls  (3.4,-.45) and (2.7,-.9) .. (2,-.9) .. controls (1.3,-.9) and (.6,-.3)  ..  (.4,-.3);
\draw [red,->]  (1.88,.45) node[left]{\scriptsize $\alpha_2$}  (2,.9) .. controls (1.9,.9) and (1.75,.7) .. (1.75,.5);  
\draw  [red] (1.75,.5) .. controls (1.75,.3) and (1.9,.1) .. (2,.1);
\draw [red,dashed] (2,.9) .. controls (2.1,.9) and (2.25,.7) .. (2.25,.5) .. controls (2.25,.3) and (2.1,.1) .. (2,.1);
\draw plot[smooth] coordinates{(1.4,.1) (1.6,0) (2,-.1) (2.4,0) (2.6,.1)};
\draw plot[smooth] coordinates{(1.6,0) (2,.1) (2.4,0)};
\draw [blue,->] (2,-.25) node{\scriptsize $\beta_2$} (2,-.4) .. controls (1.6,-.4) and (1.2,-.2)  ..  (1.2,0)  .. controls (1.2,.2) and (1.6,.4)  .. (2,.4) .. controls (2.4,.4) and (2.8,.2) .. (2.8,0) .. controls (2.8,-.2) and (2.4,-.4)  .. (2,-.4);
\draw (1.75,.76) node{\scriptsize $f$}   (1.72,.22) node{\scriptsize $e$};
\end{scope}
\end{tikzpicture}
\caption{Two Heegaard diagrams of $\rpt$}
\label{figheegrpt}\end{center}
\end{figure}
\end{example}

\subsection{First combinatorial definitions of \texorpdfstring{$\ell_2$}{l_2} and \texorpdfstring{$s_{\ell}(\CD,\matchingo)$}{sl(D,P)}}

Choose a matching $\matchingo=\{\matchm_i;i \in \underline{g}\}$ where $\matchm_i \in \alpha_{\rho^{-1}(i)} \cap \beta_i$, for a permutation $\rho$ of $\underline{g}$.
For two crossings $c$ and $d$ of $\CaC$, set
$$\tilde{\ell}_{\matchingo}(c,d)=\langle \halfbra \matchm_{\rho(i(c))},c\halfbra_{\alpha}, \halfbra \matchm_{j(d)},d\halfbra_{\beta}\rangle -\sum_{(i,j) \in \underline{g}^2}\CJ_{ji} \langle \halfbra  \matchm_{\rho(i(c))},c\halfbra_{\alpha},\beta_j \rangle \langle \alpha_i ,\halfbra \matchm_{j(d)},d\halfbra_{\beta} \rangle.$$

Then \index{ltwo@$\ell_2(\CD)$}
$$ \ell_2(\CD) = \sum_{(c,d) \in \CaC^2}\CJ_{j(c)i(d)}\CJ_{j(d)i(c)}\sigma(c)\sigma(d) \tilde{\ell}_{\matchingo}(c,d) -\sum_{c \in \CaC} \CJ_{j(c)i(c)}\sigma(c)\tilde{\ell}_{\matchingo}(c,c) $$
and \index{sl@$s_{\ell}(\CD,\matchingo)$}
$$ s_{\ell}(\CD,\matchingo) =\sum_{(c,d) \in \CaC^2}\CJ_{j(c)i(c)}\CJ_{j(d)i(d)}\sigma(c)\sigma(d) \tilde{\ell}_{\matchingo}(c,d).$$

\begin{example}
\label{exaheegtwo}
For the genus one Heegaard diagram $\CD_1$ of Figure~\ref{figheegrpt}, $\sigma(c)=1$, $\langle \alpha_1,\beta_1 \rangle_{\partial \handlebodA}=2$, $\CJ_{11}=\frac12$, choose $\{c\}$ as a matching, 
$\tilde{\ell}_{\{c\}}(c,c)=\tilde{\ell}_{\{c\}}(c,d)=\tilde{\ell}_{\{c\}}(d,c)=0,\tilde{\ell}_{\{c\}}(d,d)=\frac12-\CJ_{11}=0$ so that $\ell_2(\CD_1)=s_{\ell}(\CD_1,\{c\})=0$.
 
For the genus two Heegaard diagram $\CD_2$ of Figure~\ref{figheegrpt},  $\CJ_{11}=\frac12=-\CJ_{21}$, $\CJ_{22}=1$ and $\CJ_{12}=0$, choose $\{c,e\}$ as a matching. For any crossing $x$ of $\CD_2$, $$0=\tilde{\ell}_{\{c,e\}}(c,x)=\tilde{\ell}_{\{c,e\}}(x,c)=\tilde{\ell}_{\{c,e\}}(e,x)=\tilde{\ell}_{\{c,e\}}(x,e)=\tilde{\ell}_{\{c,e\}}(d,d),$$
$$\begin{array}{lll}\tilde{\ell}_{\{c,e\}}(f,f)&=&\frac14-\frac34\CJ_{11}-\frac14\CJ_{12}-\frac34\CJ_{21}-\frac14\CJ_{22}=0\\
\tilde{\ell}_{\{c,e\}}(d,f)&=&\frac34-\frac32\CJ_{11}-\frac12\CJ_{12}=0\\
\tilde{\ell}_{\{c,e\}}(f,d)&=&-\frac12\CJ_{11}-\frac12\CJ_{21}=0\end{array}$$
so that $\ell_2(\CD_2)=s_{\ell}(\CD_2,\{c,e\})=0$.
\end{example}

\subsection{Combinatorial definition of $e(\CD,\pointw,\matchingo)$}
\label{subdefe}

\begin{figure}[h]
\begin{center}
\begin{tikzpicture}
\useasboundingbox (-.5,-.2) rectangle (8.9,1.9);
\draw [red,->] (1,.5) arc (-90:270:.4);
\draw [red] (1,.2) node{\scriptsize $\alpha^{\prime}_1$};
\draw [red,->] (2.8,.5) arc (270:-90:.4);
\draw [red] (2.8,.2) node{\scriptsize $\alpha^{\prime\prime}_1$} (4.2,.9) node{$\dots$};
\draw [red,->] (5.6,.5) arc (-90:270:.4);
\draw [red] (5.6,.2) node{\scriptsize $\alpha^{\prime}_g$};
\draw [red,->] (7.4,.5) arc (270:-90:.4);
\draw [red] (7.4,.2) node{\scriptsize $\alpha^{\prime\prime}_g$};
\draw (-.5,-.1) -- (-.5,1.8) -- (8.9,1.8) -- (8.9,-.1) -- (-.5,-.1);
\draw [blue,thick,->] (2.55,.9) -- (2.25,.9);
\draw [blue,thick,->] (1.55,.9) -- (1.25,.9);
\draw [blue,thick,->] (7.15,.9) -- (6.85,.9);
\draw [blue,thick,->] (6.15,.9) -- (5.85,.9);
\draw [blue,thick,->] (3.05,.85) -- (3.35,.85);
\draw [blue,thick,->] (3,.75) -- (3.3,.75);
\draw [blue,thick,->] (3.25,.65) -- (2.95,.65);
\draw [blue,thick,->] (.85,.65) -- (.55,.65);
\draw [blue,thick,->] (.45,.85) -- (.75,.85);
\draw [blue,thick,->] (.5,.75) -- (.8,.75);
\draw [blue,thick,->] (5.1,.75) -- (5.4,.75);
\draw [blue,thick,->] (5.05,.85) -- (5.35,.85);
\draw [blue,thick,->] (5.35,.95) -- (5.05,.95);
\draw [blue,thick,->] (5.1,1.05) -- (5.4,1.05);
\draw [blue,thick,->] (7.6,.75) -- (7.9,.75);
\draw [blue,thick,->] (7.65,.85) -- (7.95,.85);
\draw [blue,thick,->] (7.95,.95) -- (7.65,.95);
\draw [blue,thick,->] (7.6,1.05) -- (7.9,1.05);
\draw [blue,thick,dotted] (2.25,.9) .. controls (1.9,.9) and (2.5,1.5) .. (2.8,1.5) .. controls (4,.9) and (4.9,1.05) .. (5.1,1.05) (7.9,1.05).. controls (8.4,1.05) and (8.2,1.7) .. (2.8,1.7) .. controls (2,1.7) and (1.9,.9) .. (1.55,.9);
\draw (1.6,.7) node{\scriptsize  $\matchm_1$} (2.2,.7) node{\scriptsize  $\matchm_1$} (6.2,.7) node{\scriptsize  $\matchm_g$} (6.8,.7) node{\scriptsize  $\matchm_g$} (8.6,.1) node{\scriptsize  $R_{\CD}$};
\end{tikzpicture}
\caption{The Heegaard surface cut along the $\alpha_i$ and deprived of a neighborhood of $\pointw$}
\label{figplansplit}
\end{center}
\end{figure}
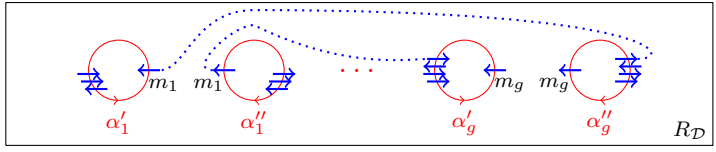

Let $\pointw$ be an exterior point of $\CD$.
The choice of $\matchingo$ being fixed, represent the Heegaard diagrams in a plane by removing from $\partial \handA$ a disk around $\pointw$ that does not intersect the diagram curves, and by cutting the surface $\partial \handA$ along the $\alpha_i$. Each $\alpha_i$ gives rise to two copies $\alpha^{\prime}_i$ and $\alpha^{\prime}_i$ of $\alpha_i$ that bound disjoint disks with opposite orientations in the plane. Locate the crossing $\matchm_i$ at the points with upward tangent vectors of $\alpha^{\prime}_i$ and $\alpha^{\prime\prime}_i$, and locate the other crossings near 
the points with downward tangent vectors as in Figure~\ref{figplansplit}. Draw the arcs of the curves $\beta_j$ so that they have horizontal tangent vectors near the crossings.

The rectangle has the standard parallelization of the plane. Then there is a map ``unit tangent vector'' from each partial projection of a beta curve $\beta_j$ in the plane to $S^1$. The total degree of this map for the curve $\beta_j$ is denoted by \index{de@$\tilth$}$\tilth(\beta_{j})$.
For a crossing $c \in \beta_j$, $\tilth(\halfbra \matchm_j,c \halfbra_{\beta}) \in \frac12 \ZZ$ denotes the degree of the restriction of this map to the arc $\halfbra \matchm_j,c \halfbra_{\beta}$.
For any $c \in \CaC$, define
\index{dec@$\tilth(c)$}$$\tilth(c)=\tilth(\halfbra \matchm_{j(c)},c \halfbra_{\beta})
-\sum_{(r,s) \in \underline{g}^2}\CJ_{sr}\langle \alpha_r,\halfbra \matchm_{j(c)},c \halfbra_{\beta}\rangle\tilth(\beta_{s}).$$

Set
\index{eDmw@$e(\CD,\pointw,\matchingo)$}$$e(\CD,\pointw,\matchingo)=\sum_{c \in \CaC} \CJ_{j(c)i(c)}\sigma(c)\tilth(c)$$
and \index{Theta@$\tilde{\Theta}$}$$\tilde{\Theta}(\CD,\pointw,\matchingo)=\ell_2(\CD)+s_{\ell}(\CD,\matchingo)-e(\CD,\pointw,\matchingo).$$

\begin{example}
\label{exaheegfour}
For the rectangular diagram of $(\CD_1,\{c\},\pointw_1)$ of Figure~\ref{figrecgenonet}, $\tilth(\halfbra c,c \halfbra_{\beta})=0$ and $\tilth(c)=0$, $\tilth(\halfbra c,d \halfbra_{\beta})=\frac12$, $\tilth(\beta_1)=2$ 
so that $\tilth(d)=-\frac12$, $e(\CD_1,\pointw_1,\{c\})=-\frac14$ and $\tilde{\Theta}(\CD_1,\pointw_1,\{c\})=\frac14$.
 
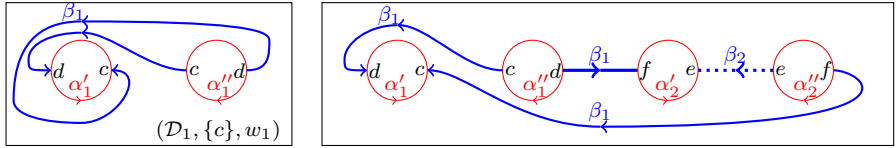
\begin{figure}[h]
\begin{center}
\begin{tikzpicture}
\useasboundingbox (0,-.2) rectangle (11.6,1.9);
\draw (3.8,.1) node[left]{\scriptsize $(\CD_1,\{c\},\pointw_1)$};
\draw [red,->] (1,.5) arc (-90:270:.4);
\draw [red] (1,.7) node{\scriptsize $\alpha^{\prime}_1$};
\draw  (.7,.9) node{\scriptsize $d$} (1.3,.9)node{\scriptsize $c$}  (3.1,.9) node{\scriptsize $d$} (2.5,.9)node{\scriptsize $c$};
\draw [blue]  (.9,1.68) node{\scriptsize $\beta_1$};
\draw [red,->] (2.8,.5) arc (270:-90:.4);
\draw [red] (2.85,.7) node{\scriptsize $\alpha^{\prime\prime}_1$};
\draw (0,-.1) -- (0,1.8) -- (3.8,1.8) -- (3.8,-.1) -- (0,-.1);
\draw [blue,thick,->]  (2.4,.9) .. controls (1.9,.9) and (1.4,1.4) .. (1,1.4);
\draw [blue,thick,->]  (1,1.4) .. controls (.6,1.4) and (.3,1.25) .. (.3,1.05) .. controls (.3,.9) .. (.6,.9);
\draw [blue,thick,->]  (3.2,.9) .. controls (3.5,.9) and (3.5,1.1) .. (3.5,1.3) .. controls (3.5,1.45) and (3,1.55) .. (1,1.55);
\draw [blue,thick,->]  (1,1.55) .. controls (.6,1.55) and (.1,1.4) .. (.1,.5) .. controls (.1,.3) and (.5,.2) .. (1,.2)  .. controls (1.5,.2) and (1.9,.9) .. (1.4,.9);
\begin{scope}[xshift=4.2cm]
\draw [red,->] (1,.5) arc (-90:270:.4);
\draw  (.7,.9) node{\scriptsize $d$} (1.3,.9)node{\scriptsize $c$}  (3.1,.9) node{\scriptsize $d$} (2.5,.9)node{\scriptsize $c$};
\draw [blue] (.9,1.63) node{\scriptsize $\beta_1$};
\draw [red,->] (2.8,.5) arc (270:-90:.4);
\draw  [red]  (1,.7) node{\scriptsize $\alpha^{\prime}_1$} (2.85,.7) node{\scriptsize $\alpha^{\prime\prime}_1$};
\draw (0,-.1) -- (0,1.8) -- (7.6,1.8) -- (7.6,-.1) -- (0,-.1);
\draw [blue,thick,->]  (2.4,.9) .. controls (1.9,.9) and (1.4,1.5) .. (1,1.5);
\draw [blue,thick,->]  (1,1.5) .. controls (.6,1.5) and (.3,1.25) .. (.3,1.05) .. controls (.3,.9) .. (.6,.9);
\draw [blue,very thick,->]  (3.7,1.1) node{\scriptsize $\beta_1$}   (4.2,.9) -- (3.7,.9) (3.2,.9) -- (3.7,.9);
\draw [blue,thick,->] (3.7,.35) node{\scriptsize $\beta_1$}  (6.8,.9) .. controls (7.2,.9) and (8.2,.15) .. (3.7,.15);
\draw [blue,thick,->] (3.7,.15) .. controls (2.6,.15) and (1.9,.9) .. (1.4,.9);
\end{scope}
\begin{scope}[xshift=7.8cm]
\draw [red,->] (1,.5) arc (-90:270:.4);
\draw  (.7,.9) node{\scriptsize $f$} (1.3,.9)node{\scriptsize $e$}  (3.1,.9) node{\scriptsize $f$} (2.5,.9)node{\scriptsize $e$};
\draw [blue]  (1.9,1.1) node{\scriptsize $\beta_2$};
\draw [red,->] (2.8,.5) arc (270:-90:.4);
\draw [red] (1,.7) node{\scriptsize $\alpha^{\prime}_2$}(2.85,.7) node{\scriptsize $\alpha^{\prime\prime}_2$};
\draw [blue,very thick,dotted,->]  (1.4,.9) -- (1.9,.9) (2.4,.9) -- (1.9,.9);
\end{scope}
\end{tikzpicture}
\caption{Rectangular diagrams of $(\CD_1,\{c\},\pointw_1)$ and $(\CD_2,\{c,e\},\pointw_2)$}
\label{figrecgenonet}
\end{center}
\end{figure}

For the rectangular diagram of $(\CD_2,\{c,e\},\pointw_2)$ of Figure~\ref{figrecgenonet}, $\tilth(c)=\tilth(e)=\tilth(\beta_1)=\tilth(\beta_2)=0$, $\tilth(d)=\tilth(f)=\frac12$,
$e(\CD_2,\pointw_2,\{c\})=\frac14$ and $\tilde{\Theta}(\CD_2,\pointw_2,\{c\})=-\frac14$.
\end{example}

\section{More on the combinatorial definition of \texorpdfstring{$\tilde{\Theta}$}{Theta} }
\setcounter{equation}{0}
\label{seccombdefmore}

In this section, we show that the quantities $\ell_2(\CD)$, $s_{\ell}(\CD,\matchingo)$ and $e(\CD,\pointw,\matchingo)$ defined in the previous section for a Heegaard diagram $\CD$ equipped with a matching $\matchingo$ and an exterior point $\pointw$ only depend on their arguments (e.g. on $\CD$, for $\ell_2(\CD)$ \dots) and not on extra data used to define them like numberings or orientations of the diagram curves. We also give alternative definitions of
$\ell_2(\CD)$ and $s_{\ell}(\CD,\matchingo)$.

\subsection{More on $e(\CD,\pointw,\matchingo)$}

Recall the notation and definitions of Subsection~\ref{subdefe} with respect to a fixed matching $\matchingo=\{\matchm_i;i \in \underline{g}\}$ where $\matchm_i \in \alpha_{\rho^{-1}(i)} \cap \beta_i$, for a permutation $\rho$ of $\underline{g}$.

\begin{lemma}
\label{leminvthetaW}
$$e(\CD,\pointw,\matchingo)=\sum_{c \in \CaC} \CJ_{j(c)i(c)}\sigma(c)\tilth(c)$$
 does not depend on our specific way of drawing the diagram with our conventions. It only depends on $\CD$, $\pointw$ and $\matchingo$.
\end{lemma}

The topological interpretation of $e(\CD,\pointw,\matchingo)$ as an Euler class given in Corollary~\ref{coreultilth} yields a conceptual proof of this lemma. We nevertheless give a purely combinatorial proof below.

We use the Kronecker symbol $\delta_{cd}$ that is $1$ if $c=d$ and $0$ otherwise.
We first prove the following lemma.

\begin{lemma}
\label{leminvthetaWprel}
A full positive twist of a curve $\alpha^{\prime}_i$ or a curve $\alpha^{\prime\prime}_i$ in Figure~\ref{figplansplit} changes $\tilth(c)$ to
$\tilth(c)+\frac12 \delta_{i(c)i}-\frac12\delta_{\rho(i)j(c)}$.
\end{lemma}
\bp
When a crossing is moved counterclockwise along a curve $\alpha$, (like along $\alpha^{\prime\prime}_i$ in Figure~\ref{figmovecross}) the degree increases (by $1$ for a full loop) when the crossing enters (the disk bounded by) $\alpha$ in Figure~\ref{figplansplit} and decreases when the crossing goes out.
Furthermore the positive crossings enter $\alpha^{\prime}_i$ and the negative ones enter $\alpha^{\prime\prime}_i$.
Then letting all the crossings make a full positive loop around
$\alpha^{\prime\prime}_i$ (resp. around $\alpha^{\prime}_i$)
changes $\tilth(\beta_s)$ to
$\tilth(\beta_s) - \langle \alpha_i, \beta_s\rangle $
(resp. to $\tilth(\beta_s) + \langle \alpha_i, \beta_s\rangle $).
Now, for a full positive loop around
$\alpha^{\prime\prime}_i$,
$\tilth(\halfbra \matchm_{j(c)},c \halfbra_{\beta})$ is changed to
$$\tilth(\halfbra \matchm_{j(c)},c \halfbra_{\beta}) - \langle \alpha_i,]\matchm_{j(c)},c[_{\beta} \rangle
 -\delta_{i(c)i}\delta_{(-1)\sigma(c)}\sigma(c)-\delta_{\rho(i)j(c)}\delta_{1\sigma(\matchm_{j(c)})}\sigma(\matchm_{j(c)}).$$
Indeed, right before $c$, $\beta_{j(c)}$ hits $\alpha^{\prime\prime}_{i}$ iff $\sigma(c)=-1$ and $i(c)=i$.
Similarly, after $\matchm_{j(c)}$, $\beta_{j(c)}$ exits
$\alpha^{\prime\prime}_{i}$ iff $\sigma(c)=1$ and $\rho^{-1}(j(c))=i$.
This expression can be rewritten as
$$\tilth(\halfbra \matchm_{j(c)},c \halfbra_{\beta}) - \langle \alpha_i,\halfbra \matchm_{j(c)},c\halfbra_{\beta} \rangle
 +\frac12 \delta_{i(c)i}-\frac12\delta_{\rho(i)j(c)}.$$
Similarly, for a full positive loop around
$\alpha^{\prime}_i$,
$\tilth(\halfbra \matchm_{j(c)},c \halfbra_{\beta})$ is changed to
$$\tilth(\halfbra \matchm_{j(c)},c \halfbra_{\beta}) + \langle \alpha_i,\halfbra \matchm_{j(c)},c\halfbra_{\beta} \rangle
 +\frac12 \delta_{i(c)i}-\frac12\delta_{\rho(i)j(c)}.$$
Now, since
$$\sum_{(r,s) \in \underline{g}^2}\CJ_{sr}\langle \alpha_r,\halfbra \matchm_{j(c)},c \halfbra_{\beta}\rangle\langle \alpha_i,\beta_{s}\rangle=\langle \alpha_i,\halfbra \matchm_{j(c)},c\halfbra_{\beta} \rangle,$$
$\tilth(c)$ is changed to
$\tilth(c)+\frac12 \delta_{i(c)i}-\frac12\delta_{\rho(i)j(c)}$ in both cases.
\eop

\noindent{\sc Proof of Lemma~\ref{leminvthetaW}:}
Note that $e(\CD,\pointw,\matchingo)$ does not depend on the numberings of the diagram curves.
We prove that $e(\CD,\pointw,\matchingo)$ does not depend on our specific way of drawing the diagram with our conventions when the orientations of the diagram curves are fixed.
When the curves $\alpha^{\prime}_i$ and $\alpha^{\prime\prime}_i$ move in the plane without being twisted, the $\tilth(c)$ stay in $\frac12 \ZZ$ and are therefore invariant.
Therefore it suffices to prove that $e(\CD,\pointw,\matchingo)$ is invariant under a full twist of a curve $\alpha^{\prime}_i$ or a curve $\alpha^{\prime\prime}_i$.
Since
$$\sum_{c \in \CaC} \CJ_{j(c)i(c)}\sigma(c)(\delta_{i(c)i}-\delta_{\rho(i)j(c)})=\sum_{c \in \alpha_i} \CJ_{j(c)i}\sigma(c)-\sum_{c \in \beta_{\rho(i)}} \CJ_{\rho(i)i(c)}\sigma(c)=1-1=0,$$
according to Lemma~\ref{leminvthetaWprel},
$e(\CD,\pointw,\matchingo)$ does not vary under these moves. It is not hard to prove that $e(\CD,\pointw,\matchingo)$ does not depend on the orientations of the curves $\beta$.  Changing the orientation of a curve $\alpha$ permutes $\alpha^{\prime}_i$ and $\alpha^{\prime\prime}_i$ and does not modify $e(\CD,\pointw,\matchingo)$ either so that the lemma is proved.
\eop

\subsection{More on \texorpdfstring{$s_{\ell}(\CD,\matchingo)$}{sl(D,P)}}

Fix a point $a_i$ inside each disk $D(\alpha_i)$ and a point $b_j$ inside each disk $D(\beta_j)$.
Then join $a_i$ to each crossing $c$ of $\alpha_i$ by a segment
$[a_i,c]_{D(\alpha_i)}$ oriented from $a_i$ to $c$ in $D(\alpha_i)$, so that these segments only meet at $a_i$ for different $c$. Similarly define segments $[c,b_{j(c)}]_{D(\beta_{j(c)})}$ from $c$ to $b_{j(c)}$ in $D(\beta_{j(c)})$.
Then for each $c$, define the {\em flow line \/}
$\gamma(c) = [a_{i(c)},c]_{D(\alpha_{i(c)})} \cup [c,b_{j(c)}]_{D(\beta_{j(c)})}$. When $\gamma(c)$ is smooth, $\gamma(c)$ is a flow line closure of a Morse function $\fMorse_M$ giving birth to $\CD$ discussed in the introduction.

For each point $a_i$ in the disk $D(\alpha_i)$ as in Subsection~\ref{subgenintro}, choose a point $a_i^+$ and a point $a_i^-$ close to $a_i$ outside $D(\alpha_i)$ so that $a_i^+$ is on the positive side of $D(\alpha_i)$ (the side of the positive normal) and $a_i^-$ is on the negative side of $D(\alpha_i)$.
Similarly fix points $b_j^+$ and $b_j^-$ close to the $b_j$ and outside the $D(\beta_j)$.

Then for a crossing $c \in \alpha_{i(c)} \cap \beta_{j(c)}$, $\gamma(c)_{\parallel}$ will denote the following chain. Consider a small meridian curve $m(c)$ of $\gamma(c)$ on $\partial \handA$, it intersects $\beta_{j(c)}$ at two points: $c_{\Asc}^+$ on the positive side of $D(\alpha_{i(c)})$ and $c_{\Asc}^-$ on the negative side of $D(\alpha_{i(c)})$. The meridian $m(c)$ also intersects $\alpha_{i(c)}$ at $c_{\Bdesc}^+$ on the positive side of $D(\beta_{j(c)})$ and $c_{\Bdesc}^-$ on the negative side of $D(\beta_{j(c)})$.
Let $[c_{\Asc}^+,c_{\Bdesc}^+]$, $[c_{\Asc}^+,c_{\Bdesc}^-]$, $[c_{\Asc}^-,c_{\Bdesc}^+]$ and $[c_{\Asc}^-,c_{\Bdesc}^-]$ denote the four quarters of $m(c)$ with the natural ends and orientations associated with the notation, as in Figure~\ref{figgammapar}.

\begin{figure}[h]
\begin{center}
\begin{tikzpicture}
\useasboundingbox (-4.2,-1.4) rectangle (4.2,1.4);
\begin{scope}[xshift=-3.5cm]
\draw [blue,->] (-.2,1.3) node{\scriptsize $\beta_j$} (0,-1.2) -- (0,1.2);
\draw [red,->] (1.1,.1) node[right]{\scriptsize $\alpha_i$} (-1.2,0) -- (1.2,0);
\draw (0,0) circle (.8);
\draw [->] (0,-.8) arc (-90:-135:.8);
\draw [->] (0,-.8) arc (-90:-45:.8);
\draw [->] (0,.8) arc (90:135:.8);
\draw [->] (0,.8) arc (90:45:.8);
\draw (.2,.2) node{$c$} (-1.8,1.3) node[left]{\scriptsize $\sigma(c)=1$} (-1,-.25) node{\scriptsize $c_{\Bdesc}^-$} (1.05,-.25) node{\scriptsize $c_{\Bdesc}^+$} (.3,-1) node{\scriptsize $c_{\Asc}^-$} (-.45,-.8) node[left]{\tiny $[c_{\Asc}^-,c_{\Bdesc}^-]$} (-.45,.8) node[left]{\tiny $[c_{\Asc}^+,c_{\Bdesc}^-]$} (.45,-.8) node[right]{\tiny $[c_{\Asc}^-,c_{\Bdesc}^+]$} (.45,.8) node[right]{\tiny $[c_{\Asc}^+,c_{\Bdesc}^+]$} (.3,1) node{\scriptsize $c_{\Asc}^+$};
\fill (0,-.8) circle (0.05) (0,.8) circle (0.05) (.8,0) circle (0.05) (-.8,0) circle (0.05) (0,0) circle (0.05);
\end{scope}
\begin{scope}[xshift=3.5cm]
\draw [blue,->] (-.2,-1.3) node{\scriptsize $\beta_j$} (0,1.2) -- (0,-1.2);
\draw [red,->] (1.1,.1) node[right]{\scriptsize $\alpha_i$} (-1.2,0) -- (1.2,0);
\draw (0,0) circle (.8);
\draw [->] (0,-.8) arc (-90:-135:.8);
\draw [->] (0,-.8) arc (-90:-45:.8);
\draw [->] (0,.8) arc (90:135:.8);
\draw [->] (0,.8) arc (90:45:.8);
\draw (.2,.2) node{$c$} (1.8,1.3) node[right]{\scriptsize $\sigma(c)=-1$} (-1,-.25) node{\scriptsize $c_{\Bdesc}^+$} (1.05,-.25) node{\scriptsize $c_{\Bdesc}^-$} (.3,-1) node{\scriptsize $c_{\Asc}^-$} (-.45,-.8) node[left]{\tiny  $[c_{\Asc}^-,c_{\Bdesc}^+]$} (-.45,.8) node[left]{\tiny $[c_{\Asc}^+,c_{\Bdesc}^+]$} (.45,-.8) node[right]{\tiny $[c_{\Asc}^-,c_{\Bdesc}^-]$} (.45,.8) node[right]{\tiny $[c_{\Asc}^+,c_{\Bdesc}^-]$} (.3,1) node{\scriptsize $c_{\Asc}^+$};
\fill (0,-.8) circle (0.05) (0,.8) circle (0.05) (.8,0) circle (0.05) (-.8,0) circle (0.05) (0,0) circle (0.05);
\end{scope}
\end{tikzpicture}
\caption{$m(c)$, $c_{\Asc}^+$, $c_{\Asc}^-$, $c_{\Bdesc}^+$ and $c_{\Bdesc}^-$}
\label{figgammapar}
\end{center}
\end{figure}

Let $\gamma_{\Asc}^+(c)$ (resp. $\gamma_{\Asc}^-(c)$) be an arc parallel to $[a_{i(c)},c]_{D(\alpha_{i(c)})}$ from $a_{i(c)}^+$ to $c_{\Asc}^+$ (resp. from $a_{i(c)}^-$ to $c_{\Asc}^-$) that does not meet $D(\alpha_{i(c)})$.
Let $\gamma_{\Bdesc}^+(c)$ (resp. $\gamma_{\Bdesc}^-(c)$) be an arc parallel to $[c,b_{j(c)}]_{D(\beta_{j(c)})}$ from $c_{\Bdesc}^+$ to $b_{j(c)}^+$ (resp. from $c_{\Bdesc}^-$ to $b_{j(c)}^-$) that does not meet $D(\beta_{j(c)})$.

$$\gamma(c)_{\parallel}=\frac12(\gamma_{\Asc}^+(c)+\gamma_{\Asc}^-(c))
+\frac14([c_{\Asc}^+,c_{\Bdesc}^+]+[c_{\Asc}^+,c_{\Bdesc}^-]+[c_{\Asc}^-,c_{\Bdesc}^+]+[c_{\Asc}^-,c_{\Bdesc}^-])
+\frac12(\gamma_{\Bdesc}^+(c)+\gamma_{\Bdesc}^-(c)).$$

Set $a_{i\parallel}=\frac12(a_i^++a_i^-)$ and $b_{j\parallel}=\frac12(b_j^++b_j^-)$. Then
$\partial \gamma(c)_{\parallel}=b_{j(c)\parallel}-a_{i(c)\parallel}$.

Recall our matching $\matchingo=\{\matchm_i;i \in \underline{g}\}$ where $\matchm_i \in \alpha_{\rho^{-1}(i)} \cap \beta_i$, for a permutation $\rho$ of $\underline{g}$, so that $\gamma_i=\gamma(\matchm_i)$ goes from $a_{\rho^{-1}(i)}$ to $b_i$.

Set \index{LD@$\Link(\CD,\matchingo)$}$$\Link(\CD,\matchingo) = \sum_{i=1}^g \gamma_i- \sum_{c\in \CaC} \CJ_{j(c)i(c)} \sigma(c) \gamma(c).$$

Note that $\Link(\CD,\matchingo)$
is a cycle since
$$ \partial \Link(\CD,\matchingo)=\sum_{i=1}^g(b_i-a_i)- \sum_{(i,j) \in \underline{g}^2} \CJ_{ji}\langle \alpha_i,\beta_j \rangle_{\partial \handA} (b_j-a_i)=0.$$
Set $\Link(\CD,\matchingo)_{\parallel} = \sum_{i=1}^g \gamma_i- \sum_{c\in \CaC} \CJ_{j(c)i(c)} \sigma(c) \gamma(c)_{\parallel}$.

In this subsection, we prove the following proposition. 

\begin{proposition} 
\label{propsllk}
For any Heegaard diagram $\CD$ equipped with a matching $\matchingo$, \index{sl@$s_{\ell}(\CD,\matchingo)$}
 $$s_{\ell}(\CD,\matchingo)=lk(\Link(\CD,\matchingo),\Link(\CD,\matchingo)_{\parallel}).$$
\end{proposition}

This proposition has the following easy corollary.
\begin{corollary}
The real number $s_{\ell}(\CD,\matchingo)$ is an invariant of the Heegaard diagram $\CD$ equipped with $\matchingo$ that does not depend on the orientations and numberings of the curves $\alpha_i$ and $\beta_j$, and that does not change when the roles of the $\alpha$-curves or the $\beta$-curves are permuted.
\end{corollary}
\eop

We first prove the following lemma that will be useful later, too.

\begin{lemma} \label{lemseifsurfgam}
For any curve $\alpha_i$ (resp. $\beta_j$), choose a basepoint $\pointp(\alpha_i)$ (resp. $\pointp(\beta_j)$). These choices being made,
for any crossing $c$ of $\CaC$, define the triangle $T_{\beta}(c)$ in the disk $D(\beta_{j(c)})$ such that
$$ \partial T_{\beta}(c) =[\pointp(\beta(c)),c ]_{\beta} + (\gamma(c) \cap \handB)- (\gamma(\pointp(\beta(c))) \cap \handB).$$
Similarly, define the triangle $T_{\alpha}(c)$ in the disk $D(\alpha_{i(c)})$ such that
$$ \partial T_{\alpha}(c) =-[ \pointp(\alpha(c)),c ]_{\alpha} + (\gamma(c) \cap \handA)- (\gamma(\pointp(\alpha(c))) \cap \handA).$$
 Let $K=\sum_{c \in \CaC}k_{c} \gamma(c)$ be  a cycle of $M$.

Let $\Sigma_T(K)=\sum_{c \in \CaC}k_{c}(T_{\alpha}(c)+T_{\beta}(c))$ and
$$\Sigma_D(K)=\sum_{(i,j,c)\in \underline{g}^2 \times \CaC}\CJ_{ji}k_c\left(\langle \halfbra  \pointp(\alpha(c)),c \halfbra_{\alpha},\beta_j\rangle D(\alpha_i)-\langle \alpha_i,\halfbra \pointp(\beta(c)),c \halfbra_{\beta} \rangle D(\beta_j)\right).$$
There exists a $2$-chain $\Sigma_{\Sigma}(K)$ in $\partial \handA$ whose boundary $ \partial \Sigma_{\Sigma}(K)$ is
$$\sum_{(i,j,c)\in \underline{g}^2 \times \CaC}\CJ_{ji}k_c\left(\langle \alpha_i,\halfbra \pointp(\beta(c)),c \halfbra_{\beta} \rangle\beta_j  -
\langle \halfbra  \pointp(\alpha(c)),c \halfbra_{\alpha},\beta_j\rangle \alpha_i\right) + \sum_{c \in \CaC}k_{c}([ \pointp(\alpha(c)),c ]_{\alpha}-[\pointp(\beta(c)),c ]_{\beta})$$
so that the boundary of $$\Sigma(K) = \Sigma_{\Sigma}(K) + \Sigma_D(K) + \Sigma_T(K)$$ is $K$.
\end{lemma}
Though it is not visible from the notation, the surfaces depend on the basepoints.

\noindent{\sc Proof of Lemma~\ref{lemseifsurfgam}:} $$\partial \Sigma_T(K) -K=\sum_{c \in \CaC}k_{c}([\pointp(\beta(c)),c ]_{\beta}-[ \pointp(\alpha(c)),c ]_{\alpha})$$ is a cycle. 
Any $1$-cycle $\sigma$ of $\partial \handA$ is homologous to 
$\sum_{(i,j)\in \underline{g}^2}\CJ_{ji}(\langle \sigma,\beta_j\rangle \alpha_i  +\langle \alpha_i,\sigma \rangle\beta_j)$. Therefore by pushing $(\partial \Sigma_T(K) -K)$ in the directions of the positive and negative normals to the $\alpha$ and the $\beta$ in $\partial \handA$, and by averaging, we see that
$(K -\partial \Sigma_T(K))$ is homologous in $\partial \handA$ to 
$$\sum_{(i,j,c)\in \underline{g}^2 \times \CaC}\CJ_{ji}k_c\left(\langle \halfbra  \pointp(\alpha(c)),c \halfbra_{\alpha},\beta_j\rangle \alpha_i-\langle \alpha_i,\halfbra \pointp(\beta(c)),c \halfbra_{\beta} \rangle\beta_j\right)$$
that bounds
$$\Sigma_D(K)=\sum_{(i,j,c)\in \underline{g}^2 \times \CaC}\CJ_{ji}k_c\left(\langle \halfbra  \pointp(\alpha(c)),c \halfbra_{\alpha},\beta_j\rangle D(\alpha_i)-\langle \alpha_i,\halfbra \pointp(\beta(c)),c \halfbra_{\beta} \rangle D(\beta_j)\right).$$
\eop

\begin{proposition}
\label{propevallk}
For any curve $\alpha_i$ (resp. $\beta_j$), choose a basepoint $\pointp(\alpha_i)$ (resp. $\pointp(\beta_j)$). These choices being fixed, set 
$$\tilde{\ell}(c,d)=\langle \halfbra \pointp(\alpha(c)),c\halfbra_{\alpha}, \halfbra \pointp(\beta(d)),d\halfbra_{\beta}\rangle -\sum_{(i,j) \in \underline{g}^2}\CJ_{ji} \langle \halfbra \pointp(\alpha(c)),c\halfbra_{\alpha},\beta_j \rangle \langle \alpha_i ,\halfbra \pointp(\beta(d)),d\halfbra_{\beta} \rangle).$$
Let $K=\sum_{c \in \CaC}k_{c} \gamma(c)$ and $L=\sum_{d \in \CaC}g_{d} \gamma(d)$ be two $1$-cycles of $M$. Then $$lk(K, L_{\parallel})=lk(L, K_{\parallel})=\sum_{(c,d) \in \CaC^2}k_{c}g_{d} \tilde{\ell}(c,d).$$
\end{proposition}
\bp
The first equality comes from the symmetry of the linking number and from the observation that $lk(K, L_{\parallel})=lk(K_{\parallel}, L)$. Compute $lk(K,L_{\parallel})$ as the intersection of $L_{\parallel}$ with the surface bounded by $K$ provided by Lemma~\ref{lemseifsurfgam}. Thus $lk(K,L_{\parallel})=\langle \Sigma_{\Sigma}(K), L_{\parallel}\rangle$.
Now, since $L=\sum_{d \in \CaC}g_{d} \gamma(d)$ is a cycle,
$$L=\sum_{d \in \CaC}g_{d} (\gamma(d) -\gamma(\pointp(\beta(d))))$$
and it suffices to prove the result when $L=\gamma(d) -\gamma(\pointp(\beta(d)))$.
For any path $[x,y]$ from a point $x$ to a point $y$ in $\partial \handA$, when $x$ and $y$ are outside  $\partial \Sigma_{\Sigma}(K)$,
$$\langle x-y,\Sigma_{\Sigma}(K)  \rangle_{\Sigma}=\langle [x,y], \partial \Sigma_{\Sigma}(K) \rangle_{\partial \handA}.$$
Thus by averaging,
$$\langle \gamma(d)_{\parallel} -\gamma(\pointp(\beta(d)))_{\parallel}, \Sigma_{\Sigma}(K)\rangle =\langle \partial \Sigma_{\Sigma}(K), \halfbra \pointp(\beta(d)),d \halfbra_{\beta} \rangle_{\partial \handA}.$$
This is $$\begin{array}{ll}\sum_{c \in \CaC}k_{c}\langle \halfbra \pointp(\alpha(c)),c \halfbra_{\alpha}, \halfbra \pointp(\beta(d)),d \halfbra_{\beta}\rangle_{\partial \handA}&\\ -\sum_{(i,j,c)\in \underline{g}^2 \times \CaC}k_c\CJ_{ji}\left(\langle \halfbra  \pointp(\alpha(c)),c \halfbra_{\alpha},\beta_j\rangle \langle \alpha_i,\halfbra \pointp(\beta(d)),d \halfbra_{\beta} \rangle_{\partial \handA}\right)&=\sum_{c \in \CaC}k_{c}\left(\tilde{\ell}(c,d)-\tilde{\ell}(c,\pointp(\beta(d)))\right).
\end{array}$$
\eop

\noindent {\sc Proof of Proposition~\ref{propsllk}:} Apply Proposition~\ref{propevallk} with the basepoints of $\matchingo$ so that $\tilde{\ell}=\tilde{\ell}_{\matchingo}$.
\eop

\subsection{More on \texorpdfstring{$\ell_2(\CD)$}{l_2(D)}}

\begin{proposition}
\label{propevalbis}
For any curve $\alpha_i$ (resp. $\beta_j$), choose a basepoint $\pointp(\alpha_i)$ (resp. $\pointp(\beta_j)$). These choices being made, for two crossings $c$ and $d$ of $\CaC$, set
$$\ell(c,d)=\langle [\pointp(\alpha(c)),c\halfbra_{\alpha}, [\pointp(\beta(d)),d\halfbra_{\beta}\rangle -\sum_{(i,j) \in \underline{g}^2}\CJ_{ji} \langle [\pointp(\alpha(c)),c\halfbra_{\alpha},\beta_j \rangle \langle \alpha_i ,[\pointp(\beta(d)),d\halfbra_{\beta} \rangle$$
and 
$\tilde{\ell}(c,d)=\langle \halfbra \pointp(\alpha(c)),c\halfbra_{\alpha}, \halfbra \pointp(\beta(d)),d\halfbra_{\beta}\rangle -\sum_{(i,j) \in \underline{g}^2}\CJ_{ji} \langle \halfbra \pointp(\alpha(c)),c\halfbra_{\alpha},\beta_j \rangle \langle \alpha_i ,\halfbra \pointp(\beta(d)),d\halfbra_{\beta} \rangle.$
Then, for any $2$--cycle $\twocycG=\sum_{(c,d) \in \CaC^2}g_{c d} (\gamma(c) \times \gamma(d)_{\parallel})$ of $M^2$,
\index{ltwo@$\ell^{(2)}$}$$\ell^{(2)}(\twocycG)=\sum_{(c,d) \in \CaC^2}g_{c d}
\ell(c,d)=\sum_{(c,d) \in \CaC^2}g_{c d}
\tilde{\ell}(c,d).$$
Furthermore, $\ell^{(2)}(\twocycG)$ is independent of the choices of the basepoints $\pointp(\alpha_i)$ or $\pointp(\beta_j)$, and of the numberings and orientations of the curves $\alpha_i$ or $\beta_j$.
\end{proposition}
\bp
Let $\ell^{\prime}$ be defined as 
$\ell$ except that the basepoint $\pointp_i=\pointp(\alpha_i)$ of $\alpha_i$ is changed to a basepoint $\pointq_i$.
When $c \in \alpha_i \setminus [\pointp_i,\pointq_i[_{\alpha}$, 
\begin{equation}\label{eq}
\ell^{\prime}(c,d)-\ell(c,d)=-\langle [\pointp_i,\pointq_i[_{\alpha}, [\pointp(\beta(d)),d\halfbra_{\beta}\rangle +\sum_{(r,j) \in \underline{g}^2}\CJ_{jr} \langle [\pointp_i,\pointq_i[_{\alpha},\beta_j \rangle \langle \alpha_r ,[\pointp(\beta(d)),d\halfbra_{\beta} \rangle
\end{equation}
When $c \in [\pointp_i,\pointq_i[_{\alpha}$, $[\pointq_i,c|_{\alpha} \setminus [\pointp_i,c|_{\alpha}= [\pointq_i,\pointp_i[_{\alpha}= \alpha_i \setminus [\pointp_i,\pointq_i[_{\alpha}$.
Since $$\langle \alpha_i, [\pointp(\beta(d)),d\halfbra_{\beta}\rangle=\sum_{(r,j) \in \underline{g}^2}\CJ_{jr} \langle \alpha_i,\beta_j \rangle \langle \alpha_r ,[\pointp(\beta(d)),d\halfbra_{\beta} \rangle,$$
$\ell^{\prime}(c,d)-\ell(c,d)$ is given by Formula~\ref{eq} that does not depend on $c \in \alpha_i$ in this case either.
Then $$\sum_{(c,d) \in \CaC^2}g_{c d}
(\ell^{\prime}(c,d)-\ell(c,d))=\sum_{(c,d) \in \alpha_i \times \CaC}g_{c d}
(\ell^{\prime}(c,d)-\ell(c,d)).$$
Since $$\partial (\gamma(c) \times \gamma(d)_{\parallel})= (b_{j(c)}-a_{i(c)}) \times \gamma(d)_{\parallel} -\gamma(c) \times (b_{j(d)\parallel}-a_{i(d)\parallel}),$$ for any $d \in \CaC$,
$\sum_{c \in \alpha_i}g_{c d}=0$.
Since the right-hand side of Formula~\ref{eq} does not depend on $c\in \alpha_i$, this shows that $\sum_{(c,d) \in \CaC^2}g_{c d}
\ell(c,d)$ does not depend on the basepoint choice on $\alpha_i$. Similarly, it does not depend on the choices of the basepoints on the $\beta_j$.

Similarly, $\sum_{(c,d) \in \CaC^2}g_{c d}
\ell(c,d)=\sum_{(c,d) \in \CaC^2}g_{c d}
\tilde{\ell}(c,d)$.

Using $\tilde{\ell}$, changing the orientation of $\alpha_i$ changes $\halfbra \pointp(\alpha(c)),c\halfbra_{\alpha}$ to $-\alpha_i + \halfbra \pointp(\alpha(c)),c\halfbra_{\alpha}$ for $c \in \alpha_i$, and does not change $\sum_{(c,d) \in \CaC^2}g_{c d}
\tilde{\ell}(c,d)$.
\eop

\begin{remark}
Let $[S]$ be the homology class of $\{x\} \times \partial B_x$ in $M^2\setminus \mbox{diagonal}$, where $B_x$ is a ball of $M$ and $x$ is a point inside $B_x$. Then $H_2(M^2\setminus \mbox{diagonal};\QQ)=\QQ[S]$, and it is proved in \cite[Proposition~3.4]{lesHC} that the class of a $2$--cycle $$\twocycG=\sum_{(c,d) \in \CaC^2}g_{c d} (\gamma(c) \times \gamma(d)_{\parallel})$$ in $H_2(M^2\setminus \mbox{diagonal};\QQ)$ is $\ell^{(2)}(\twocycG)[S]$.
Furthermore, for two disjoint one-cycles $K$ and $L$ of $M$, the class of $K \times L$ in $H_2(M^2\setminus \mbox{diagonal};\QQ)$ is $lk(K,L)[S]$ so that Proposition~\ref{propevallk} provides an alternative proof of \cite[Proposition~3.4]{lesHC} when $\twocycG$ is the product of two one-cycles.
This is the needed case to produce combinatorial expressions of linking numbers involved in the variations of $p_1$ that we are going to study later.
\end{remark}

\begin{proposition}
\label{proplambdahMorse}
Set \index{GD@$\twocycG(\CD)$}$$\twocycG(\CD)=\sum_{(c,d) \in \CaC^2}\CJ_{j(c)i(d)}\CJ_{j(d)i(c)}\sigma(c)\sigma(d) (\gamma(c) \times \gamma(d)_{\parallel}) -\sum_{c \in \CaC} \CJ_{j(c)i(c)}\sigma(c)(\gamma(c) \times \gamma(c)_{\parallel}).$$
Then $\twocycG(\CD)$ is a $2$--cycle of $M^2$. Let \index{ltwo@$\ell_2(\CD)$}$\ell_2(\CD)=\ell^{(2)}(\twocycG(\CD))$.
Then $\ell_2(\CD)$ is an invariant of the Heegaard diagram that does not depend on the orientations and numberings of the curves $\alpha_i$ and $\beta_j$. It not change when the roles of the $\alpha$-curves or the $\beta$-curves are permuted either.
\end{proposition}
\bp
Let us first prove that $\twocycG(\CD)$ is a $2$-cycle.
Note that, for any $j$,
$$\sum_{c\in \beta_j}\CJ_{j(d)i(c)}\sigma(c)=\sum_{i=1}^g\CJ_{j(d)i}\langle \alpha_i,\beta_j \rangle=\delta_{jj(d)}=\left\{\begin{array}{ll}1 & \mbox{if}\; j=j(d)\\
0 &\mbox{if}\; j\neq j(d)
\end{array}\right.$$
and, for any $i$, $\sum_{c\in \alpha_i}\CJ_{j(c)i(d)}\sigma(c)=\sum_{j=1}^g\CJ_{ji(d)}\langle \alpha_i,\beta_j \rangle=\delta_{ii(d)}$.
Therefore, for any $d \in \CaC$,
$$\partial\left(\sum_{c \in \CaC}\CJ_{j(c)i(d)}\CJ_{j(d)i(c)}\sigma(c)\gamma(c)\right)=\CJ_{j(d)i(d)}(b_{j(d)}-a_{i(d)})=\CJ_{j(d)i(d)}\partial \gamma(d)$$ 
so that
$$\partial \left(\sum_{(c,d) \in \CaC^2}\CJ_{j(c)i(d)}\CJ_{j(d)i(c)}\sigma(c)\sigma(d) (\gamma(c) \times \gamma(d)_{\parallel}) \right)=$$
$$\sum_{d \in \CaC}\sigma(d)\CJ_{j(d)i(d)} (\partial \gamma(d)) \times \gamma(d)_{\parallel} -\sum_{c \in \CaC}\CJ_{j(c)i(c)}\sigma(c) \gamma(c) \times \partial \gamma(c)_{\parallel}$$
and  $\partial \twocycG(\CD)=0$.

Permuting the roles of the $\alpha_i$ and the $\beta_j$ reverses the orientation of $\partial \handA$ and changes $\CJ$ to the transposed matrix. It does not change $\ell_2(\CD)$ because of the symmetry in the definition of $\ell^{(2)}$.
\eop

\section{The \texorpdfstring{$\infty$}{infty}-combings \texorpdfstring{$\ComX(\pointw,\matchingo)$}{X(\pointw,P)} and their \texorpdfstring{$p_1$}{p1}}
\setcounter{equation}{0}
\label{seccombproof}

\subsection{On the \texorpdfstring{$\infty$}{infty}-combing $\ComX(\pointw,\matchingo)$}
\label{subcombX}

In order to finish our description of $\ComX(\pointw,\matchingo)$ started in the introduction, we need to describe the vector field that replaces the gradient field $\ComX_{\fMorse_M}$ in regular neighborhoods $N(\gamma_i=\gamma(\matchm_i))$ of
the flow lines $\gamma_i$ associated with a matching $\matchingo$ of $\CD$.
Up to renumbering and reorienting the $\beta_j$, assume that $\matchm_i \in \alpha_i \cap \beta_i$ to simplify notation.

Choose a natural trivialization $(X_1,X_2,X_3)$ of $T\check{M}$ on a regular neighborhood $N(\gamma_i)$ of $\gamma_i$, such that:
\begin{itemize}
\item $\gamma_i$ is directed by $X_1$,
\item the other flow lines never have $X_1$ as an oriented tangent vector,
\item $(X_1,X_2)$ is tangent to the ascending manifold $\Asc_i$ of $a_i$ (except on the parts of $\Asc_i$ near $b_i$ that
come from other crossings of $\alpha_i \cap \beta_i$), and $(X_1,X_3)$ is tangent to the descending manifold $\Bdesc_i$ of $b_i$ (except on the parts of $\Bdesc_i$ near $a_i$ that
come from other crossings of $\alpha_i \cap \beta_i$).
\end{itemize}
This parallelization identifies the unit tangent bundle $UN(\gamma_i)$ of $N(\gamma_i)$ with $S^2 \times N(\gamma_i)$.

\medskip
\noindent There is a homotopy  $\homotop\colon [0,1]\times (N(\gamma_i) \setminus \gamma_i) \rightarrow S^2$, such that
\begin{itemize}
\item
$\homotop(0,.)$ is the unit vector with the same direction as the gradient vector of the underlying Morse function $\fMorse_M$,
\item
$\homotop(1,.)$ is the constant map to $(-X_1)$ and \item
$\homotop(t,y)$ goes from $\homotop(0,y)$ to $(-X_1)$ along the shortest geodesic arc $[\homotop(0,y),-X_1]$ of $S^2$ from $\homotop(0,y)$ to $(-X_1)$.
\end{itemize}
Let $2\eta$ be the distance between $\gamma_i$ and $\partial N(\gamma_i)$ and
let $\ComX(y) = \homotop(\mbox{max}(0,1-d(y,\gamma_i)/\eta),y)$ on $N(\gamma_i)\setminus \gamma_i$, and
$\ComX=-\ComX_1$ along $\gamma_i$.

\begin{figure}[h]
\begin{center}
\begin{tikzpicture}
\useasboundingbox (-1.4,-1.4) rectangle (1.4,1.4);
\draw (-1.3,0) -- (1.3,0) (0,-1.3) -- node[very near end, right]{$\linh_i(B)$} (0,1.3);
\draw [very thick,->] (0,-.5) -- (0,-.2);
\draw [very thick,->] (0,.5) -- (0,.2);
\draw [very thick,->] (-.2,0) -- (-.5,0);
\draw [very thick,->] (.2,0) -- (.5,0);
\draw [very thick,->] (.2,.4) -- (.4,.2);
\draw [very thick,->] (.2,-.4) -- (.4,-.2);
\draw [very thick,->] (-.2,.4) -- (-.4,.2);
\draw [very thick,->] (-.2,-.4) -- (-.4,-.2);
\draw [very thick,->] (2.2,.6) -- (2.2,1.1);
\draw [very thick,->] (2.2,.6) -- (2.7,.6);
\draw (1.3,.3) node{$\linh_i(A)$} (2.7,.6) node[right]{$X_2$} (2.2,1.1) node[right]{$X_3$};
\end{tikzpicture}
\caption{Projection of $\ComX$}
\label{figprojX}
\end{center}
\end{figure}
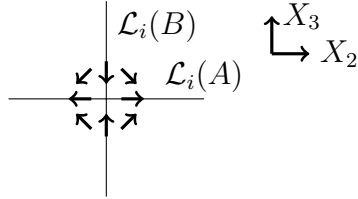

Note that $\ComX$ is tangent to $\Asc_i$ on $N(\gamma_i)$ (except on the parts of $\Asc_i$ near $b_i$ that
come from other crossings of $\alpha_i \cap \beta_i$), and that $\ComX$ is tangent to $\Bdesc_i$ on $N(\gamma_i)$ (except on the parts of $\Bdesc_i$ near $a_i$ that
come from other crossings of $\alpha_i \cap \beta_i$).
More generally, project the normal bundle of $\gamma_i$ to $\RR^2$ in the $\ComX_1$--direction by sending $\gamma_i$ to $0$, $\Asc_i$ to an axis $\linh_i(A)$ and $\Bdesc_i$ to an axis $\linh_i(B)$. Then the projection of $\ComX$ goes towards $0$ along $\linh_i(B)$ and starts from $0$ along $\linh_i(A)$, it has the direction of $s_a(y)$ at a point $y$ of $\RR^2$ near $0$, where
$s_a$ is the planar reflexion that fixes $\linh_i(A)$ and reverses $\linh_i(B)$. See Figure~\ref{figprojX}.

Then $\ComX(y)$ is on the half great circle that contains $s_a(y)$, $X_1$ and $(-X_1)$.
In Figure~\ref{figgammai}, $\gamma_i$ is a vertical segment, all the other flow lines corresponding to crossings involving $\alpha_i$ go upward from $a_i$, and $\ComX$ is simply the upward vertical field.

\begin{figure}[h]
\begin{center}
\begin{tikzpicture}
\useasboundingbox (1,1.4) rectangle (5,3.3);
\draw (1.3,1.5) .. controls (2,1.5) and (2,3.2) .. (3,3.2);
\draw [->] (4.7,1.5) .. controls (4,1.5) and (4,3.2) .. (3,3.2) ;
\draw [thick,->,fill=gray!20] (2.5,2) arc (180:0:.5);
\draw (3,2) circle (.5);
\draw (3.5,2) node[right]{\scriptsize $\beta_i$};
\draw [thick,->] (3,3.2) .. controls (2.9,3.2) and (2.8,3) .. (2.8,2.85);
\draw [thick] (3,2.5) .. controls (2.9,2.5) and (2.8,2.7) .. (2.8,2.85) node[left=-.2]{\scriptsize $\alpha_i$};
\draw [thick,dashed] (3,3.2) .. controls (3.1,3.2) and (3.2,3) .. (3.2,2.85) .. controls (3.2,2.7) and (3.1,2.5) .. (3,2.5);
\draw [->] (3,2.85) -- node[very near end,right=.1]{\scriptsize $\gamma_i$} (3,2);
\end{tikzpicture}
\caption{$\gamma_i$}
\label{figgammai}
\end{center}
\end{figure}

\subsection{On $p_1(\ComX(\pointw,\matchingo))$}
\label{subpone}

Invariants $p_1$ of $\infty$-combings --or of \emph{combings} that are homotopy classes of sections of the unit tangent bundle-- of rational homology spheres $M$ valued in $\QQ$
have been introduced and studied in \cite{lescomb} as extensions of a relative first Pontrjagin class from parallelizations to combings.

For a combing that extends to a parallelization $\tau$, the map $p_1$ coincides with the Hirzebruch defect (or Pontrjagin number) of the parallelization $\tau$, studied in \cite{hirzebruchEM, km, lesconst, lesmek}. 
For a parallelization $\tau\colon M \times \RR^3 \rightarrow TM$ of a $3$-manifold $M$ that bounds a connected oriented $4$-dimensional manifold $W$ with signature $0$, $p_1(\tau)$ is defined as the evaluation at the fundamental class of $[W,\partial W]$ of the relative first Pontrjagin class of $TW$ equipped with the trivialization of $TW_{|\partial W}$ that is the stabilization by the ``outward normal exterior first'' of $\tau$.
For $\infty$-combings that extend to parallelizations standard near $\infty$, $p_1$ is defined similarly by replacing $W$ by a connected oriented signature $0$ cobordism $W_c$ with corners between $B(1)$ and the rational homology ball $B_M$. 
A neighborhood of the boundary
$$\partial W_c=-B(1)\cup_{\partial B(1) \sim 0 \times B(1)}  (-[0,1] \times \partial B(1))  \cup_{\partial B_M \sim 1 \times \partial B(1)} B_M,$$  of such a cobordism
is naturally identified with an open subspace of one of the products $[0,1[ \times B(1)$ or $]0,1] \times B_M$ near $\partial W_c$, so that the standard parallelization of $\RR^3$ and $\tau$ induce a trivialization of $TW_{c|\partial W_c}$ by stabilizing by the ``tangent vector to $[0,1]$ first''.
For more details, see \cite[Section 1.5]{lesconst}.

Recall that any smooth compact oriented $3$-manifold $M$
can be equipped with a parallelization $\tau$.
When such a parallelization $\tau$ of $M$ is given, two sections $X$ and $Y$ of $U\check{M}$ or $U{M}$ induce a map $(X,Y)\colon \check{M} \rightarrow S^2 \times S^2$. Such sections are said to be {\em transverse\/} if the graphs of the induced maps $(X,Y)$ and $(X,-Y)$ are transverse to $\check{M} \times \mbox{diag}(S^2 \times S^2)$ in $\check{M} \times S^2 \times S^2$. This is generic and independent of $\tau$.
For two transverse sections $X$ and $Y$, let $L_{X=Y}$ be the preimage of the diagonal of $S^2$ under the map $(X,Y)$. Thus $L_{X=Y}$ is an oriented link that is cooriented by the fiber of the normal bundle to the diagonal of $(S^2)^2$. It can be assumed to avoid $\infty \in M$, generically.
In \cite[Theorem 1.2]{lescomb}, we proved that our extensions $p_1$ satisfy the following property that finishes defining them unambiguously.

\begin{theorem}
\label{thmlescomb}
When $X$ and $Y$ are two transverse $\infty$-combings (resp. combings) of a rational homology sphere $M$,
$$p_1(Y)-p_1(X)=4lk(L_{X=Y},L_{X=-Y}).$$
\end{theorem}

In \cite[Section 4.3]{lescomb}, we also proved that for combings, $p_1$ coincides with the invariant $\theta_G$ defined by Gompf in \cite[Section 4]{Go}.

The following properties of $p_1$ are easy to deduce from its definition.

\begin{proposition}
\begin{itemize}
\item A constant nonzero section $N$ of $T\RR^3$ is an $\infty$-combing of $S^3$ such that
$p_1(N)=0$.
\item Let $M$ and $M^{\prime}$ be rational homology spheres equipped with $\infty$-combings $X$ and $X^{\prime}$.
Assume that $X^{\prime}$ coincides with a constant section $N$ of $B(1)$ on $\partial B_{M^{\prime}}$ and that there is 
a standard ball $B(1)$ embedded in $M$ where $X$ coincides with $N$.
Replacing this embedded ball $(B(1),N)$ by $(B_{M^{\prime}},X^{\prime})$ gives rise to an $\infty$-combing $X^{\prime \prime}$ of the obtained manifold
whose $p_1$ is $p_1(X)+p_1(X^{\prime})$.
\item Changing the orientation of $M$ changes $p_1(X)$ to $-p_1(X)$.
\end{itemize}
\end{proposition}
\eop

Let $\xi$ be an oriented plane bundle over a compact oriented surface $S$ and let $\sigma$ be a nowhere vanishing section of $\xi$ on $\partial S$. The \emph{relative Euler number} $e(\xi,S,\sigma)$ of $\sigma$ is the algebraic intersection
of an extension of $\sigma$ to $S$ with the zero section of $\xi$. When $S$ is connected, it is the obstruction to extending $\sigma$ as a nowhere vanishing section of $\xi$.
The following proposition is a direct corollary of consequences of Theorem~\ref{thmlescomb} derived in \cite{lescomb}. 

\begin{proposition}
\label{propvarpontcpprel}
Let $\matchingo$ and $\matchingo^{\prime}$ be two matchings of $\CD$. Let $\Link(\matchingo^{\prime},\matchingo )= \Link(\CD,\matchingo^{\prime}) - \Link(\CD,\matchingo)$, and let $\Sigma(\Link(\matchingo^{\prime},\matchingo ))$ be a compact oriented surface bounded by $\Link(\matchingo^{\prime},\matchingo )$ in $M \setminus (S^3 \setminus B(1))$. 
Consider the four following fields $\ComY^{++}$, $\ComY^{+-}$, $(\ComY^{-+}=-\ComY^{+-})$ and $(\ComY^{--}=-\ComY^{++})$ in a neighborhood of the $\gamma(c)$.
$\ComY^{++}$ and $\ComY^{+-}$ are positive normals for $\Asc_i$ (that is oriented like $D(\alpha_i)$) on 
$\Asc_i \cap \fMorse_M^{-1}(]-\infty,3])$, and
$\ComY^{++}$ and $\ComY^{-+}$ are positive normals for $\Bdesc_j$ on $\Bdesc_j \cap \fMorse_M^{-1}([3,+\infty[)$. These four fields are orthogonal to $\ComX(\pointw,\matchingo)$ over $\Link(\matchingo^{\prime},\matchingo )$ and they define parallels $\Link(\matchingo^{\prime},\matchingo )_{\parallel Y^{\varepsilon, \eta}}$ of $\Link(\matchingo^{\prime},\matchingo )$ obtained by pushing in the $Y^{\varepsilon, \eta}$-direction.
Then
$$p_1(\ComX(\pointw,\matchingo^{\prime})) -p_1(\ComX(\pointw,\matchingo))=-\sum_{(\varepsilon,\eta) \in \{+,-\}^2}lk(\Link(\matchingo^{\prime},\matchingo ), \Link(\matchingo^{\prime},\matchingo )_{\parallel Y^{\varepsilon, \eta}}) + \Euhler(\pointw,\matchingo^{\prime},\matchingo )$$
where
$$\Euhler(\pointw,\matchingo^{\prime},\matchingo )= - \sum_{(\varepsilon,\eta) \in \{+,-\}^2} e(\ComX(\pointw,\matchingo)^{\perp},\Sigma(\Link(\matchingo^{\prime},\matchingo )),Y^{\varepsilon, \eta})$$
\end{proposition}
\bp Set $L=\Link(\matchingo^{\prime},\matchingo )$. Construct a cable $L_2$ of $L$ locally obtained by pushing one copy of $L$ in each direction normal to the $\Bdesc_j$, except near the $a_i$ where $L_2$ sits in $\Asc_i$. Define the field $Z$ over $L_2$ such that at a point $k$ of $L_2$, $Z$ has the direction of the vector from the closest point to $k$ on $L$ towards $k$.
Thus $\ComX(\pointw,\matchingo^{\prime})=D(\ComX(\pointw,\matchingo),L,L_2,Z,-1)$ with the notation Proposition~4.21 in \cite{lescomb}.

Then $((L_{\parallel Y^{+, +}}, L_{\parallel Y^{-, -}}), (Y^{+, +}, Y^{-, -}))$ is obtained from $(L_2,\ComZ)$ by some half-twists and $$((L_{\parallel Y^{+, -}}, L_{\parallel Y^{-, +}}), (Y^{+, -}, Y^{-, +}))$$ is obtained from $(L_2,\ComZ)$ by the opposite half-twists.
Then according to Proposition~4.21 in \cite{lescomb}, with the notation of \cite[Definition~4.16]{lescomb},
$$p_1(\ComX(\pointw,\matchingo^{\prime}))=\frac12 \left(p_1(D(\ComX(\pointw,\matchingo),L,L_{\parallel Y^{+, +}},Y^{+, +},-1)) + p_1(D(\ComX(\pointw,\matchingo),L,L_{\parallel Y^{+, -}},Y^{+, -},-1))\right).$$
Thus $p_1(\ComX(\pointw,\matchingo^{\prime}))= \frac14 \sum_{(\varepsilon,\eta) \in \{+,-\}^2}p_1(D(\ComX(\pointw,\matchingo),L,L_{\parallel Y^{\varepsilon,\eta}},Y^{\varepsilon,\eta},-1))$
and, according to \cite[Proposition~4.18 and Lemma~4.14]{lescomb},
$$\begin{array}{ll}p_1(\ComX(\pointw,\matchingo^{\prime}))- p_1(\ComX(\pointw,\matchingo))=&-\sum_{(\varepsilon,\eta) \in \{+,-\}^2}lk(\Link(\matchingo^{\prime},\matchingo ), \Link(\matchingo^{\prime},\matchingo )_{\parallel Y^{\varepsilon, \eta}})\\& - \sum_{(\varepsilon,\eta) \in \{+,-\}^2} e(\ComX(\pointw,\matchingo)^{\perp},\Sigma(\Link(\matchingo^{\prime},\matchingo )),Y^{\varepsilon, \eta}).\end{array}$$
\eop

The following theorem will be proved in Subsection~\ref{subpfbousillete}.
\begin{theorem}
\label{thmbousillete}
Let $L(\pointw^{\prime},\pointw)$ be the union of the closures of the flow line through $\pointw^{\prime}$ and the reversed flow line through $\pointw$.
 $$p_1(\ComX(\pointw^{\prime},\matchingo)) -p_1(\ComX(\pointw,\matchingo))=8 lk(\Link(\CD,\matchingo),L(\pointw^{\prime},\pointw))$$
\end{theorem}

\section{On the variations of $p_1(\ComX(\pointw,\matchingo))$}
\setcounter{equation}{0}
\label{secpXWCP}

\subsection{More on the variation of $p_1$ when $\matchingo$ changes}
\label{subpfpropvarpont}

\begin{lemma} \label{lemeultilth}
 Let $K=\sum_{c \in \CaC}k_{c} \gamma(c)$ be  a cycle of $M$, and let $\Sigma(K)$ be a surface bounded by $K$ in $\check{M}$. For $(\varepsilon,\eta) \in \{+,-\}^2$, let $Y^{\varepsilon, \eta}$ be the field
defined in Proposition~\ref{propvarpontcpprel} along the $\gamma(c)$.
Then $$\sum_{(\varepsilon,\eta) \in \{+,-\}^2} e(\ComX(\pointw,\matchingo)^{\perp},\Sigma(K),Y^{\varepsilon, \eta})=-4\sum_{c \in \CaC}k_{c}\tilth(c)$$
where $\tilth$ is defined before Lemma~\ref{leminvthetaW} with respect to our initial data that involve $(\pointw,\matchingo)$.
\end{lemma}
\bp Set $\ComX(\matchingo) =\ComX(\pointw,\matchingo)$. Since $M$ is a rational homology sphere, $e(\ComX(\matchingo)^{\perp},\Sigma(K),Y^{\varepsilon, \eta})$ does not depend on the surface $\Sigma(K)$. Choose the surface constructed in Lemma~\ref{lemseifsurfgam} with the points of $\matchingo$ as basepoints.
After removing the neighborhood $N(\gamma(\pointw))$ of the flow line through $\pointw$, $\fMorse_M^{-1}(]-\infty,0])$ behaves as a product by the rectangle $R_{\CD}$ of Figure~\ref{figplansplit} and has the product parallelization induced by the vertical vector field and the parallelization of $R_{\CD}$.
This parallelization extends to the one-handles of $\handA$ as the standard parallelization of $\RR^3$ in Figure~\ref{figgammai} so that it naturally extends to $\fMorse_M^{-1}(]-\infty,3])$, it furthermore extends to the neighborhood of the favourite flow lines in Figure~\ref{figgammai}. The first vector of this parallelization is $\ComX(\matchingo)$ and its second vector is everywhere orthogonal to $D(\alpha_i)$. It can be chosen to be $Y^{\varepsilon, \eta}$.
In a symmetric way, $\ComX(\matchingo)^{\perp}$ has a unit section that coincides with the second vector of the above parallelization on the neighborhoods of the favourite flow lines in Figure~\ref{figgammai} and that is orthogonal to $D(\beta_i)$ on $\fMorse_M^{-1}([4,\infty[) \setminus N(\gamma(\pointw))$.
Thus $e(\ComX(\matchingo)^{\perp},\Sigma(K),Y^{\varepsilon, \eta})$ reads
$$\sum_{c \in \CaC}k_{c}
e(\ComX(\matchingo)^{\perp},\halfbra \matchm_{j(c)},c \halfbra_{\beta}\times [3,4],\ComY^{\varepsilon,\eta})
-\sum_{(i,j,c) \in \underline{g}^2 \times \CaC}\CJ_{ji}k_{c}\langle \alpha_i,\halfbra \matchm_{j(c)},c\halfbra_{\beta}\rangle e(\ComX(\matchingo)^{\perp},\beta_j \times [3,4],\ComY^{\varepsilon,\eta})$$
where 
$\tilth(\halfbra \matchm_{j(c)},c \halfbra_{\beta})=-\frac14\sum_{(\varepsilon,\eta) \in \{+,-\}^2}e(\ComX(\matchingo)^{\perp},\halfbra \matchm_{j(c)},c \halfbra_{\beta}\times [3,4],\tilde{\ComY}^{\varepsilon,\eta})$ and
$$\tilth(\beta_{s})=-\frac14\sum_{(\varepsilon,\eta) \in \{+,-\}^2}e(\ComX(\matchingo)^{\perp},\beta_{s}\times [3,4],\tilde{\ComY}^{\varepsilon,\eta})$$ with respect to our partial extensions $\tilde{\ComY}^{\varepsilon,\eta}$ of ${\ComY}^{\varepsilon,\eta}$. (See \cite[Lemma~7.5]{lesHC} for more details.)
\eop

We get the following proposition as a direct corollary of Lemma~\ref{lemeultilth}.

\begin{proposition}
\label{propvarpontcp}
Under the hypotheses of Proposition~\ref{propvarpontcpprel}, if $\matchingo^{\prime}=\{d_j\}_{j\in \underline{g}}$, then
$$ \Euhler(\pointw,\matchingo^{\prime},\matchingo )=4 \sum_{j=1}^g \tilth(d_j)$$
where $\tilth$ is defined with respect to our initial data that involve $(\pointw,\matchingo)$.
\end{proposition}
\eop

Note that Lemma~\ref{leminvthetaWprel} independently implies that $\sum_{j=1}^g \tilth(d_j)$ only depends on $(\pointw,\matchingo,\matchingo^{\prime})$.

Lemma~\ref{lemeultilth} also yields the following second corollary that is \cite[Proposition~7.2]{lesHC}, that in turn yields Corollary~\ref{coreindepalphbet}.

\begin{corollary}
\label{coreultilth} Let $\Sigma(L(\CD,\matchingo))$ be a surface bounded by $L(\CD,\matchingo)$ in $\check{M}$.
$$e(\CD,\pointw,\matchingo)= \frac14\sum_{(\varepsilon,\eta) \in \{+,-\}^2} e(\ComX(\pointw,\matchingo)^{\perp},\Sigma(L(\CD,\matchingo)),Y^{\varepsilon, \eta})$$
\end{corollary}
\eop

\begin{corollary}
\label{coreindepalphbet}
$e(\CD,\pointw,\matchingo)$ is unchanged when the roles of the curves $\alpha$ and the curves $\beta$ are permuted.
\end{corollary}
\bp Permuting the roles of the curves $\alpha$ and the curves $\beta$ reverses the orientation of $L(\CD,\matchingo)$ and changes $\ComX(\pointw,\matchingo)$ to its opposite while the set $\{Y^{\varepsilon, \eta}\}_{(\varepsilon,\eta) \in \{+,-\}^2}$ is preserved.
\eop

\subsection{Associating a closed combing to a combing}
\label{subclscomb}

The Heegaard surface $\fMorse_M^{-1}(0)$ of our Morse function $\fMorse_M$ is obtained by gluing the complement $D_R$ of a rectangle in a sphere $S^2$ to the boundary of the rectangle $R_{\CD}$ of Figure~\ref{figplansplit}.
Let $D_R \times [-2,7]$ denote the intersection of $\fMorse_M^{-1}([-2,7])$ with the flow lines through $D_R$ so that $\fMorse_M$ is the projection to $[-2,7]$ on $D_R \times [-2,7]$ and the flow lines read $\{x\}\times [-2,7]$ there. Similarly, our Morse function $\fMorse_M$ reads as the projection on the interval on $$\fMorse_M^{-1}([-2,0] \cup [6,8] )=\left(S^2\times [-2,0]\right) \cup \left(S^2\times [6,8]\right)$$ while $\fMorse_M^{-1}([-3,-2])$ and $\fMorse_M^{-1}([7,9])$ are balls centered at a minimum and a maximum mapped to $-3$ and $9$, respectively.

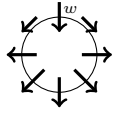
\begin{figure}[h]
\begin{center}
\begin{tikzpicture}
\useasboundingbox (-.8,-.8) rectangle (.8,.8);
\draw (0,0) circle (.5);
\draw [very thick,->] (0,-.3) -- (0,-.7);
\draw [very thick,->] (0,.7) -- (0,.3);
\draw [very thick,->] (.3,0) -- (.7,0);
\draw [very thick,->] (-.3,0) -- (-.7,0);
\draw [very thick,->] (-.2,-.2) -- (-.5,-.5);
\draw [very thick,->] (.2,-.2) -- (.5,-.5);
\draw [very thick,->] (.3,.5) -- (.5,.3);
\draw [very thick,->] (-.3,.5) -- (-.5,.3);
\draw (.15,.6) node{\tiny $\pointw$};
\end{tikzpicture}
\caption{The vector field near a minimum in a planar section of $\fMorse_M^{-1}([-3,-2])$}
\label{figXphimin}
\end{center}
\end{figure}

The combing $\ComX(\pointw,\matchingo)$ of Subsection~\ref{subcombX} of $B_M$ can be extended as a closed combing $\ComX(M,\pointw,\matchingo)$ that is obtained from the tangent $\ComX_{\phi}$ to the flow lines outside $B_M$ by reversing it along the line $\overline{\{\pointw\}\times ]-3,9[}$ as follows:

Let us first describe $\ComX(M,\pointw,\matchingo)$ on $D_R \times [-2,8]$. Let $D$ be a small disk of $D_R$ centered at $\pointw$. Reverse the flow on $\{\pointw\}\times [-2,8]$ so that it coincides with the tangent $\ComX_{\phi}$ to the flow outside $D \times [-2,8]$, and so that on a ray of $D\times \{t\}$ directed by a vector $\ComZ$ from the center, it describes the half great circle $[-\ComX_{\phi},\ComX_{\phi}]_{\ComZ}$ from $(-\ComX_{\phi})$ to $\ComX_{\phi}$ through $\ComZ$, if $t \in [-2,7]$.
Then on $S^2\times \{-2\}$, $\ComX$ is naturally homotopic to the restriction to the boundary of a constant field of $B^3$. See Figure~\ref{figXphimin} for a vertical section of the ball centered at the minimum where the constant vector field points downward. We extend it as such.

\begin{figure}[h]
\begin{center}
\begin{tikzpicture}
\useasboundingbox (-4.8,-1.5) rectangle (.8,1);
\draw (-4,0) circle (.7);
\draw [very thick,<-] (-4,-.5) -- (-4,-.9);
\draw [very thick,<-] (-4,.9) -- (-4,.5);
\draw [very thick,<-] (-3.5,0) -- (-3.1,0);
\draw [very thick,<-] (-4.5,0) -- (-4.9,0);
\draw [very thick,<-] (-4.3,-.45) -- (-4.75,-.75);
\draw [very thick,<-] (-3.7,-.45) -- (-3.25,-.75);
\draw [very thick,->] (-3.55,.75) -- (-3.25,.45);
\draw [very thick,->] (-4.45,.75) -- (-4.75,.45);
\draw (.2,.3) node{\tiny $\pointw$};
\draw (-4,-1.1) node{On $S^2 \times\{7\}$};
\draw (0,0) circle (.5);
\draw [very thick,<-] (0,-.3) -- (0,-.7);
\draw [very thick,<-] (0,.7) -- (0,.3);
\draw [very thick,<-] (.3,0) -- (.7,0);
\draw [very thick,<-] (-.3,0) -- (-.7,0);
\draw [very thick,<-] (-.2,-.2) -- (-.5,-.5);
\draw [very thick,<-] (.2,-.2) -- (.5,-.5);
\draw [very thick,<-] (.3,.5) -- (.5,.3);
\draw [very thick,<-] (-.3,.5) -- (-.5,.3);
\draw (.2,.3) node{\tiny $\pointw$};
\draw (0,-1.1) node{On $S^2 \times\{8\}$};
\end{tikzpicture}
\caption{The vector field near a maximum}
\label{figXphimax}
\end{center}
\end{figure}
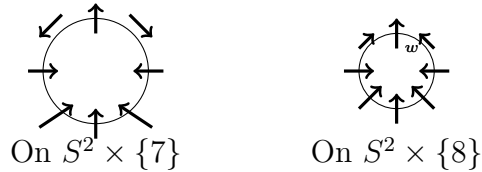

Now, on $S^2\times \{7\}$, $\ComX$ looks like in Figure~\ref{figXphimax}. It would naturally look as the restriction to the boundary of a constant field of $B^3$ if the half great circle $[-\ComX_{\phi},\ComX_{\phi}]_{\ComZ}$ from $(-\ComX_{\phi})$ to $\ComX_{\phi}$ through $\ComZ$ went through $(-\ComZ)$. 
Let $\rho_{\ComX_{\phi},\theta}$ denote the rotation with axis $\ComX_{\phi}$ and with angle $\theta$
For $t\in [7,8]$, on a ray of $D \times \{t\}$ directed by a vector $\ComZ$ from the center, let $\ComX$  describe the half great circle $[-\ComX_{\phi},\ComX_{\phi}]_{\rho_{\ComX_{\phi},(t-7)\pi}(\ComZ)}$ from $(-\ComX_{\phi})$ to $\ComX_{\phi}$ through $\rho_{\ComX_{\phi},(t-7)\pi}(\ComZ)$.
Now, we extend $\ComX$ as the constant field of $B^3$ that we see near the maximum, to obtain the closed combing $\ComX(M,\pointw,\matchingo)$.

\begin{lemma}
\label{lemindclos}
 $\left(p_1(\ComX(M,\pointw,\matchingo)) -p_1(\ComX(\pointw,\matchingo)\right))$ is a constant independent of $M$, $\pointw$ and $\matchingo$.
\end{lemma}
\bp Since the combing in the outside ball is unambiguously defined, and since it extends to a parallelization, this result follows from the definition of $p_1$ that includes the variation formula of Theorem~\ref{thmlescomb}. (See \cite{lescomb} for more details.) \eop 

\subsection{An abstract expression for the variation of $p_1$ when $\pointw$ varies}

This section is devoted to the proof of the following proposition that describes the variation of the Pontrjagin class $p_1(\ComX(\pointw,\matchingo))$ when $\pointw$ varies.
\begin{proposition}
\label{propchgW}
Let $\pointw$ and $\pointw^{\prime}$ be two exterior points of $\CD$.
Let $[\pointw,\pointw^{\prime}]_{\alpha}$ be a path on $\partial \handA$ from $\pointw$ to $\pointw^{\prime}$ disjoint from the $\alpha_i$ and let $[\pointw^{\prime},\pointw]_{\beta}$ be a path on $\partial \handA$ from $\pointw^{\prime}$ to $\pointw$ disjoint from the $\beta_j$. Set $[\pointw,\pointw^{\prime}]_{\beta}=-[\pointw^{\prime},\pointw]_{\beta}$.
Assume that the tangent vectors of $[\pointw,\pointw^{\prime}]_{\alpha}$ and $[\pointw,\pointw^{\prime}]_{\beta}$ at $\pointw$ and $\pointw^{\prime}$ coincide.
Let $$L([\pointw,\pointw^{\prime}]_{\alpha},[\pointw^{\prime},\pointw]_{\beta})=\left([\pointw,\pointw^{\prime}]_{\alpha}\times\{2\}\right) \cup \left(\{\pointw^{\prime}\}\times [2,4] \right) \cup \left([\pointw^{\prime},\pointw]_{\beta}\times\{4\}\right) \cup\left(\{\pointw\}\times [4,2] \right).$$
Let $\varepsilon=\pm 1$. Let $\ComY$ be a vector field defined on $L([\pointw,\pointw^{\prime}]_{\alpha},[\pointw^{\prime},\pointw]_{\beta})$ that is tangent to the Morse levels $\partial \handA \times \{t\}$ and that is a $\varepsilon$-normal (positive if $\varepsilon=1$ and negative otherwise) to $[\pointw,\pointw^{\prime}]_{\alpha}$ and a $(-\varepsilon)$ normal to $[\pointw^{\prime},\pointw]_{\beta}$. Let $L([\pointw,\pointw^{\prime}]_{\alpha},[\pointw^{\prime},\pointw]_{\beta})_{\parallel \ComY}$ be the induced parallel of $L([\pointw,\pointw^{\prime}]_{\alpha},[\pointw^{\prime},\pointw]_{\beta})$.
Let $\Sigma$ be a surface bounded by $L([\pointw,\pointw^{\prime}]_{\alpha},[\pointw^{\prime},\pointw]_{\beta})$.
Then $$\begin{array}{ll}p_1(\ComX(\pointw^{\prime},\matchingo))) -p_1(\ComX(\pointw,\matchingo)))=&4e(\ComX({\pointw},\matchingo)^{\perp},\Sigma,\ComY)\\ &
- 4lk(L([\pointw,\pointw^{\prime}]_{\alpha},[\pointw^{\prime},\pointw]_{\beta}),L([\pointw,\pointw^{\prime}]_{\alpha},[\pointw^{\prime},\pointw]_{\beta})_{\parallel \ComY}).\end{array}$$
\end{proposition}
\bp First note that $\ComX({\pointw},\matchingo)$ directs $\{\pointw^{\prime}\}\times [2,4]$ and $\{\pointw\}\times [4,2]$ so that the right-hand side of the equality above is independent of the field $\ComY$ that satisfies the conditions of the statement.
Let $L(\pointw^{\prime},\pointw)$ be the knot of $M$ that is the union of the closures of $\{\pointw^{\prime}\}\times ]-3,9[$ and $\{\pointw\}\times (-]-3,9[)$.
Let $\tilde{\ComX}(M,\pointw^{\prime},\matchingo)$ be obtained from $\ComX(M,\pointw,\matchingo)$ by reversing $\ComX(M,\pointw,\matchingo)$, that is tangent to $L(\pointw^{\prime},\pointw)$, along $L(\pointw^{\prime},\pointw)$.
In this situation, there is a standard way of reversing (namely the one that was used along $\{\pointw\} \times [-2,7]$ in Subsection~\ref{subclscomb}) by choosing a framing that determines both the parallel and the orthogonal field.

Proposition~\ref{propchgW} is the direct consequence 
of Lemma~\ref{lemindclos} and 
of the following three lemmas.
\begin{lemma}
\label{lemchgW1}
There exists a constant $C_0$ independent of $(M,\pointw,\pointw^{\prime},\matchingo)$ such that $$\begin{array}{ll}p_1(\tilde{\ComX}(M,\pointw^{\prime},\matchingo))-p_1(\ComX(M,\pointw,\matchingo))=& 4 e(\ComX({\pointw},\matchingo)^{\perp},\Sigma,\ComY)+ 4C_0 \\& -4 lk(L([\pointw,\pointw^{\prime}]_{\alpha},[\pointw^{\prime},\pointw]_{\beta}),L([\pointw,\pointw^{\prime}]_{\alpha},[\pointw^{\prime},\pointw]_{\beta})_{\parallel \ComY}).\end{array}$$
\end{lemma}

\begin{lemma}
\label{lemchgW2}
There exists a constant $C_1$ independent of $(M,\pointw,\pointw^{\prime},\matchingo)$ such that $$p_1(\tilde{\ComX}(M,\pointw^{\prime},\matchingo))-p_1(\ComX(M,\pointw^{\prime},\matchingo))=4 C_1.$$
\end{lemma}

\begin{lemma}
\label{lemchgW3}
 $$C_0-C_1=0.$$
\end{lemma}
\eop

\noindent{\sc Proof of Lemma~\ref{lemchgW1}:}
Let $T([\pointw,\pointw^{\prime}]_{\alpha})$ be the (closure of the) past of $[\pointw,\pointw^{\prime}]_{\alpha} \times\{2\}$ under the flow. This is a triangle and we can assume that it is smoothly embedded (near the minimum).
Similarly, let $T([\pointw^{\prime},\pointw]_{\beta})$ be the future of $[\pointw^{\prime},\pointw]_{\beta} \times\{4\}$ under the flow, assume without loss that it intersects $S^2\times \{7\}$ as a half-great circle, so that it intersects $\fMorse_M^{-1}([7,9])$ as a hemidisk denoted by $T_7([\pointw^{\prime},\pointw]_{\beta})$.  Orient $T([\pointw,\pointw^{\prime}]_{\alpha})$ and $T([\pointw^{\prime},\pointw]_{\beta})$ so that $$\partial (\Sigma + T([\pointw,\pointw^{\prime}]_{\alpha}) +T([\pointw^{\prime},\pointw]_{\beta}))= L(\pointw^{\prime},\pointw).$$

Then $\ComY$ extends to $T([\pointw,\pointw^{\prime}]_{\alpha})$ as the $(\varepsilon)$-normal on $T([\pointw,\pointw^{\prime}]_{\alpha})$ that is in $\ComX(M,\pointw,\matchingo)^{\perp}$. Similarly, $\ComY$ extends to $T([\pointw^{\prime},\pointw]_{\beta})$ as the $(\varepsilon)$-normal on $T([\pointw^{\prime},\pointw]_{\beta})$, it is a unit vector field that is in $\ComX(M,\pointw,\matchingo)^{\perp}$ outside the interior of $T_7([\pointw^{\prime},\pointw]_{\beta})$. Use $\ComY$ to frame $L(\pointw^{\prime},\pointw)$. 
Then, according to \cite[Proposition~4.18 and Lemma~4.14]{lescomb} where $\eta=1$, 
$$p_1(\tilde{\ComX}(M,\pointw^{\prime},\matchingo))-p_1(\ComX(M,\pointw,\matchingo))$$
$$= 4e(\ComX(M,\pointw,\matchingo)^{\perp},\Sigma +T_7([\pointw^{\prime},\pointw]_{\beta}),\ComY )-4lk(L(\pointw^{\prime},\pointw),L(\pointw^{\prime},\pointw)_{\parallel \ComY} )$$
where
$$ e(\ComX(M,\pointw,\matchingo)^{\perp},T_7([\pointw^{\prime},\pointw]_{\beta},\ComY )=C_0$$
for a constant $C_0$ independent of $(M,\pointw,\pointw^{\prime},\matchingo)$, and
$$lk(L(\pointw^{\prime},\pointw),L(\pointw^{\prime},\pointw)_{\parallel \ComY})=lk(L([\pointw,\pointw^{\prime}]_{\alpha},[\pointw^{\prime},\pointw]_{\beta}),L([\pointw,\pointw^{\prime}]_{\alpha},[\pointw^{\prime},\pointw]_{\beta})_{\parallel \ComY}).$$
\eop

\noindent{\sc Proof of Lemma~\ref{lemchgW2}:}
Recall that $D$ is a small disk of $\partial \handA$ centered at $\pointw$. The vector fields $\tilde{\ComX}(M,\pointw^{\prime},\matchingo)$ and $\ComX(M,\pointw^{\prime},\matchingo)$ coincide outside $\fMorse_M^{-1}([-3,-2]\cup[7,9]) \cup D\times [-2,7]$. This is a ball where the definition of these fields is unambiguous and independent of $(M,\pointw,\pointw^{\prime},\matchingo)$.
\eop

\noindent{\sc Proof of Lemma~\ref{lemchgW3}:}
According to the previous lemmas, for any
$(M,\pointw,\pointw^{\prime},\matchingo)$, 
$$\begin{array}{ll} p_1(\ComX(M,\pointw^{\prime},\matchingo))-p_1(\ComX(M,\pointw,\matchingo))= & -4lk(L([\pointw,\pointw^{\prime}]_{\alpha},[\pointw^{\prime},\pointw]_{\beta}),L([\pointw,\pointw^{\prime}]_{\alpha},[\pointw^{\prime},\pointw]_{\beta})_{\parallel \ComY})\\& + 4e(\ComX({\pointw},\matchingo)^{\perp},\Sigma,\ComY)+4(C_0-C_1).\end{array}$$
When $M$ is $S^3$ equipped with a Morse function with $2$ extrema and no other critical points, and when $\pointw$ and $\pointw^{\prime}$ are two points of $S^2$ related by a geodesic arc $[\pointw,\pointw^{\prime}]_{\alpha}=-[\pointw^{\prime},\pointw]_{\beta}$, it is easy to check that the terms of the formula are zero, except for $(C_0-C_1)$ that is therefore also $0$.
\eop

\subsection{A combinatorial formula for the variation of $p_1$ when $\pointw$ varies}
Now, we give an explicit formula for the right-hand side of Proposition~\ref{propchgW}.

\begin{figure}[h]
\begin{center}
\begin{tikzpicture}
\useasboundingbox (0,-1.4) rectangle (8.4,1.9);
\draw (0,-1.3) rectangle (8.4,1.8) (3,1.6) node{\scriptsize $\pointw$};
\draw [very thick,->]  (3.2,1.8) .. controls (3.2,1) and (4.5,-.8) .. (5.5,-.8);
\draw [very thick] (5.5,-.8) .. controls (10,-.8) and (5.5,4) .. (5.5,.2) (5.5,-1) node{\scriptsize $[\pointw,\pointw^{\prime}]_{\alpha}$} (4.5,1.2) node{\scriptsize $[\pointw,\pointw^{\prime}]_{\beta}$};
\draw [dashed,->]  (3.2,1.8) .. controls (3.2,1.3) and (4,1.5) .. (4.5,1.5);
\draw [dashed] (6.2,.7) .. controls (5.9,.7) and (5.5,.6) .. (5.5,.2) (5.5,0) node{\scriptsize $\pointw^{\prime}$};
\end{tikzpicture}
\caption{$[\pointw,\pointw^{\prime}]_{\alpha}$ and $[\pointw,\pointw^{\prime}]_{\beta}$}
\label{figWWprime}
\end{center}
\end{figure}
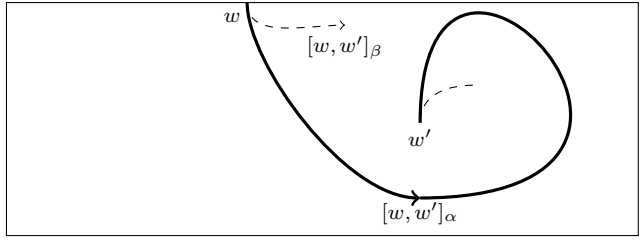

\begin{proposition}
\label{propformchgW}
Assume that $\pointw$ is on the upper side of the rectangle $R_{\CD}$ of Figure~\ref{figplansplit}.
Assume that $[\pointw,\pointw^{\prime}]_{\alpha}$ and $[\pointw,\pointw^{\prime}]_{\beta}=-[\pointw^{\prime},\pointw]_{\beta}$ point downward near $\pointw$ and $\pointw^{\prime}$ and that $[\pointw,\pointw^{\prime}]_{\beta}$ 
is on the same side of $[\pointw,\pointw^{\prime}]_{\alpha}$ near $\pointw$ and $\pointw^{\prime}$ as in Figure~\ref{figWWprime}.
Let $\tilth^{(\pointw)}([\pointw,\pointw^{\prime}]_{\alpha})$ be the degree of the tangent map to $[\pointw,\pointw^{\prime}]_{\alpha}$ on the rectangle $R_{\CD}$ of Figure~\ref{figplansplit}.
 Let $\tilth^{(\pointw)}([\pointw,\pointw^{\prime}]_{\beta})$ be the degree of the tangent map to $[\pointw,\pointw^{\prime}]_{\beta}$ on $R_{\CD}$, where $[\pointw,\pointw^{\prime}]_{\beta}$ intersects the $\alpha^{\prime}_j$ and the $\alpha^{\prime \prime}_j$ on their vertical portions opposite to the crossings of $\matchingo$, with horizontal tangencies. 
Then $$p_1(\ComX(\pointw^{\prime},\matchingo))) - p_1(\ComX(\pointw,\matchingo)))=p_1^{\prime}(\matchingo;\pointw,\pointw^{\prime})$$
where $$\begin{array}{ll} p_1^{\prime}(\matchingo;\pointw,\pointw^{\prime})=&
 4\tilth^{(\pointw)}([\pointw,\pointw^{\prime}]_{\alpha})-4\tilth^{(\pointw)}([\pointw,\pointw^{\prime}]_{\beta})\\&+4\sum_{(i,j)\in \underline{g}^2} \CJ_{ji}\langle \alpha_i,[\pointw,\pointw^{\prime}]_{\beta}\rangle \tilth^{(\pointw)}(\beta_j)\\&-4\langle ]\pointw,\pointw^{\prime}[_{\alpha}, ]\pointw,\pointw^{\prime}[_{\beta}\rangle
\\&+4\sum_{(i,j)\in \underline{g}^2} \CJ_{ji}\langle \alpha_i,[\pointw,\pointw^{\prime}]_{\beta}\rangle \langle [\pointw,\pointw^{\prime}]_{\alpha},\beta_j\rangle
\end{array}$$
\end{proposition}

\bp Define the field $\ComY$ of Proposition~\ref{propchgW} along 
$\{\pointw^{\prime}\}\times [2,4]$ and $\{\pointw\}\times [4,2]$, as the field pointing to the right in Figure~\ref{figWWprime} that is preserved by the flow along $\{\pointw^{\prime}\}\times [2,4]$ and $\{\pointw\}\times [4,2]$, so that it is always normal to $[\pointw,\pointw^{\prime}]_{\alpha} \times [2,4]$ or $[\pointw,\pointw^{\prime}]_{\beta} \times [2,4]$ along 
$\{\pointw^{\prime}\}\times [2,4]$ and $\{\pointw\}\times [4,2]$.
Let $L=L([\pointw,\pointw^{\prime}]_{\alpha},[\pointw^{\prime},\pointw]_{\beta})$ and let $L_{\parallel}=L_{\parallel\ComY}$. 
The proposition follows by applying Proposition~\ref{propchgW}, with the computations of Lemmas~\ref{lemlkpar} and \ref{lemobstchgW} below (replacing $[\pointw,\pointw^{\prime}]_{\beta}=-[\pointw^{\prime},\pointw]_{\beta}$). \eop

\begin{lemma}
\label{lemlkpar}
$$
lk\left(L,L_{\parallel }\right)=-\langle ]\pointw,\pointw^{\prime}[_{\alpha}, ]\pointw^{\prime},\pointw[_{\beta}\rangle+\sum_{(i,j)\in \underline{g}^2} \CJ_{ji}\langle \alpha_i,[\pointw^{\prime},\pointw]_{\beta}\rangle \langle [\pointw,\pointw^{\prime}]_{\alpha},\beta_j\rangle.$$
\end{lemma}

In order to prove Lemma~\ref{lemlkpar}, we shall use the following lemma.

\begin{lemma}
\label{lemsigthree}
There is a surface $\Sigma([\pointw,\pointw^{\prime}]_{\alpha}, [\pointw^{\prime},\pointw]_{\beta})$ in $\partial \handA \setminus \mathring{D}_R$ such that
$$\begin{array}{ll}\partial \Sigma([\pointw,\pointw^{\prime}]_{\alpha}, [\pointw^{\prime},\pointw]_{\beta})=&[\pointw,\pointw^{\prime}]_{\alpha} - \sum_{(i,j)\in \underline{g}^2} \CJ_{ji}\langle[\pointw,\pointw^{\prime}]_{\alpha},\beta_j\rangle \alpha_i
   \\&  + [\pointw^{\prime},\pointw]_{\beta} - \sum_{(i,j)\in \underline{g}^2} \CJ_{ji}\langle \alpha_i,[\pointw^{\prime},\pointw]_{\beta}\rangle \beta_j.
  \end{array}
$$
Let $\pointw^{\prime}_E$ be a point very close to $\pointw^{\prime}$ on its right-hand side. Then
$$\begin{array}{ll}\langle \Sigma([\pointw,\pointw^{\prime}]_{\alpha}, [\pointw^{\prime},\pointw]_{\beta}), \pointw^{\prime}_E\rangle_{\partial \handA} =&-\langle ]\pointw,\pointw^{\prime}[_{\alpha}, ]\pointw^{\prime},\pointw[_{\beta}\rangle
\\&+\sum_{(i,j)\in \underline{g}^2} \CJ_{ji}\langle \alpha_i,[\pointw^{\prime},\pointw]_{\beta}\rangle \langle [\pointw,\pointw^{\prime}]_{\alpha},\beta_j\rangle.\end{array}$$

\end{lemma}
\bp Since the prescribed boundary $\partial \Sigma([\pointw,\pointw^{\prime}]_{\alpha}, [\pointw^{\prime},\pointw]_{\beta})$ is a cycle that does not intersect the $\alpha_i$ and the $\beta_j$, algebraically, the surface $\Sigma([\pointw,\pointw^{\prime}]_{\alpha}, [\pointw^{\prime},\pointw]_{\beta})$ exists.
Let $\pointw_E$ be a point very close to $\pointw$ on its right-hand side. Along a path $]\pointw_E,\pointw^{\prime}_E[_{\alpha}$ parallel to $]\pointw,\pointw^{\prime}[_{\alpha}$, the intersection of a point with $\Sigma([\pointw,\pointw^{\prime}]_{\alpha}, [\pointw^{\prime},\pointw]_{\beta})$ starts with the value $0$ and varies when the path meets $\partial \Sigma([\pointw,\pointw^{\prime}]_{\alpha}, [\pointw^{\prime},\pointw]_{\beta})$ so that
$$\begin{array}{lll}\langle \Sigma([\pointw,\pointw^{\prime}]_{\alpha}, [\pointw^{\prime},\pointw]_{\beta}), \pointw^{\prime}_E\rangle_{\partial \handA} &=&-\langle]\pointw_E,\pointw^{\prime}_E[_{\alpha},\partial \Sigma([\pointw,\pointw^{\prime}]_{\alpha}, [\pointw^{\prime},\pointw]_{\beta})\rangle \\
&=&-\langle ]\pointw_E,\pointw^{\prime}_E[_{\alpha}, ]\pointw^{\prime},\pointw[_{\beta}\rangle\\&&+\sum_{(i,j)\in \underline{g}^2} \CJ_{ji}\langle \alpha_i,[\pointw^{\prime},\pointw]_{\beta}\rangle \langle ]\pointw_E,\pointw^{\prime}_E[_{\alpha},\beta_j\rangle.\end{array}$$
\eop

\noindent{\sc Proof of Lemma~\ref{lemlkpar}:}
$L$ bounds $$
\begin{array}{ll}\Sigma_0=&\Sigma([\pointw,\pointw^{\prime}]_{\alpha}, [\pointw^{\prime},\pointw]_{\beta}) + ([\pointw,\pointw^{\prime}]_{\alpha}\times [2,3])- ([\pointw^{\prime},\pointw]_{\beta} \times [3,4]) \\&+ \sum_{(i,j)\in \underline{g}^2} \CJ_{ji}\langle[\pointw,\pointw^{\prime}]_{\alpha},\beta_j\rangle D(\alpha_i)  +\sum_{(i,j)\in \underline{g}^2} \CJ_{ji}\langle \alpha_i,[\pointw^{\prime},\pointw]_{\beta}\rangle D(\beta_j).\end{array}$$

The link $L_{\parallel \ComY}=L([\pointw,\pointw^{\prime}]_{\alpha},[\pointw^{\prime},\pointw]_{\beta})_{\parallel \ComY}$ does not meet the $D(\alpha_i)$ and the $D(\beta_j)$.
Therefore its intersection with
$\Sigma_0$ is the intersection of $\pointw^{\prime}_E$ with $\Sigma([\pointw,\pointw^{\prime}]_{\alpha}, [\pointw^{\prime},\pointw]_{\beta})$ so that Lemma~\ref{lemsigthree} yields the conclusion.
\eop

\begin{lemma}
\label{lemobstchgW}
Let $\Sigma_1$ be a surface bounded by $L$ in $M$. Then
$$\begin{array}{ll}e(\ComX({\pointw},\matchingo)^{\perp},\Sigma_1,{\ComY} )=&\tilth^{(\pointw)}([\pointw,\pointw^{\prime}]_{\alpha})-\tilth^{(\pointw)}([\pointw,\pointw^{\prime}]_{\beta})\\&-\sum_{(i,j)\in \underline{g}^2} \CJ_{ji}\langle \alpha_i,[\pointw^{\prime},\pointw]_{\beta}\rangle \tilth^{(\pointw)}(\beta_j).\end{array}$$
\end{lemma}
\bp
Let $\Sigma_2=\Sigma([\pointw,\pointw^{\prime}]_{\alpha}, [\pointw^{\prime},\pointw]_{\beta}) \times \{2\}$ with the surface $\Sigma([\pointw,\pointw^{\prime}]_{\alpha}, [\pointw^{\prime},\pointw]_{\beta}) \subset \partial \handA$ of Lemma~\ref{lemsigthree}. The link $L$ bounds
$$\begin{array}{ll}\Sigma_1=&\Sigma_2 - [\pointw^{\prime},\pointw]_{\beta} \times [2,4] + \sum_{(i,j)\in \underline{g}^2} \CJ_{ji}\langle[\pointw,\pointw^{\prime}]_{\alpha},\beta_j\rangle D_{\leq 2}(\alpha_i)\\&  +\sum_{(i,j)\in \underline{g}^2} \CJ_{ji}\langle \alpha_i,[\pointw^{\prime},\pointw]_{\beta}\rangle D_{\geq 2}(\beta_j)
\end{array}$$
where $D_{\leq 2}(\alpha_i)=D(\alpha_i) \cap \fMorse_M^{-1}([-3,2])$ and $D_{\geq 2}(\beta_j)=D(\beta_j) \cap \fMorse_M^{-1}([2,9])$.
First extend $\ComY$ to the product by $[2,4]$ of a short vertical segment $[\pointw,\pointw^{(S)}]$ from $\pointw$ to some point $\pointw^{(S)}$ below $\pointw$, such that $\ComX({\pointw},\matchingo)$ directs $\pointw \times [4,2]$ and $\pointw^{(S)} \times [2,4]$, and $\ComX({\pointw},\matchingo)$ is tangent to $[\pointw,\pointw^{(S)}]\times [2,4]$.
Truncate the rectangle of Figure~\ref{figWWprime} so that $\pointw^{(S)}$ is on its boundary and $\pointw^{(S)}$ replaces $\pointw$ in the right-hand side of the equality of the statement without change.
Now, $\ComX({\pointw},\matchingo)$ is orthogonal to this rectangle, and the restriction to $\{\pointw^{(S)},\pointw^{\prime}\} \times [2,4]$ of the field $\ComY$ extends as the field $\ComY_E$ that points East or to the right in Figures~\ref{figplansplit}, \ref{figWWprime} and \ref{figgammai} so that is normal to the $D(\alpha_i)$. This defines the standard extension associated with Figure~\ref{figplansplit}
of $\ComY_E$ on $\fMorse_M^{-1}([0,2])\setminus (D_R \times [0,2])$ that is extended to $[\pointw^{\prime},\pointw]_{\beta} \times \{4\}$ so that it is normal to $[\pointw^{\prime},\pointw]_{\beta} \times [2,4]$ along $[\pointw^{\prime},\pointw]_{\beta} \times \{4\}$.
The field $\ComY_E$ can be extended independently to $D_{\geq 2}(\beta_j)$ and to $[\pointw^{\prime},\pointw]_{\beta} \times [2,4]$ as a field normal to these surfaces.
Then the Euler class of $\ComY_E$ with respect to $\Sigma_1$ can be computed by comparing this extension to the standard one above on 
$[\pointw,\pointw^{\prime}]_{\beta} \times \{2\} + \sum_{(i,j)\in \underline{g}^2} \CJ_{ji}\langle \alpha_i,[\pointw^{\prime},\pointw]_{\beta}\rangle (\beta_j\times \{2\})$,
$$e( \ComX({\pointw},\matchingo)^{\perp},\Sigma_1,\ComY_E )=-\tilth^{(\pointw)}([\pointw,\pointw^{\prime}]_{\beta})-\sum_{(i,j)\in \underline{g}^2} \CJ_{ji}\langle \alpha_i,[\pointw^{\prime},\pointw]_{\beta}\rangle \tilth^{(\pointw)}(\beta_j).$$
The fields $\ComY_E$ and $\ComY$ coincide on $L \setminus ([\pointw,\pointw^{\prime}]_{\alpha} \times \{2\})$ and
$$e(\ComX({\pointw},\matchingo)^{\perp},\Sigma_1,{\ComY})=e( \ComX({\pointw},\matchingo)^{\perp},\Sigma_1,\ComY_E )+\tilth^{(\pointw)}([\pointw,\pointw^{\prime}]_{\alpha}).$$
\eop

\subsection{Proof of Theorem~\ref{thmbousillete}}
\label{subpfbousillete}

Thanks to Proposition~\ref{propformchgW},
in order to prove Theorem~\ref{thmbousillete}, we are left with the proof that 
 $$ p_1^{\prime}(\matchingo;\pointw,\pointw^{\prime})=8 lk(\Link(\CD,\matchingo),L(\pointw^{\prime},\pointw))$$
where $p_1^{\prime}(\matchingo;\pointw,\pointw^{\prime})$ is defined in the statement of Proposition~\ref{propformchgW} and $L(\pointw^{\prime},\pointw)$ is the union of the closures of the flow line through $\pointw^{\prime}$ and the reversed flow line through $\pointw$.
In order to prove this, fix an exterior point $\pointw_0$ of $\CD$, and define $p_1^{\prime\prime}(\pointw)$, for any point exterior point $\pointw$ of $\CD$, as

$$p_1^{\prime\prime}(\pointw)= p_1^{\prime}(\matchingo;\pointw_0,\pointw)-8 lk(\Link(\CD,\matchingo),L(\pointw,\pointw_0))$$

\begin{lemma}
 $p_1^{\prime\prime}$ satisfies the following properties:
\begin{itemize}
 \item $p_1^{\prime\prime}(\pointw)$ only depends on the connected component of $\pointw$ in the complement of the $\alpha_i$ and the $\beta_j$ in the closed surface $\partial \handA$,
\item $p_1^{\prime\prime}(\pointw_0)=0$,
\item For any $4$ points $\pointw$, $S$, $E$, $N$ located around a crossing $d \notin \matchingo$,
as in Figure~\ref{figneardsimp}
$$ p_1^{\prime\prime}(N)+p_1^{\prime\prime}(S)=p_1^{\prime\prime}(\pointw)+p_1^{\prime\prime}(E).$$
\begin{figure}[h]
\begin{center}
\begin{tikzpicture}
\useasboundingbox (-4.2,-1) rectangle (4.2,1);
\draw [very thick,->] (0,1) -- (0,-1) node[right]{$\alpha_{i(d)}$};
\draw [thick] (-2,0) -- (2,0) node[below,right]{$\beta_{j(d)}$};
\fill (-.5,.5) circle (.05) (.5,.5) circle (.05) (.5,-.5) circle (.05) (-.5,-.5) circle (.05);
\draw (-.8,.6) node{\scriptsize $\pointw$} (.8,-.4) node{\scriptsize $E$} (.8,.6) node{\scriptsize $N$} (-.8,-.4) node{\scriptsize $S$};
\end{tikzpicture}
\caption{Near d}
\label{figneardsimp}
\end{center}
\end{figure}
\end{itemize}
\end{lemma}
\bp
The first two properties come from the definition.
Let us prove the third one. 

Set
$$D=\left(p_1^{\prime\prime}(N)+p_1^{\prime\prime}(S)-(p_1^{\prime\prime}(\pointw)+p_1^{\prime\prime}(E))\right)$$
Note that $D$ is independent of $\pointw_0$, thanks to Proposition~\ref{propformchgW}, and that it reads
$D=D_1-8D_2$
with
$$D_1=p_1^{\prime}(\matchingo;\pointw,N)+ p_1^{\prime}(\matchingo;\pointw,S)- p_1^{\prime}(\matchingo;\pointw,E)\;\;
\mbox{and}\;\; D_2=lk(L(N,\pointw)+L(S,E),\Link(\CD,\matchingo))$$
and
$$\begin{array}{ll} p_1^{\prime}(\matchingo;\pointw,\pointw^{\prime})=&
 4\tilth^{(\pointw)}([\pointw,\pointw^{\prime}]_{\alpha})-4\tilth^{(\pointw)}([\pointw,\pointw^{\prime}]_{\beta})\\&+4\sum_{(i,j)\in \underline{g}^2} \CJ_{ji}\langle \alpha_i,[\pointw,\pointw^{\prime}]_{\beta}\rangle \tilth^{(\pointw)}(\beta_j)
\\&-4\langle \Sigma([\pointw,\pointw^{\prime}]_{\alpha}, [\pointw^{\prime},\pointw]_{\beta}), \pointw^{\prime}_E\rangle_{\partial \handA}
\end{array}$$
according to Proposition~\ref{propformchgW} and Lemma~\ref{lemsigthree}.

We are going to prove that $$D_1=8D_2=-8\sigma(d)\CJ_{j(d)i(d)}.$$
Let us first compute $D_1$.
Its computation involves paths $[\pointw,\pointw^{\prime}]_{\alpha}$ and $[\pointw,\pointw^{\prime}]_{\beta}$ starting 
from $\pointw$ on the upper side of the rectangle $R_{\CD}$ of Figure~\ref{figplansplit}, before reaching a point $\pointw^{\prime}=N$, $S$ or $E$. We assume that all these paths begin by following
a first path $[\pointw,\tilde{\pointw}]$ that connects $\pointw$ to a point $\tilde{\pointw}$ near $d$ in the complement of the curves $\alpha_i$ and $\beta_j$ in $R_{\CD}$ and that this path $[\pointw,\tilde{\pointw}]$ has tangent vectors pointing downward at its ends.
The degree of the path $[\pointw,\tilde{\pointw}]$ does not matter since it is counted twice with opposite sign
in $(\tilth^{(\pointw)}([\pointw,\pointw^{\prime}]_{\alpha})-\tilth^{(\pointw)}([\pointw,\pointw^{\prime}]_{\beta}))$.
Thus we may change $\pointw$ to $\tilde{\pointw}$ in $p_1^{\prime}(\matchingo;\pointw,\pointw^{\prime})$ or equivalently assume that $\pointw$
arises near $d$ as in Figure~\ref{figneardsimp} split along $\alpha_{i(d)}$ and embedded in 
Figure~\ref{figplansplit} as soon as we translate our initial conventions for tangencies near the boundaries.
Now (keeping the first composition by $[\pointw,\tilde{\pointw}]$ in mind) we can draw our paths $[\tilde{\pointw},S]_{\alpha}$ and $[\tilde{\pointw},N]_{\beta}$ in Figure~\ref{figneard}
below where $\tilde{\pointw}$ is denoted by $\pointw$.
These paths together with the other drawn paths $[N,E]_{\alpha}$ and $[S,E]_{\beta}$ bound a ``square`` $C$ around $d$.
In Figure~\ref{figneard}, there are also dashed paths $[\pointw,N]_{\alpha}$ and $[\pointw,S]_{\beta}$ that may be complicated outside the pictured neighborhood of our square but that meet this neighborhood as in the figure.
We choose $[\pointw,E]_{\alpha}$ (resp. $[\pointw,E]_{\beta}$) to be the path composition of $[\pointw,N]_{\alpha}$ and $[N,E]_{\alpha}$ (resp. $[\pointw,S]_{\beta}$ and $[S,E]_{\beta}$).

\begin{figure}[h]
\begin{center}
\begin{tikzpicture}
\useasboundingbox (-6.2,-2) rectangle (6.2,2);
\draw [very thick,->] (0,1.7) -- (0,-2) node[right]{\scriptsize $\alpha_{i(d)}$};
\draw [thick] (-4,0) -- (3,0)  (-4,0) node[below,left]{\scriptsize $\beta_{j(d)}$};
\fill (-1.2,.9) circle (.05) (1.2,.9) circle (.05) (1.2,-1.7) circle (.05) (-1.2,-1.7) circle (.05) (1.65,-1.7) circle (.05);
\draw [->] (-1.2,.9) -- (-1.2,-.3);
\draw (-1.2,-1.7) -- (-1.2,-.3);
\draw [->] (1.2,.9) -- (1.2,.2);
\draw (1.2,-1.7) -- (1.2,.2);
\draw [dashed,->] (-1.2,.9) .. controls (-1.2,-.2) and (-.4,1.9) .. (0,1.9);
\draw [dashed] (0,1.9) .. controls (.4,1.9) and (1.2,1.9) .. (1.2,.9);
\draw [->] (-1.2,.9) .. controls (-1.2,.1) and (.6,1.2) .. (1.1,1.3);
\draw (1.1,1.3) .. controls (1.4,1.3) and (1.2,1.2) .. (1.2,.9);
\begin{scope}[yshift=-2.6cm]
\draw [->] (-1.2,.9) .. controls (-1.2,.3) and (.6,1.2) .. (1.1,1.3);
\draw (1.1,1.3) .. controls (1.4,1.3) and (1.2,1.2) .. (1.2,.9);
\end{scope}
\draw [dashed,->] (-1.2,.9) .. controls (-1.2,.4) and (.4,1.8) .. (1.7,1.8);
\draw [dashed] (1.7,1.8) .. controls (2.4,1.8) and (3.2,.4) .. (3.2,0) .. controls (3.2,-2) and (-1.2,.6) .. (-1.2,-1.7);
\draw (-1.4,1) node{\scriptsize $\pointw$} (1.15,1.3) node[right]{\scriptsize $[\pointw,N]_{\beta}$} (-1.15,-.3) node[left]{\scriptsize $[\pointw,S]_{\alpha}$} (1.15,.3) node[right]{\scriptsize $[N,E]_{\alpha}$} (-1.4,-1.7) node{\scriptsize $S$} (1.35,-1.7) node{\scriptsize $E$} (1.95,-1.7) node{\scriptsize $E_E$} (1,.8) node{\scriptsize $N$} (1.15,-1.2) node[left]{\scriptsize $[S,E]_{\beta}$} (-.1,1.9) node[left]{\scriptsize $[\pointw,N]_{\alpha}$} (1.8,1.9) node[right]{\scriptsize $[\pointw,S]_{\beta}$};
\end{tikzpicture}
\caption{Near d}
\label{figneard}
\end{center}
\end{figure}
With these choices, the contribution to $D_1$ of the parts
$$\tilth^{(\pointw)}([\pointw,\pointw^{\prime}]_{\alpha})-\tilth^{(\pointw)}([\pointw,\pointw^{\prime}]_{\beta})+\sum_{(i,j)\in \underline{g}^2} \CJ_{ji}\langle \alpha_i,[\pointw,\pointw^{\prime}]_{\beta}\rangle \tilth^{(\pointw)}(\beta_j)$$ cancel.
When $\pointw^{\prime}$ is $E$, $N$ or $S$, let $\Sigma(\pointw^{\prime})=\Sigma([\pointw,\pointw^{\prime}]_{\alpha}, [\pointw^{\prime},\pointw]_{\beta})$, with the notation of Lemma~\ref{lemsigthree}. 
Then
$$\begin{array}{ll}D_1&=4\langle \Sigma(E), E_E \rangle-4\langle \Sigma(N), N_E \rangle -4\langle \Sigma(S), S_E \rangle  \\&=-4\langle \Sigma(N)+\Sigma(S)-\Sigma(E), E_E \rangle -4 \langle [N_E,E_E], \partial \Sigma(N) \rangle-4 \langle [S_E,E_E],  \partial \Sigma(S) \rangle \end{array}$$
where $\partial(\Sigma(N)+\Sigma(S)-\Sigma(E))=[\pointw,S]_{\alpha}-[N,E]_{\alpha} +[N,\pointw]_{\beta} -[E,S]_{\beta}$ so that
$(\Sigma(N)+\Sigma(S)-\Sigma(E))$ is our square and $\langle \Sigma(N)+\Sigma(S)-\Sigma(E), E_E \rangle=0$,

$$\begin{array}{ll}\langle [N_E,E_E]_{\alpha}, \partial \Sigma(N) \rangle&=
\langle [N_E,E_E]_{\alpha}, [N,\pointw]_{\beta} -\sum_{(i,j)\in \underline{g}^2} \CJ_{ji} \langle \alpha_i,[N,\pointw]_{\beta}\rangle \beta_j \rangle
\\&=  \sigma(d)\CJ_{j(d)i(d)}\end{array}$$

$$\begin{array}{ll}\langle  [S_E,E_E]_{\beta},  \partial \Sigma(S) \rangle&=
\langle -[\pointw,S]_{\alpha} +\sum_{(i,j)\in \underline{g}^2} \CJ_{ji} \langle [\pointw,S]_{\alpha},\beta_j\rangle \alpha_i,[S_E,E_E]_{\beta} \rangle
\\&= \sigma(d)\CJ_{j(d)i(d)}\end{array}$$
Then $D_1=-8\sigma(d)\CJ_{j(d)i(d)}$.

In order to compute $D_2$, construct a Seifert surface for $L(N,\pointw)+L(S,E)$ made of 
\begin{itemize}
 \item two triangles parallel to the $D(\beta)$ with bottom boundaries $[\pointw,N]_{\beta}$
and $[E,S]_{\beta}$,
\item two triangles parallel to the $D(\alpha)$ with top edges $[S,\pointw]_{\alpha}$
and $[N,E]_{\alpha}$,
\item our square $C$ bounded by $([N,\pointw]_{\beta} \cup [\pointw,S]_{\alpha} \cup [S,E]_{\beta} \cup
[E,N]_{\alpha})$ that is a meridian of $\gamma(d)$.
\end{itemize}

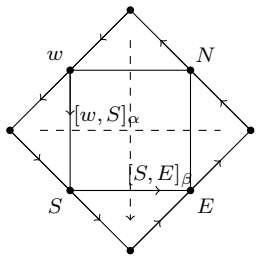
\begin{figure}[h]
\begin{center}
\begin{tikzpicture}
\useasboundingbox (-4.2,-2) rectangle (4.2,2);
\draw [dashed,->] (0,1.2) -- (0,-1.2);
\draw [dashed] (-1.2,0) -- (1.2,0);
\fill (-.8,.8) circle (.05) (.8,.8) circle (.05) (.8,-.8) circle (.05) (-.8,-.8) circle (.05);
\fill (1.6,0) circle (.05) (-1.6,0) circle (.05) (0,1.6) circle (.05) (0,-1.6) circle (.05);
\draw (-1,1) node{\scriptsize $\pointw$} (1,-1) node{\scriptsize $E$} (1,1) node{\scriptsize $N$} (-1,-1) node{\scriptsize $S$} (.4,-.6) node{\scriptsize $[S,E]_{\beta}$} (-.9,.2) node[right]{\scriptsize $[\pointw,S]_{\alpha}$};
\draw (1.6,0) -- (0,1.6) -- (-1.6,0) -- (0,-1.6) -- (1.6,0);
\draw (-.8,-.8) -- (-.8,.8) -- (.8,.8) -- (.8,-.8) -- (-.8,-.8);
\draw [->] (1.6,0) -- (1.2,.4);
\draw [->] (1.2,.4) -- (.4,1.2);
\draw [->] (0,1.6) -- (-.4,1.2);
\draw [->] (-.4,1.2) -- (-1.2,.4);
\draw [->] (-1.6,0) -- (-1.2,-.4);
\draw [->] (-1.2,-.4) -- (-.4,-1.2);
\draw [->] (0,-1.6) -- (.4,-1.2);
\draw [->] (.4,-1.2) -- (1.2,-.4);
\draw [->] (-.8,.8) -- (-.8,.2);
\draw [->] (-.8,-.8) -- (.4,-.8);
\end{tikzpicture}
\caption{A Seifert surface of $L(N,\pointw)+L(S,E)$}
\label{figseifsur}
\end{center}
\end{figure}

Therefore $D_2=-\sigma(d)\CJ_{j(d)i(d)}.$
\eop

Now, we conclude as follows. According to the above lemma, the variation of $p_1^{\prime\prime}$ across a curve $\alpha_i$ or $\beta_j$
is constant so that the variation of $p_1^{\prime\prime}$ along a path $\gamma$ reads 
$$\sum_i v_i \langle \gamma,\alpha_i \rangle + \sum_j w_j \langle \gamma,\beta_j \rangle$$
Since this is zero for any loop $\gamma$, the $v_i$ and the $w_j$ vanish, and the function $p_1^{\prime\prime}$ is constant. Then it is identically zero and Theorem~\ref{thmbousillete} is proved. \eop

\section{Behaviour of \texorpdfstring{$\tilde{\Theta}$}{Theta} when $\pointw$ and $\matchingo$ vary}
\setcounter{equation}{0}
\label{secvartilWCP}

In this section, we compute the variations of $\tilde{\Theta}(\pointw,\matchingo)$ when $\pointw$ and $\matchingo$ change, and we find that these variations
coincide with the 
variations of $\frac14p_1(\ComX(\pointw,\matchingo))$ computed in the previous section.
Thus we prove that $\left(\tilde{\Theta}(\pointw,\matchingo)-\frac14p_1(\ComX(\pointw,\matchingo))\right)$ is independent of $(\pointw,\matchingo)$.

\subsection{Changing $\pointw$}
\label{subsecvareW}

Let us first prove the following proposition that is similar to Proposition~\ref{propevallk}.

\begin{proposition}
\label{propevallinkbis}
Let $\pointw$ and $\pointw^{\prime}$ be two exterior points of $\CD$. Let $L(\pointw^{\prime},\pointw)$ be the union of the closures of the flow line through $\pointw^{\prime}$ and the reversed flow line through $\pointw$, let $[\pointw,\pointw^{\prime}]_{\beta}$ be a path from $\pointw$ to $\pointw^{\prime}$ outside the $\beta_j$ and let $[\pointw,\pointw^{\prime}]_{\alpha}$ be a path from $\pointw$ to $\pointw^{\prime}$ outside the $\alpha_i$. For any curve $\alpha_i$ (resp. $\beta_j$), choose a basepoint $\pointp(\alpha_i)$ (resp. $\pointp(\beta_j)$).
For any $1$--cycle $K=\sum_{c \in \CaC}k_{c} \gamma(c)$,
$$\begin{array}{lll}lk(K,L(\pointw^{\prime},\pointw)) 
&= &\sum_{c \in \CaC} k_c\langle [\pointw,\pointw^{\prime}]_{\alpha},[ \pointp(\beta(c)),c \halfbra_{\beta}\rangle 
\\&&-\sum_{(j,i) \in \underline{g}^2}\sum_{c \in \CaC}k_c\CJ_{ji}\langle \alpha_i, [ \pointp(\beta(c)),c \halfbra_{\beta} \rangle \langle [\pointw,\pointw^{\prime}]_{\alpha},\beta_j \rangle \end{array}
$$
where $\beta(c)=\beta_{j(c)}$.
\end{proposition}
\bp
As in Lemma~\ref{lemseifsurfgam}, $K$ bounds a chain 
$$\begin{array}{ll}\Sigma(K)=&\Sigma_{\Sigma}(K)+\sum_{c \in \CaC} k_c(T_{\beta}(c) + T_{\alpha}(c)) \\
&-\sum_{(j,i) \in \underline{g}^2}\sum_{c \in \CaC}k_c\CJ_{ji} \left(\langle \alpha_i, \halfbra \pointp(\beta(c)),c \halfbra_{\beta} \rangle D(\beta_j)
-\langle \halfbra \pointp(\alpha(c)),c \halfbra_{\alpha},\beta_j\rangle D(\alpha_i)\right)\end{array}$$
where $\Sigma_{\Sigma}(K)$ is a chain of $\partial \handA\setminus \{\pointw\}$ with boundary
$$\begin{array}{ll}\partial \Sigma_{\Sigma}(K)=&\sum_{c \in \CaC} k_c(\halfbra \pointp(\alpha(c)),c \halfbra_{\alpha}-\halfbra \pointp(\beta(c)),c \halfbra_{\beta}) \\
&+\sum_{(j,i) \in \underline{g}^2}\sum_{c \in \CaC}k_c\CJ_{ji} \left(\langle \alpha_i, \halfbra \pointp(\beta(c)),c \halfbra_{\beta} \rangle \beta_j
-\langle \halfbra \pointp(\alpha(c)),c \halfbra_{\alpha},\beta_j\rangle \alpha_i\right).\end{array}$$
Now, $lk(L,L(\pointw^{\prime},\pointw))$ is the intersection of $\pointw^{\prime}$ and $\Sigma_{\Sigma}(K)$, that is $ \langle -[\pointw,\pointw^{\prime}]_{\alpha}, \partial \Sigma_{\Sigma}(K)\rangle$.
\eop

\begin{figure}[h]
\begin{center}
\begin{tikzpicture}[scale=.6]
\useasboundingbox (-7.2,-2) rectangle (7.2,2);
\draw (-7,-1.9) rectangle (-1,1.9) (-4,1.9) -- (-4,.6) (-3.5,1.6) node{\scriptsize $\pointw$} (-3.7,.3) node{\scriptsize $\pointw^{\prime}$};
\draw (1,-1.9) rectangle (7,1.9) (4.3,-.75) node{\scriptsize $\pointw$} (4,1.6) node{\scriptsize $\pointw^{\prime}$};
\fill (-4,1.9) circle (.05) (4,1.9) circle (.05) (-4,.6) circle (.05) (4,1.9) circle (.05) (4,-1) circle (.05);
\draw  [->] (-5,.9) -- (-3,.9);
\draw  [->] (-3,1.2) -- (-5,1.2);
\draw  [->] (-.5,0) -- (.5,0);
\draw [->] (3,.9) .. controls (4,.9) and (3.8,1.8) .. (2.5,1.8) .. controls (1.2,1.8) and (1.1,1) .. (1.1,0) .. controls (1.1,-1) and (2,-1.8) .. (4,-1.8);
\draw [->] (4,-1.8) .. controls (6,-1.8) and (6.9,-1) .. (6.9,0) .. controls (6.9,1) and (6.8,1.8) .. (5.5,1.8) .. controls (4.2,1.8) and (4,.9) .. (5,.9);
\draw [<-] (3,1.2) .. controls (3.6,1.2) and (3.5,1.6) .. (2.5,1.6);
\draw (2.5,1.6) .. controls (1.3,1.6) and (1.3,.9) .. (1.3,0);
\draw [-<] (1.3,0) .. controls (1.3,-.9) and (2.2,-1.6) .. (4,-1.6);
\draw (4,-1.6) .. controls (5.8,-1.6) and (6.7,-.9) .. (6.7,0) .. controls (6.7,.9) and (6.7,1.6) .. (5.5,1.6) .. controls (4.5,1.6) and (4.4,1.2) .. (5,1.2);
\end{tikzpicture}
\caption{Changing $\pointw$ to $\pointw^{\prime}$}
\label{figchangW}
\end{center}
\end{figure}
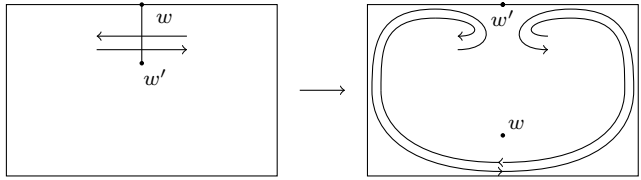

\begin{lemma}
\label{lemchangeW}
Let $\pointw$ and $\pointw^{\prime}$ be two exterior points of $\CD$. Let $L(\pointw^{\prime},\pointw)$ be the union of the closures of the flow line through $\pointw^{\prime}$ and the reversed flow line through $\pointw$.
Let $[\pointw,\pointw^{\prime}]_{\alpha}$ be a path of $\Sigma \setminus (\cup_{i=1}^g \alpha_i)$ from $\pointw$ to $\pointw^{\prime}$.
Set $$\tilde{\Theta}^{\prime}=
\tilde{\Theta}(\pointw^{\prime},\matchingo)-\tilde{\Theta}(\pointw,\matchingo)=e(\CD,\pointw,\matchingo)-e(\CD,\pointw^{\prime},\matchingo).$$
Then
$$\tilde{\Theta}^{\prime}=2\sum_{c \in \CaC} \CJ_{j(c)i(c)}\sigma(c)\left(\sum_{(r,s) \in \underline{g}^2}\CJ_{sr}\langle \alpha_r,\halfbra \matchm_{j(c)},c \halfbra_{\beta}\rangle \langle[\pointw,\pointw^{\prime}] _{\alpha},\beta_{s} \rangle - \langle[\pointw,\pointw^{\prime}]_{\alpha} ,\halfbra \matchm_{j(c)},c \halfbra_{\beta}  \rangle
\right).$$
\end{lemma}
\bp
Pick a vertical path $[\pointw,\pointw^{\prime}]_{\alpha}$  from a point $\pointw$ in the boundary of the rectangle of Figure~\ref{figplansplit} to the point $\pointw^{\prime}$ that cuts horizontal parts of the $\beta$ curves. When $\pointw$ is changed to $\pointw^{\prime}$, the portions or arcs near the intersection points
with $[\pointw,\pointw^{\prime}]_{\alpha}$ are transformed to arcs that turn around the whole picture of Figure~\ref{figplansplit}. This operation adds $2$ to the degree of an arc oriented from left to right. See Figure~\ref{figchangW}.

Therefore
$$\tilth^{(\pointw^{\prime})}(\beta_{s})-\tilth^{(\pointw)}(\beta_{s})=2\langle[\pointw,\pointw^{\prime}]_{\alpha} ,\beta_{s} \rangle$$
and $$\tilth^{(\pointw^{\prime})}(\halfbra \matchm_{j(c)},c\halfbra_{\beta})-\tilth^{(\pointw)}(\halfbra \matchm_{j(c)},c\halfbra_{\beta})=2\langle[\pointw,\pointw^{\prime}]_{\alpha} ,\halfbra \matchm_{j(c)},c\halfbra_{\beta}  \rangle.$$ 
\eop

\begin{corollary}
\label{corchangext}
$$\tilde{\Theta}(\pointw^{\prime},\matchingo)-\tilde{\Theta}(\pointw,\matchingo)=2lk(\Link(\CD,\matchingo),L(\pointw^{\prime},\pointw))=\frac14 p_1(\ComX(\pointw^{\prime},\matchingo))-\frac14p_1(\ComX(\pointw,\matchingo)).$$
\end{corollary}
This follows from Lemma~\ref{lemchangeW}, Proposition~\ref{propevallinkbis} and Theorem~\ref{thmbousillete}. \eop

\subsection{Changing $\matchingo$}
\label{secchangcro}

Let $\matchingo^{\prime}=\{d_{i} \in \alpha_i \cap \beta_{\rho^{-1}(i)}\}$ be another matching for a permutation $\rho$. The matching $\matchingo^{\prime}$ replaces our initial matching $\matchingo$ of positive crossings $\matchm_{i} \in \alpha_i \cap \beta_i$.

Set $\Link(\matchingo)=\Link(\CD,\matchingo)$ and $\Link(\matchingo^{\prime})=\Link(\CD,\matchingo^{\prime})$.

Let $\Link(\matchingo^{\prime},\matchingo )= \Link(\matchingo^{\prime}) - \Link(\matchingo) = \sum_{i=1}^g (\gamma(d_i) -\gamma_i)$.

This subsection is devoted to the proof the following proposition.
\begin{proposition}\label{propchangmatc}
 Under the assumptions above,
 $$\tilde{\Theta}(\pointw,\matchingo^{\prime})-\tilde{\Theta}(\pointw,\matchingo)=\frac14p_1(\ComX(\pointw,\matchingo^{\prime}))-\frac14p_1(\ComX(\pointw,\matchingo))$$
\end{proposition}

This proposition is a direct corollary of Propositions~\ref{propvarpontcpprel}, \ref{propvarpontcp} and \ref{propchangcro} so that we are left with the proof of Proposition~\ref{propchangcro} below.

\begin{proposition}
\label{propchangcro}
Under the assumptions above,
$$\begin{array}{ll}\tilde{\Theta}(\pointw,\matchingo^{\prime})-\tilde{\Theta}(\pointw,\matchingo)&=
lk(L(\matchingo^{\prime}), L(\matchingo^{\prime})_{\parallel})-lk(L(\matchingo), L(\matchingo)_{\parallel}) + e(\CD,\pointw,\matchingo) - e(\CD,\pointw,\matchingo^{\prime})\\
&=\sum_{i=1}^g\tilth(d_i) -lk(\Link(\matchingo^{\prime},\matchingo ),\Link(\matchingo^{\prime},\matchingo )_{\parallel}).\end{array}$$
Here, $\tilth$ is defined with respect to our initial data that involve $\pointw$ and $\matchingo$.
\end{proposition}

Proposition~\ref{propchangcro} is a direct consequence of Lemma~\ref{lemvarlkcro} below and Lemma~\ref{lemvarthetacro} that will be proved at the end of this subsection.

\begin{lemma}
\label{lemvarlkcro}
$$lk(\Link(\matchingo^{\prime}),\Link(\matchingo^{\prime})_{\parallel})-lk(\Link(\matchingo),\Link(\matchingo)_{\parallel})=2 lk(\Link(\matchingo^{\prime},\matchingo ),\Link(\matchingo^{\prime})_{\parallel})-lk(\Link(\matchingo^{\prime},\matchingo ),\Link(\matchingo^{\prime},\matchingo )_{\parallel})$$
\end{lemma}
\bp
Use
the symmetry of the linking number, and replace $\Link(\matchingo)=\Link(\matchingo^{\prime})-\Link(\matchingo^{\prime},\matchingo )$.
\eop

\begin{lemma}
\label{lemvarthetacro}
$$e(\CD,\pointw,\matchingo^{\prime})-e(\CD,\pointw,\matchingo)=2 lk(\Link(\matchingo^{\prime},\matchingo ),\Link(\matchingo^{\prime})_{\parallel})-\sum_{j=1}^g\tilth(d_j)
$$
where $\tilth(d_{\rho(j)})=\tilth(\halfbra \matchm_{j}, d_{\rho(j)}\halfbra_{\beta})- \sum_{s=1}^g\sum_{i=1}^g\CJ_{si}\langle \alpha_i,\halfbra \matchm_{j},d_{\rho(j)}\halfbra_{\beta}\rangle\tilth(\beta_{s})$.
\end{lemma}

\begin{lemma}
\label{lemlPP}
$$\begin{array}{ll}lk(\Link(\matchingo^{\prime},\matchingo ),\Link(\matchingo^{\prime})_{\parallel})=&\sum_{i=1}^g\sum_{c \in \CaC}\sigma(c) \CJ_{j(c)i(c)}\left(\sum_{(s,r) \in \underline{g}^2}\CJ_{sr}\langle\halfbra \matchm_i,d_i \halfbra_{\alpha}, \beta_s \rangle \langle \alpha_r ,\halfbra d_{\rho(j(c))},c\halfbra_{\beta} \rangle \right)\\
&-\sum_{i=1}^g\sum_{c \in \CaC}\sigma(c) \CJ_{j(c)i(c)}\left(\langle \halfbra \matchm_i,d_i \halfbra_{\alpha}, \halfbra d_{\rho(j(c))},c\halfbra_{\beta} \rangle \right).\\
\end{array}$$
\end{lemma}
\bp
Use Proposition~\ref{propevalbis} with $\pointp(\alpha_i)=\matchm_i$, $\pointp(\beta_j)=d_{\rho(j)}$, and $\tilde{\ell}$.
\eop

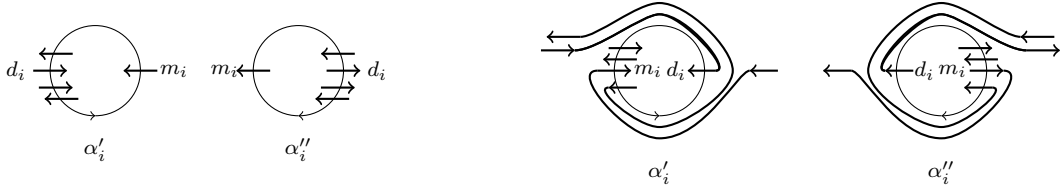
\begin{figure}[h]
\begin{center}
\begin{tikzpicture}[scale=1.5]
\useasboundingbox (-2,-.2) rectangle (6.4,1.9);
\begin{scope}[xshift=-3cm]
\draw [->] (1,.5) arc (-90:270:.4);
\draw (1,.2) node{\scriptsize $\alpha^{\prime}_i$};
\draw [->] (2.8,.5) arc (270:-90:.4);
\draw (2.8,.2) node{\scriptsize $\alpha^{\prime\prime}_i$};
\draw [thick,->] (2.55,.9) -- (2.25,.9);
\draw (1.7,.9) node{\scriptsize $\matchm_i$} (2.15,.9) node{\scriptsize $\matchm_i$} (3.5,.9) node{\scriptsize $d_i$} (.3,.9) node{\scriptsize $d_i$};
\draw [thick,->] (1.55,.9) -- (1.25,.9);
\draw [thick,->] (3.05,.9) -- (3.35,.9);
\draw [thick,->] (3,.75) -- (3.3,.75);
\draw [thick,->] (3.3,1.05) -- (3,1.05);
\draw [thick,->] (3.25,.65) -- (2.95,.65);
\draw [thick,->] (.85,.65) -- (.55,.65);
\draw [thick,->] (.45,.9) -- (.75,.9);
\draw [thick,->] (.5,.75) -- (.8,.75);
\draw [thick,->] (.8,1.05) -- (.5,1.05);
\end{scope}
\begin{scope}[xshift=2cm]
\draw [->] (1,.5) arc (-90:270:.4);
\draw (1,0) node{\scriptsize $\alpha^{\prime}_i$};
\draw (.9,.9) node{\scriptsize $\matchm_i$} (1.15,.9) node{\scriptsize $d_i$};
\draw [thick,->] (.8,1) -- (.5,1);
\draw [thick,->] (.55,1.1) -- (.85,1.1);
\draw [thick] (.25,1.08).. controls (.45,1.08) and (.85,1.4) .. (1,1.4) .. controls (1.2,1.4) and (1.65,.9) .. (1.5,.9);

\draw [thick,->] (1.5,.9) -- (1.25,.9);
\draw [thick,->] (-.05,1.08) -- (.25,1.08);
\draw [thick] (.25,1.08).. controls (.45,1.08) and (.85,1.4) .. (1,1.4) .. controls (1.25,1.4) and (1.65,.9) .. (1.5,.9);

\draw [thick,->] (2.05,.9) -- (1.8,.9);
\draw [thick,->] (.45,.9) -- (.75,.9);
\draw [thick] (1.8,.9).. controls (1.7,.9) and (1.4,.3) .. (1,.3) .. controls (.7,.3) and (.2,.9) .. (.45,.9);

\draw [thick,->] (.3,1.2) -- (0,1.2);
\draw [thick,->] (.8,.75) -- (.55,.75);
\draw [thick] (.3,1.2).. controls (.45,1.2) and (.8,1.5) .. (1,1.5) .. controls (1.25,1.5) and (1.65,1.1) .. (1.65,.9)
.. controls (1.65,.7) and (1.2,.4) .. (1,.4) .. controls (.8,.4) and (.4,.75) .. (.55,.75);
\end{scope}
\begin{scope}[xscale=-1,xshift=-6.5cm]
\draw [->] (1,.5) arc (-90:270:.4);
\draw (1,0) node{\scriptsize $\alpha^{\prime\prime}_i$};
\draw (.9,.9) node{\scriptsize $\matchm_i$} (1.15,.9) node{\scriptsize $d_i$};
\draw [thick,<-] (.8,1) -- (.5,1);
\draw [thick,<-] (.55,1.1) -- (.85,1.1);
\draw [thick] (.25,1.08).. controls (.45,1.08) and (.85,1.4) .. (1,1.4) .. controls (1.2,1.4) and (1.65,.9) .. (1.5,.9);

\draw [thick,<-] (1.5,.9) -- (1.25,.9);
\draw [thick,<-] (-.05,1.08) -- (.25,1.08);
\draw [thick] (.25,1.08).. controls (.45,1.08) and (.85,1.4) .. (1,1.4) .. controls (1.25,1.4) and (1.65,.9) .. (1.5,.9);

\draw [thick,<-] (2.05,.9) -- (1.8,.9);
\draw [thick,<-] (.45,.9) -- (.75,.9);
\draw [thick] (1.8,.9).. controls (1.7,.9) and (1.4,.3) .. (1,.3) .. controls (.7,.3) and (.2,.9) .. (.45,.9);

\draw [thick,<-] (.3,1.2) -- (0,1.2);
\draw [thick,<-] (.8,.75) -- (.55,.75);
\draw [thick] (.3,1.2).. controls (.45,1.2) and (.8,1.5) .. (1,1.5) .. controls (1.25,1.5) and (1.65,1.1) .. (1.65,.9)
.. controls (1.65,.7) and (1.2,.4) .. (1,.4) .. controls (.8,.4) and (.4,.75) .. (.55,.75);
\end{scope}
\end{tikzpicture}
\caption{Making the crossings move around}
\label{figmovecross}
\end{center}
\end{figure}

\noindent{\sc Proof of Lemma~\ref{lemvarthetacro}:}
Move the crossings of $[\matchm_i,d_{i}]$, counterclockwise along $\alpha^{\prime\prime}_i$ and clockwise along $\alpha^{\prime}_i$ as in Figure~\ref{figmovecross} so that $\matchm_i$ and $d_{i}$ make half a loop and the crossings of $]\matchm_i,d_{i}[$ make a (almost) full loop until they reach the standard position with respect to $d_{i}$.

As in the proof of Lemma~\ref{leminvthetaW}, on both sides of each crossing $c$ of $]\matchm_i,d_{i}[$ the degree is incremented by $(-\sigma(c))$, and it is incremented by $(-\sigma(\matchm_i)/2)$ on both sides of $\matchm_i$ and by $(-\sigma(d_i)/2)$ on both sides of $d_{i}$
so that after this modification the degree $\tilth^{\prime}(\beta_{j})$ of $\beta_j$
reads
$$\tilth^{\prime}(\beta_{j})=\tilth(\beta_{j})-2\sum_{i=1}^g\langle \halfbra \matchm_i,d_{i}\halfbra_{\alpha},\beta_j\rangle.$$

Before this modification, the degree of the tangent to $\beta_j$ from $d_{\rho(j(c))}$ to $c$ was $$\tilth(\halfbra d_{\rho(j(c))},c\halfbra_{\beta})= \left\{\begin{array}{ll}
\tilth(\halfbra d_{\rho(j(c))},\matchm_{j(c)} \halfbra_{\beta})) +\tilth(\halfbra \matchm_{j(c)},c \halfbra_{\beta}) \; & \mbox{if} \; c \in [\matchm_{j(c)},d_{\rho(j(c))}[_{\beta}\\
\tilth(\halfbra d_{\rho(j(c))},\matchm_{j(c)} \halfbra_{\beta})) +\tilth(\halfbra \matchm_{j(c)},c \halfbra_{\beta}) -\;\tilth(\beta_{j(c)}) & \mbox{if} \; c \in [d_{\rho(j(c))},\matchm_{j(c)}[_{\beta}
\end{array}\right.$$
After the modification, it reads 
$$\tilth^{\prime}(\halfbra d_{\rho(j(c))},c \halfbra_{\beta})= \tilth(\halfbra d_{\rho(j(c))},c\halfbra_{\beta})-2\sum_{i=1}^g\langle \halfbra \matchm_i,d_{i}\halfbra_{\alpha},\halfbra d_{\rho(j(c))},c\halfbra_{\beta}\rangle.$$
Now $$e(\CD,\pointw,\matchingo^{\prime})=\sum_{c \in \CaC} \CJ_{j(c)i(c)}\sigma(c)\tilth^{\prime}(c)$$
where $\tilth^{\prime}(c)=\tilth^{\prime}(\halfbra d_{\rho(j(c))},c \halfbra_{\beta}) -
\sum_{(r,s) \in \underline{g}2}\CJ_{sr}\langle \alpha_r,\halfbra  d_{\rho(j(c))},c \halfbra_{\beta}\rangle\tilth^{\prime}(\beta_{s}).$
Thus $$e(\CD,\pointw,\matchingo^{\prime})-e(\CD,\pointw,\matchingo)=e_1(\pointw,\matchingo,\matchingo^{\prime}) - e_2(\pointw,\matchingo,\matchingo^{\prime})$$
where
$$e_1(\pointw,\matchingo,\matchingo^{\prime})=\sum_{c \in \CaC} \CJ_{j(c)i(c)}\sigma(c)\left(\tilth^{\prime}(\halfbra d_{\rho(j(c))},c \halfbra_{\beta})-\tilth(\halfbra \matchm_{j(c)},c \halfbra_{\beta})\right)$$
and
$$e_2(\pointw,\matchingo,\matchingo^{\prime})=\sum_{c \in \CaC} \CJ_{j(c)i(c)}\sigma(c)\sum_{(r,s) \in \underline{g}^2}\CJ_{sr}\left(\langle \alpha_r,\halfbra  d_{\rho(j(c))},c \halfbra_{\beta}\rangle\tilth^{\prime}(\beta_{s})-\langle \alpha_r,\halfbra  \matchm_{j(c)},c \halfbra_{\beta}\rangle\tilth(\beta_{s})\right).$$

$$\begin{array}{lll}e_1(\pointw,\matchingo,\matchingo^{\prime})
&=&\sum_{c \in \CaC} \CJ_{j(c)i(c)}\sigma(c)\tilth(\halfbra d_{\rho(j(c))},\matchm_{j(c)} \halfbra_{\beta})-\sum_{c \in [d_{\rho(j(c))},\matchm_{j(c)}[_{\beta}} \CJ_{j(c)i(c)}\sigma(c)\tilth(\beta_{j(c)})\\
&&-2\sum_{c \in \CaC} \CJ_{j(c)i(c)}\sigma(c)\sum_{i=1}^g\langle \halfbra \matchm_i,d_{i}\halfbra_{\alpha},\halfbra d_{\rho(j(c))},c\halfbra_{\beta}\rangle\\
&=&\sum_{j=1}^g\tilth(\halfbra d_{\rho(j)},\matchm_{j} \halfbra_{\beta})- \sum_{j=1}^g\sum_{i=1}^g\CJ_{ji}\langle \alpha_i,[d_{\rho(j)},\matchm_{j}[_{\beta}\rangle\tilth(\beta_{j})\\
&&-2\sum_{c \in \CaC} \CJ_{j(c)i(c)}\sigma(c)\sum_{i=1}^g\langle \halfbra \matchm_i,d_{i}\halfbra_{\alpha},\halfbra d_{\rho(j(c))},c\halfbra_{\beta}\rangle. \end{array}$$
Since $$\begin{array}{lll}\langle \alpha_r,\halfbra  d_{\rho(j(c))},c \halfbra_{\beta}\rangle\tilth^{\prime}(\beta_{s})-\langle \alpha_r,\halfbra  \matchm_{j(c)},c \halfbra_{\beta}\rangle\tilth(\beta_{s})&=&
 -2\langle \alpha_r,\halfbra  d_{\rho(j(c))},c \halfbra_{\beta}\rangle \sum_{i=1}^g\langle \halfbra \matchm_i,d_{i}\halfbra_{\alpha},\beta_s \rangle \\&&+ \langle \alpha_r,\halfbra d_{\rho(j(c))}, \matchm_{j(c)}\halfbra_{\beta}\rangle\tilth(\beta_{s}) \\&&-\chi_{[d_{\rho(j(c))},\matchm_{j(c)}[_{\beta}}(c) \langle \alpha_r,\beta_{j(c)}\rangle\tilth(\beta_{s}),\end{array}
$$ 
where $\chi_{[d_{\rho(j(c))},\matchm_{j(c)}[_{\beta}}(c)=\left\{\begin{array}{ll}1\;&\mbox{if}\; c \in[d_{\rho(j(c))},\matchm_{j(c)}[_{\beta} \\
0 \; &\mbox{otherwise,} \end{array}\right.$

$$\begin{array}{lll}e_2(\pointw,\matchingo,\matchingo^{\prime})
&=& -2\sum_{c \in \CaC} \CJ_{j(c)i(c)}\sigma(c)\sum_{(r,s,i) \in \underline{g}^3}\CJ_{sr} \langle \halfbra \matchm_i,d_{i}\halfbra_{\alpha},\beta_s \rangle \langle \alpha_r,\halfbra  d_{\rho(j(c))},c \halfbra_{\beta}\rangle\\
&&+\sum_{(r,s,j) \in \underline{g}^3}\CJ_{sr}\langle \alpha_r,\halfbra d_{\rho(j)}, \matchm_{j}\halfbra_{\beta}\rangle\tilth(\beta_{s})\\
&&-\sum_{(r,s,j) \in \underline{g}^3}\CJ_{sr} \sum_{i=1}^g \CJ_{ji}\langle \alpha_i,[d_{\rho(j)},\matchm_{j}[_{\beta}\rangle \langle \alpha_r,\beta_{j}\rangle\tilth(\beta_{s})\\
&=& -2\sum_{c \in \CaC} \CJ_{j(c)i(c)}\sigma(c)\sum_{(r,s,i) \in \underline{g}^3}\CJ_{sr} \langle \halfbra \matchm_i,d_{i}\halfbra_{\alpha},\beta_s \rangle \langle \alpha_r,\halfbra  d_{\rho(j(c))},c \halfbra_{\beta}\rangle\\
&&+\sum_{(r,s,j) \in \underline{g}^3}\CJ_{sr}\langle \alpha_r,\halfbra d_{\rho(j)}, \matchm_{j}\halfbra_{\beta}\rangle\tilth(\beta_{s})\\
&&-\sum_{(i,j) \in \underline{g}^2} \CJ_{ji}\langle \alpha_i,[d_{\rho(j)},\matchm_{j}[_{\beta}\rangle \tilth(\beta_{j}).
\end{array}$$
Therefore, according to Lemma~\ref{lemlPP}, $e(\CD,\pointw,\matchingo^{\prime})-e(\CD,\pointw,\matchingo)=2lk(\Link(\matchingo^{\prime},\matchingo ),\Link(\matchingo^{\prime})_{\parallel})+V$
where $$\begin{array}{lll}V&=&\sum_{j=1}^g\tilth(\halfbra d_{\rho(j)},\matchm_{j} \halfbra_{\beta})- \sum_{(r,s,j) \in \underline{g}^3}\CJ_{sr}\langle \alpha_r,\halfbra d_{\rho(j)}, \matchm_{j}\halfbra_{\beta}\rangle\tilth(\beta_{s})\\
&=& -\sum_{j=1}^g\tilth(\halfbra \matchm_{j}, d_{\rho(j)}\halfbra_{\beta})+ \sum_{(r,s,j) \in \underline{g}^3}\CJ_{sr}\langle \alpha_r,\halfbra \matchm_{j},  d_{\rho(j)}\halfbra_{\beta}\rangle\tilth(\beta_{s})\\
&&+\sum_{j=1}^g\tilth(\beta_j)-\sum_{(r,s,j) \in \underline{g}^3}\CJ_{sr}\langle \alpha_r,\beta_j\rangle\tilth(\beta_{s})
        \end{array}
$$
Since the last line vanishes, we get the result.
\eop 

Corollary~\ref{corchangext} and Proposition~\ref{propchangmatc} allow us to define the function $\tilde{\lambda}$ of Heegaard diagrams
\index{lambdatilde@$\tilde{\lambda}(\CD)$}
$$\tilde{\lambda}(\CD)=\frac{\tilde{\Theta}(\CD,\pointw,\matchingo)}{6}-\frac{p_1(\ComX(\pointw,\matchingo))}{24}$$
that does not depend on the orientations and numberings of the curves $\alpha_i$ and $\beta_j$ and that is also unchanged by permuting the roles of the $\alpha_i$ and $\beta_j$, thanks to 
Corollary~\ref{coreindepalphbet}.

\section{Invariance of \texorpdfstring{$\tilde{\lambda}$}{lambda} }
\setcounter{equation}{0}
\label{secinvcomb}

In this section, we are first going to prove that $\tilde{\lambda}$ only depends on the Heegaard
decomposition induced by $\CD$ of $M$, and not on the curves $\alpha_i$ and $\beta_j$.
Then it will be easily observed that $\tilde{\lambda}$ is additive under connected sum of Heegaard decompositions and that $\tilde{\lambda}$ maps the genus one Heegaard decomposition of $S^3$ to $0$. Since according to the so-called Reidemeister-Singer theorem, two Heegaard decompositions of a $3$-manifold become diffeomorphic after some connected sums with this Heegaard decomposition of $S^3$, we will conclude that $\tilde{\lambda}$ is an invariant of rational homology spheres, that is additive under connected sum.

\subsection{Systems of meridians of a handlebody}

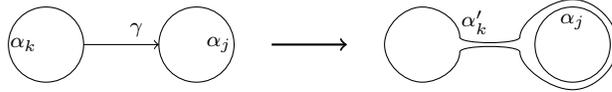
\begin{figure}[h]
\begin{center}
\begin{tikzpicture}
\useasboundingbox (-5.4,-.7) rectangle (5.4,.7);
\draw (-3.5,-.5) arc (-90:270:.5);
\draw (-3.8,0) node{\scriptsize $\alpha_k$};
\draw [->] (-3,0) -- (-2,0);
\draw (-1.5,-.5) arc (-90:270:.5);
\draw (-1.2,0) node{\scriptsize $\alpha_j$};
\draw (-2.3,.2) node{\scriptsize $\gamma$};
\draw [thick,->] (-.5,0) -- (.5,0);
\draw (1,0) .. controls (1,-.25) and (1.25,-.5) .. (1.5,-.5) .. controls (1.75,-.5) and (2,-.3) .. (2,-.1) .. controls (2,0) and (2.8,0) .. (2.8,-.1) .. controls (2.8,-.3) and (3.2,-.6) .. (3.5,-.6) .. controls (3.8,-.6) and (4.1,-.3) .. (4.1,0) .. controls (4.1,.3) and (3.8,.6) .. (3.5,.6) .. controls (3.2,.6) and (2.8,.3) .. (2.8,.1) .. controls (2.8,0) and (2,0) .. (2,.1) .. controls (2,.3) and (1.75,.5) .. (1.5,.5) .. controls (1.25,.5) and  (1,.25) .. (1,0);
\draw (3.5,-.5) arc (-90:270:.5);
\draw (3.5,.3) node{\scriptsize $\alpha_j$};
\draw (2.2,.3) node{\scriptsize $\alpha^{\prime}_k$};
\end{tikzpicture}
\caption{Handle slide in $\partial \handA$}
\label{fighandslid}
\end{center}
\end{figure}

A \emph{handle slide} in a system $\{\alpha_i\}_{i \in \underline{g}}$ of meridians of a curve $\alpha_k$ across a curve $\alpha_j$, with $j \neq k$, is defined as follows:
Choose a path $\gamma$ in $\partial \handA$ from a point $\gamma(0) \in \alpha_k$ to a point $\gamma(1) \in \alpha_j$ such that $\gamma(]0,1[)$ does not meet $\cup_{i \in \underline{g}}\alpha_i$ and change $\alpha_k$ to the band sum $\alpha^{\prime}_k$ of $\alpha_k$ and a parallel of $\alpha_j$ on the $\gamma$-side as in Figure~\ref{fighandslid}.

A \emph{right-handed Dehn twist} about a simple closed curve $K(S^1)$ of a surface $F$ is a homeomorphism of $F$ that fixes the exterior of a collar $K(S^1) \times [-\pi,\pi]$ of $K$ in $F$ pointwise, and that maps $(K(\exp(i\theta)),t)$ to $(K(\exp(i(\theta +t +\pi))),t)$.

In order to prove that $\tilde{\lambda}$ only depends on the Heegaard
decompositions and not on the chosen systems $\{\alpha_i\}_{i \in \underline{g}}$ and $\{\beta_j\}_{j \in \underline{g}}$ of meridians of $\handA$ and $\handB$ we shall use the following
standard theorem.

\begin{theorem}
Up to isotopy, renumbering of meridians, orientation reversals of meridians, two meridian systems of a handlebody are obtained from one another by a finite number of handle slides.
\end{theorem}
\bp Let $\{\alpha_i\}_{i \in \underline{g}}$ and $\{\alpha^{\prime}_i\}_{i \in \underline{g}}$
be two systems of meridians of $\handA$. There exists an orientation-preserving diffeomorphism of $\handA$ that maps the first system to the second one.

\begin{figure}[h]
\begin{center}
\begin{tikzpicture}
\useasboundingbox (-2,-1.2) rectangle (7,1.2);
\draw [thick] (0,0) circle (.55);
\draw (-.8,0) circle (.1);
\draw [thick] (0,0) circle (1.1);
\draw  [thick, ->] (0,.8) arc (-270:90:.8);
\draw (1.4,.1) node{\scriptsize $\alpha_2$};
\draw [thick,->] (.35,0) -- (1.3,0);
\draw (.05,.95) node{\scriptsize $a_2$};
\draw (-1.5,.3) node[left]{\scriptsize foot of $H_1$};
\draw [->] (-1.5,.3) -- (-.85,.05);
\draw [thick,->] (2,0) -- (3,0);
\draw [->]  (5,.8) arc (-270:90:.8);
\draw (.05,.95) node{\scriptsize $a_2$};
\draw (4.2,0) circle (.1);
\draw [thick] (5.35,0) .. controls (5.7,0) .. (5.7,0.07) .. controls (5.7,.35) and (5.35,.7) .. (5,.7);
\begin{scope}[xshift=5cm]
\draw [thick,->] (0,.7) .. controls (-.35,.7) and (-.65,.35) .. (-.65,0) .. controls (-.65,-.15) and (-.65,-.3) .. (-.75,-.3) .. controls (-.85,-.3) and (-1,-.2) .. (-1,0) .. controls (-1,.45) and (-.47,.95) .. (0,.95) .. controls  (.47,.95) and (.95,.47) .. (.95,.1) .. controls (.95,0) .. (1.3,0);
\draw (1.3,.1) node[right]{\scriptsize $\theta_{12}(\alpha_2$)};
\draw  (0,0) circle (.55);
\draw (0,0) circle (1.1);
\end{scope}
\end{tikzpicture}
\caption{Action of $\theta_{12}$ on $\alpha_2$}
\label{figthetasuz}
\end{center}
\end{figure}
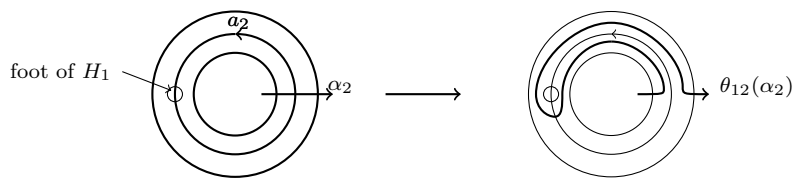

See $\handA$ as the unit ball $B(1)$ of $\RR^3$ with embedded handles $D(\alpha_i) \times [0,1]$ attached along $D(\alpha_i) \times \partial [0,1]$, so that there is a rotation $\rho$ of angle $\frac{2\pi}{g}$ of $\RR^3$ that maps $\handA$ to itself and that permutes the handles, cyclically.
See the meridians disks bounded by the $\alpha_i$ as disks $D(\alpha_i)=D(\alpha_i) \times \{\frac12 \}$ that cut the handles.
Let $H_i$ denote the handle of $\alpha_i$.
In \cite[Theorem 4.1]{Suzuki}, Suzuki proves that the group of isotopy classes of orientation-preserving diffeomorphisms of $\handA$ is generated by $6$ generators represented by the following diffeomorphisms
\begin{itemize}
 \item the rotation $\rho$ above of \cite[3.1]{Suzuki} that permutes the $\alpha_i$, cyclically, 
\end{itemize}
and the remaining $5$-diffeomorphisms that fix all the handles $H_i$, for $i>2$, pointwise,
\begin{itemize}
\item the knob interchange $\rho_{12}$ of \cite[3.4]{Suzuki}, that exchanges $H_1$ and $H_2$ and maps $\alpha_1$ to $\alpha_2$ and $\alpha_2$ to $\alpha_1$,
\item the knob twist $\omega_1$ of \cite[3.2]{Suzuki} that fixes $H_2$ pointwise and that maps $\alpha_1$ to the curve with opposite orientation, (it is the final time of an ambient isotopy of $\RR^3$ that performs a half-twist on a disk of $\handA$ that contains the two feet ($D(\alpha_1) \times \{0\}$ and $D(\alpha_1) \times \{1\}$) of the handle $H_1$),
\item the right-handed Dehn twist $\tau_1^{-1}$ of \cite[3.3]{Suzuki} along a curve parallel to $\alpha_1$,
\item the sliding $\xi_{12}$ of \cite[3.5 and 3.9]{Suzuki}, that is the final time of an ambient isotopy of $\RR^3 \times \RR$ that fixes the handles $H_i$, for $i>2$, pointwise, and that lets one foot of $H_1$ slide along a circle parallel to $\alpha_2$ once,
\item the sliding $\theta_{12}$ of \cite[3.5 and 3.8]{Suzuki}, that is the final time of an ambient isotopy of $\RR^3$ that fixes the handles $H_i$, for $i>2$, pointwise, and that lets one foot of $H_1$ slide along a circle $a_2$ that cuts $\alpha_2$ once and that does not meet the interiors of the $H_i$, for $i \neq 2$.
\end{itemize}
All these generators are described more precisely in \cite[Section 3]{Suzuki}.
All of them
except $\theta_{12}$ fix the set of curves $\alpha_i$ seen as unoriented curves, while $\theta_{12}$ fixes all the curves $\alpha_i$, for $i \neq 2$ pointwise.
When the foot of $H_1$ moves along the circle $a_2$, the curves that cross $a_2$ move with it, so that the meridian $\alpha_2$ is changed as in Figure~\ref{figthetasuz} that is a figure of a handle slide of $\alpha_2$ across $\alpha_1$.
\eop

\subsection{Isotopies of systems of meridians}

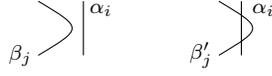
\begin{figure}[h]
\begin{center}
\begin{tikzpicture}[scale=.6]
\useasboundingbox (-2,-.2) rectangle (7,1.4);
\draw (0,0) .. controls (1,.6) .. (0,1.2)
(1,0) -- (1,1.2)
(.9,1) node[right]{\scriptsize $\alpha_i$} (.1,0) node[left]{\scriptsize $\beta_j$};
\draw (4,0) .. controls (5,.6) .. (4,1.2)
(4.5,0) -- (4.5,1.2)
(4.5,1) node[right]{\scriptsize $\alpha_i$} (4.1,0) node[left]{\scriptsize $\beta^{\prime}_j$};
\end{tikzpicture}
\caption{Birth of a bigon}
\label{figbigon}
\end{center}
\end{figure}

When the $\alpha_i$ are fixed on $\partial \handA$, and when the $\beta_j$ vary by isotopy, the only generic encountered accidents are the births or deaths of bigons that modify the Heegaard diagram as in Figure~\ref{figbigon} that represents
the birth of a bigon between an arc of $\alpha_i$ and an arc of $\beta_j$.

Therefore, in order to prove that $\tilde{\lambda}$ is invariant when the $\beta_j$ (or the $\alpha_i$) are moved by an isotopy, it
is enough to prove the following proposition:

\begin{proposition}
 \label{propbigon}
For any Heegaard diagram $\CD$ and $\CD^{\prime}$ such that
$\CD^{\prime}$ is obtained from $\CD$ by a birth of a bigon as above.
$$\tilde{\lambda}(\CD^{\prime})=\tilde{\lambda}(\CD).$$
\end{proposition}

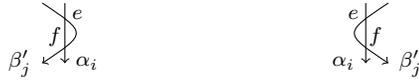
\begin{figure}[h]
\begin{center}
\begin{tikzpicture}
\useasboundingbox (-3,0) rectangle (3,1);
\draw [->] (-2,.8) -- (-2,0);
\draw [->] (-2.3,.8) ..  controls (-1.7,.4) .. (-2.3,0);
\draw [->] (2,.8) -- (2,0);
\draw [->] (2.3,.8) ..  controls (1.7,.4) .. (2.3,0);
\draw (-2,0) node[right]{\scriptsize $\alpha_i$} (-2.3,0) node[left]{\scriptsize $\beta^{\prime}_j$} (2,0) node[left]{\scriptsize $\alpha_i$} (2.3,0) node[right]{\scriptsize $\beta^{\prime}_j$} (-2.05,.65) node[right]{\scriptsize $e$} (-1.9,0.35) node[left]{\scriptsize $f$} (2.05,.65) node[left]{\scriptsize $e$} (1.9,0.35) node[right]{\scriptsize $f$};
\end{tikzpicture}
\caption{The considered two bigons}
\label{figtwobigons}
\end{center}
\end{figure}

Since we know that changing the orientation of $\alpha_i$ does not modify $\tilde{\lambda}$, we assume that our born bigon is one
of the two bigons shown in Figure~\ref{figtwobigons}, with two arcs going from a crossing $e$ to a crossing $f$, without loss.

We fix a matching $\matchingo$ for $\CD=((\alpha_i),(\beta_j))$ and the same one for $\CD^{\prime}$, and an exterior point $\pointw$ of $\CD^{\prime}$ outside the bigon so that $\pointw$ is also an exterior point of $\CD$.

\begin{lemma}
 $$p_1(\ComX(\CD,\pointw,\matchingo))=p_1(\ComX(\CD^{\prime},\pointw,\matchingo))$$
\end{lemma}
\bp The two fields $\ComX(\CD,\pointw,\matchingo)$ and $\ComX(\CD^{\prime},\pointw,\matchingo)$ may be assumed to coincide outside a ball that contains the past and the future in $\fMorse_M^{-1}([-2,7])$ of a disk of $\handA$ around the bigon, with respect to a flow associated with $\CD^{\prime}$. Since both fields are positive normals to the level surfaces
of $\fMorse_M$ on this ball they are homotopic.
\eop

Now, Proposition~\ref{propbigon} is a direct consequence of Lemmas~\ref{lemaddbigonCh} and \ref{lemaddbigonth}.

\begin{lemma}
\label{lemaddbigonCh}
$$\ell_2(\CD^{\prime})=\ell_2(\CD) +\CJ_{ji}/2.$$
$$s_{\ell}(\CD^{\prime},\matchingo)=s_{\ell}(\CD,\matchingo)$$
\end{lemma}
\bp
Let $\CaC$ be the set of crossings of $\CD$. Note that $\sigma(f)=-\sigma(e)$.
Then

$$\begin{array}{ll}\twocycG(\CD^{\prime})-\twocycG(\CD)=&\sum_{c \in \CaC}\CJ_{j(c)i}\CJ_{ji(c)}\sigma(c)\sigma(f) \gamma(c) \times (\gamma(f)-\gamma(e))_{\parallel} \\
&+\sum_{d \in \CaC}\CJ_{ji(d)}\CJ_{j(d)i}\sigma(d)\sigma(f) (\gamma(f)-\gamma(e)) \times \gamma(d)_{\parallel}\\
&+\CJ_{ji}^2 (\gamma(f)-\gamma(e)) \times (\gamma(f)-\gamma(e))_{\parallel}\\
&-\CJ_{ji}\sigma(f)(\gamma(f) \times \gamma(f)_{\parallel}-\gamma(e) \times \gamma(e)_{\parallel}).\end{array}$$

Use Proposition~\ref{propevalbis} to compute $\ell^{(2)}(\twocycG(\CD^{\prime})-\twocycG(\CD))$ with the basepoints of $\matchingo$,
so that for any $c \in \CaC$,
$$\ell(c,f)-\ell(c,e)= \langle [\pointp(\alpha(c)),c\halfbra_{\alpha},\halfbra e,f\halfbra_{\beta}\rangle -\sum_{(k,\ell) \in \underline{g}^2}\CJ_{\ell k} \langle [\pointp(\alpha(c)),c\halfbra_{\alpha},\beta_{\ell} \rangle \langle \alpha_k ,\halfbra e,f\halfbra_{\beta} \rangle=0$$
since $\langle \alpha_k ,\halfbra e,f\halfbra_{\beta} \rangle=0$ for any $k$, and $\langle [\pointp(\alpha(c)),c\halfbra_{\alpha},\halfbra e,f\halfbra_{\beta}\rangle=0$ for any $c  \in \CaC$.
Similarly, for any $d  \in \CaC$, $\ell(e,d)=\ell(f,d)$
and $$\ell(f-e,f-e)=\langle \halfbra e,f\halfbra_{\alpha},\halfbra e,f\halfbra_{\beta}\rangle=0.$$
Finally, $$\ell_2(\CD^{\prime})-\ell_2(\CD)=-\CJ_{ji}\sigma(f)(\ell(f,f)-\ell(e,e))$$
where
$$\ell(f,f)-\ell(e,e)= \langle [e,f\halfbra_{\alpha},[e,f\halfbra_{\beta}\rangle -  \langle [e,e\halfbra_{\alpha},[e,e\halfbra_{\beta}\rangle =\sigma(e) + \frac14 \sigma(f)-\frac14 \sigma(e) =-\frac12 \sigma(f)$$
so that $\ell_2(\CD^{\prime})-\ell_2(\CD)=\frac12\CJ_{ji}$.
Similarly, $s_{\ell}(\CD^{\prime},\matchingo)=s_{\ell}(\CD,\matchingo)$.
\eop

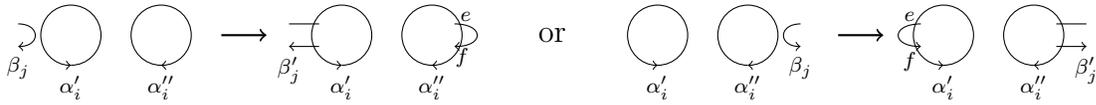
\begin{figure}[h]
\begin{center}
\begin{tikzpicture}
\useasboundingbox (-7,-.8) rectangle (7,.6);
\draw [->] (-6.4,-.4) arc (-90:270:.4);
\draw (-6.4,-.7) node{\scriptsize $\alpha^{\prime}_i$};
\draw [->] (-5.2,-.4) arc (270:-90:.4);
\draw (-5.2,-.7) node{\scriptsize $\alpha^{\prime\prime}_i$};
\draw [->] (-7.1,.15) ..  controls (-6.8,.15) and  (-6.8,-.15) .. (-7.1,-.15) node[below]{\scriptsize $\beta_j$};
\draw [thick,->] (-4.4,0) -- (-3.8,0);
\begin{scope}[xshift=3.6cm]
\draw [->] (-6.4,-.4) arc (-90:270:.4);
\draw (-6.4,-.7) node{\scriptsize $\alpha^{\prime}_i$};
\draw [->] (-5.2,-.4) arc (270:-90:.4);
\draw (-5.2,-.7) node{\scriptsize $\alpha^{\prime\prime}_i$};
\draw [->] (-7.1,.15) -- (-6.7,.15) (-6.7,-.15) -- (-7.1,-.15) node[below]{\scriptsize $\beta^{\prime}_j$};
\draw [->] (-4.9,.15) ..  controls (-4.5,.15) and  (-4.5,-.15) .. (-4.9,-.15);
\draw (-4.75,.25) node{\scriptsize $e$} (-4.8,-.3) node{\scriptsize $f$};
\end{scope}
\draw (0,0) node{or};
\begin{scope}[xshift=7.8cm]
\draw [->] (-6.4,-.4) arc (-90:270:.4);
\draw (-6.4,-.7) node{\scriptsize $\alpha^{\prime}_i$};
\draw [->] (-5.2,-.4) arc (270:-90:.4);
\draw (-5.2,-.7) node{\scriptsize $\alpha^{\prime\prime}_i$};
\draw [->] (-4.5,.15) ..  controls (-4.8,.15) and  (-4.8,-.15) .. (-4.5,-.15) node[below]{\scriptsize $\beta_j$};
\draw [thick,->] (-4,0) -- (-3.4,0);
\end{scope}
\begin{scope}[xshift=11.6cm]
\draw [->] (-6.4,-.4) arc (-90:270:.4);
\draw (-6.4,-.7) node{\scriptsize $\alpha^{\prime}_i$};
\draw [->] (-5.2,-.4) arc (270:-90:.4);
\draw (-5.2,-.7) node{\scriptsize $\alpha^{\prime\prime}_i$};
\draw [->] (-4.5,.15) -- (-4.9,.15) (-4.9,-.15) -- (-4.5,-.15) node[below]{\scriptsize $\beta^{\prime}_j$};
\draw [->] (-6.7,.15) ..  controls (-7.1,.15) and  (-7.1,-.15) .. (-6.7,-.15);
\draw (-6.85,.25) node{\scriptsize $e$} (-6.85,-.35) node{\scriptsize $f$};
\end{scope}
\end{tikzpicture}
\caption{Adding a bigon}
\label{figplansplitbig}
\end{center}
\end{figure}

\begin{lemma}
\label{lemaddbigonth}
$$e(\CD^{\prime},\pointw,\matchingo)=e(\CD,\pointw,\matchingo) + \CJ_{ji}/2.$$
\end{lemma}
\bp Adding a bigon changes Figure~\ref{figplansplit} as in Figure~\ref{figplansplitbig}.

In particular, the $\tilth(\beta_s)$ of Subsection~\ref{subdefe} are unchanged, and so are the $\tilth(c)$, for $c \in \CaC$.
Then $e(\CD^{\prime},\pointw,\matchingo)-e(\CD,\pointw,\matchingo)=\CJ_{ji}\sigma(f)\tilth(\halfbra e,f\halfbra_{\beta})$ that is $\frac12\CJ_{ji}$ according to Figure~\ref{figplansplitbig}.
\eop

\begin{remark}
If the two arcs of the bigon did not begin at the same vertex, then $\CJ_{ji}$ would be replaced by $-\CJ_{ji}$ in the results of 
Lemmas~\ref{lemaddbigonCh} and \ref{lemaddbigonth}.
\end{remark}

\subsection{Handle slides}

\begin{figure}[h]
\begin{center}
\begin{tikzpicture}
\useasboundingbox (-5.4,-.7) rectangle (5.4,.7);
\draw [blue,->] (-4,0) arc (-180:180:.5);
\draw [blue] (-3.8,0) node{\scriptsize $\beta_1$};
\draw [->] (-3,0) -- (-2,0);
\draw [blue,->] (-1,0) arc (0:360:.5);
\draw [blue] (-1.2,0) node{\scriptsize $\beta_2$};
\draw (-2.3,.2) node{\scriptsize $\gamma$};
\draw [red,thick] (-1.7,-.075) -- (-2.3, -.225);
\draw [red] (-2.2, -.15) node[below]{\scriptsize $\alpha_{i(e)}$};
\draw [thick,->] (-.5,0) -- (.5,0);
\draw [blue,->] (1,0) .. controls (1,-.25) and (1.25,-.5) .. (1.5,-.5) .. controls (1.75,-.5) and (2,-.3) .. (2,-.1) .. controls (2,0) and (2.8,0) .. (2.8,-.1) .. controls (2.8,-.3) and (3.2,-.6) .. (3.5,-.6) .. controls (3.8,-.6) and (4.1,-.3) .. (4.1,0) .. controls (4.1,.3) and (3.8,.6) .. (3.5,.6) .. controls (3.2,.6) and (2.8,.3) .. (2.8,.1) .. controls (2.8,0) and (2,0) .. (2,.1) .. controls (2,.3) and (1.75,.5) .. (1.5,.5) .. controls (1.25,.5) and  (1,.25) .. (1,0);
\draw [blue,->] (3.5,.5) arc (-270:90:.5);
\draw [red,thick] (3.3,-.075) -- (2.5, -.275);
\draw (3.15,-.2) node{\scriptsize $e$};
\draw (2.75,-.35) node{\scriptsize $e^+$};
\draw [blue] (3.5,.3) node{\scriptsize $\beta_2$};
\draw [blue] (1.2,0) node{\scriptsize $\beta^{\prime}_1$};
\end{tikzpicture}
\caption{The considered handle slide}
\label{fighandslidtwo}
\end{center}
\end{figure}
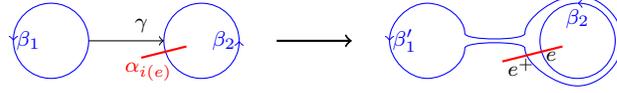

This section is devoted to proving that $\tilde{\lambda}$ is invariant under handle slide.
Since $\tilde{\lambda}$ depends neither on the orientations of the curves $\alpha_i$ and $\beta_j$, nor on their numberings, and since permuting the roles of the $\alpha_i$ and $\beta_j$ does not change $\tilde{\lambda}$, it is sufficient to study a handle slide that transforms $\CD$ to a diagram $\CD^{\prime}$ by changing $\beta_1$ to a band sum $\beta^{\prime}_1$ of $\beta_1$ and the parallel $\beta_2^+$ of $\beta_2$ (on its positive side) as in Figure~\ref{fighandslidtwo}. Up to the isotopies treated in the previous section, we may assume that the path $\gamma$ from $\beta_1$ to $\beta_2$ does not meet the curves $\alpha_i$, without loss, and we do.
The first crossing on $\beta_2^+$ will be called $e^+$. It corresponds to a crossing $e \in \alpha_{i(e)} \cap \beta_2$ as in Figure~\ref{fighandslidtwo}.

Fix $\pointw$ outside a neighborhood of the path $\gamma$ and $\beta_2$ so that it makes sense to say that $\pointw$ is the same for $\CD$ and $\CD^{\prime}$. Fix a matching $\matchingo$ for $\CD$. Assume $\matchingo = \{\matchm_i\}_{i \in \underline{g}}$ and $\matchm_i \in \alpha_i \cap \beta_i$ (by renumbering the $\alpha$ curves if necessary). The set $\CaC^{\prime}$ of crossings of $\CD^{\prime}$ contains $\CaC$ so that $\matchingo$ is also a matching for $\CD^{\prime}$.

Under these assumptions, we are going to prove that $\tilde{\lambda}(\CD^{\prime})=\tilde{\lambda}(\CD)$ by proving the following lemmas.

\begin{lemma}
\label{lemhandleslidepone} $$p_1(\ComX(\CD,\pointw,\matchingo))=p_1(\ComX(\CD^{\prime},\pointw,\matchingo)).$$\end{lemma}

\begin{lemma}
\label{lemhandleslideCh}
$$\begin{array}{ll}\ell_2(\CD^{\prime})-\ell_2(\CD)&=
\sum_{c \in \beta_2, d \in [e,c\halfbra_{\beta}} \sigma(c)\sigma(d)\CJ_{1i(c)}\CJ_{2i(d)}\\&\stackrel{\mbox{\scriptsize def}}{=}\sum_{c \in \beta_2, d \in [e,c[_{\beta}} \sigma(c)\sigma(d)\CJ_{1i(c)}\CJ_{2i(d)} + \frac12\sum_{c \in \beta_2}\CJ_{1i(c)}\CJ_{2i(c)} .\end{array}$$
\end{lemma}

\begin{lemma}
\label{lemhandleslidelk}
$$s_{\ell}(\CD^{\prime},\matchingo)-
s_{\ell}(\CD,\matchingo)= \sum_{d \in \beta_2, c \in [e,d\halfbra_{\beta}}\sigma(c)\sigma(d)\CJ_{1i(c)}\CJ_{2i(d)} -\sum_{c \in [e,\matchm_2\halfbra_{\beta}}\sigma(c)\CJ_{1i(c)}.$$
\end{lemma}

\begin{lemma}
\label{lemhandleslidetheta}
$$e(\CD^{\prime},\pointw,\matchingo)-e(\CD,\pointw,\matchingo)=\sum_{c \in \halfbra \matchm_2,e[_{\beta}}\sigma(c)\CJ_{1i(c)}.$$
\end{lemma}

Since 
$$\sum_{c \in \halfbra \matchm_2,e[_{\beta}}\sigma(c)\CJ_{1i(c)} + \sum_{c \in [e,\matchm_2\halfbra_{\beta}}\sigma(c)\CJ_{1i(c)}=\sum_{c \in \beta_2}\sigma(c)\CJ_{1i(c)}=\sum_{i=1}^g \CJ_{1i}\langle \alpha_i,\beta_2\rangle=0$$
and
$$\sum_{c \in \beta_2, d \in [e,c\halfbra_{\beta}} \sigma(c)\sigma(d)\CJ_{1i(c)}\CJ_{2i(d)}+\sum_{d \in \beta_2, c \in [e,d\halfbra_{\beta}}\sigma(c)\sigma(d)\CJ_{1i(c)}\CJ_{2i(d)}=\sum_{(c,d) \in \beta_2^2}\sigma(c)\sigma(d)\CJ_{1i(c)}\CJ_{2i(d)}=0,$$
these four lemmas imply that $\tilde{\lambda}(\CD^{\prime})=\tilde{\lambda}(\CD)$. \eop

\noindent{\sc Proof of Lemma~\ref{lemhandleslidepone}:}
Let $\ComX=\ComX(\CD,\pointw,\matchingo)$ and $\ComX^{\prime}=\ComX(\CD^{\prime},\pointw,\matchingo)$. First note that $\ComX$ and $\ComX^{\prime}$ coincide in $\handA$.
We describe a homotopy $(Y_t)_{t \in [0,1]}$ from  $Y_0=(-\ComX)$ and $Y_1=(-\ComX^{\prime})$ on $\handB$.

\begin{figure}[h]
\begin{center}
\begin{tikzpicture}
\useasboundingbox (-10,-2.7) rectangle (10,-.2);
\draw (2.5,-2.7) .. controls (2.05,-2.7) and (1,-.8) .. (0,-.8) .. controls (-1,-.8) and (-2.05,-2.7) .. (-2.5,-2.7);
\draw [thick,fill=gray!20] (-.5,-2.2) .. controls (-.8,-2.2) and (-.3,-1.25) .. (0,-1.25) .. controls (.3,-1.25) and (.8,-2.2) .. (.5,-2.2) arc (0:180:.5);
\draw (0,-2.2) circle (.5);
\draw [red,very thick,->] (.5,-2.2) arc (0:150:.5);
\draw [red,very thick] (.5,-2.2) arc (0:180:.5);
\draw [red] (.1,-2) node[left]{\scriptsize $\alpha_i$};
\draw [blue] (0,-.8) .. controls (-.1,-.8) and (-.2,-1) .. (-.2,-1.25);
\draw [blue,->] (0,-1.7) .. controls (-.1,-1.7) and (-.2,-1.35) .. (-.2,-1.25);
\draw [blue] (-.05,-1.15) node[left]{\scriptsize $\beta_i$};
\draw [blue,dashed] (0,-.8) .. controls (.1,-.8) and (.2,-1.05) .. (.2,-1.25) .. controls (.2,-1.45) and (.1,-1.7) .. (0,-1.7);
\draw [>->]  (0,-2.2) -- (0,-1.25);
\draw (.2,-2.2) node{\scriptsize $\gamma_i$};
\end{tikzpicture}
\caption{The front part of the disk $D_i$ where the field points inward the surface}
\label{figgammaidisk}
\end{center}
\end{figure}
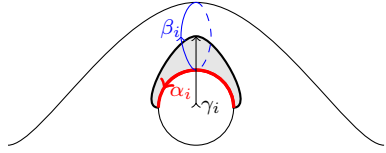

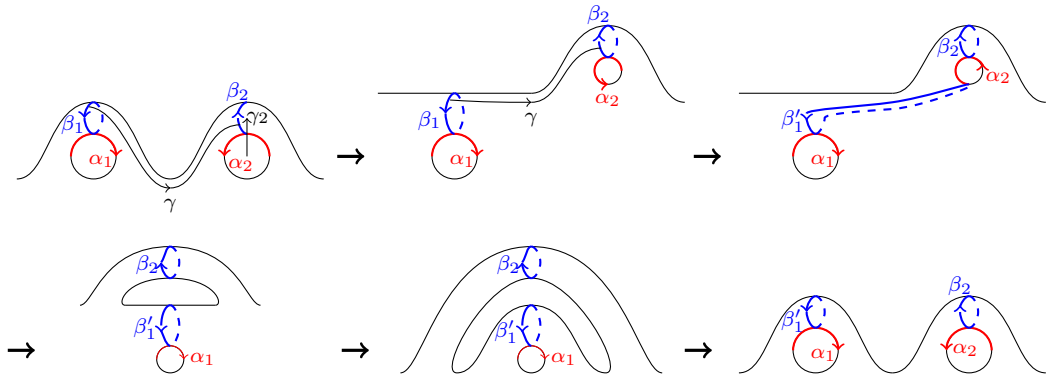
\begin{figure}[h]
\begin{center}
\begin{tikzpicture}[scale=.6]
\useasboundingbox (-10,-3) rectangle (10,5);
\begin{scope}[xshift=-9.3cm]
\draw (1.3,1.5) .. controls (2,1.5) and (2,3.2) .. (3,3.2);
\draw (4.7,1.5) .. controls (4,1.5) and (4,3.2) .. (3,3.2);
\draw [red] (2.35,1.8) node[right]{\scriptsize $\alpha_2$};
\draw (3,2) circle (.5);
\draw [red,thick,->] (3.5,2) arc (0:180:.5);
\draw [blue] (2.8,3.44) node{\scriptsize $\beta_2$};
\draw [blue,thick] (3,2.5) .. controls (2.9,2.5) and (2.8,2.7) .. (2.8,2.85);
\draw [blue,thick,-<] (3,3.2) .. controls (2.9,3.2) and (2.8,3) .. (2.8,2.85);
\draw (3.25,2.85) node{\scriptsize $\gamma_2$};
\draw [->] (3,2) -- (3,2.85);
\end{scope}
\begin{scope}[xshift=-12.7cm]
\draw (1.3,1.5) .. controls (2,1.5) and (2,3.2) .. (3,3.2);
\draw (4.7,1.5) .. controls (4,1.5) and (4,3.2) .. (3,3.2) ;
\draw (3,2) circle (.5);
\draw [red,thick,<-] (3.5,2) arc (0:180:.5);
\draw [red] (3.7,1.9) node[left]{\scriptsize $\alpha_1$};
\draw [blue,thick,->] (3,3.2) .. controls (2.9,3.2) and (2.8,3) .. (2.8,2.85);
\draw [blue,thick] (3,2.5) .. controls (2.9,2.5) and (2.8,2.7) .. (2.8,2.85);
\draw [blue] (3.05,2.7) node[left]{\scriptsize $\beta_1$};
\draw [blue,thick,dashed] (3,3.2) .. controls (3.1,3.2) and (3.2,3) .. (3.2,2.85) .. controls (3.2,2.7) and (3.1,2.5) .. (3,2.5);
\end{scope}
\begin{scope}[xshift=-9.3cm]
\draw [->] (-.5,3.1) .. controls (.3,3) and (.6,1.3) .. (1.3,1.3) node[below]{\scriptsize $\gamma$};
\draw (1.3,1.3) .. controls (2,1.3) and (2,2.7) .. (2.85,2.7);
\draw [very thick,->] (5,2) -- (5.6,2);
\end{scope}
\begin{scope}[xshift=-1.3cm,yshift=1.7cm]
\draw (1.3,1.7) .. controls (2,1.7) and (2,3.2) .. (3,3.2);
\draw (4.7,1.5) .. controls (4,1.5) and (4,3.2) .. (3,3.2);
\draw [red] (3,1.55) node{\scriptsize $\alpha_2$};
\draw (3,2.2) circle (.3);
\draw [red,thick,->] (3.3,2.2) arc (0:270:.3);
\draw [blue] (2.8,3.45) node{\scriptsize $\beta_2$};
\draw [blue,thick] (3,2.5) .. controls (2.9,2.5) and (2.8,2.7) .. (2.8,2.85);
\draw [blue,thick,-<] (3,3.2) .. controls (2.9,3.2) and (2.8,3) .. (2.8,2.85);
\draw [blue,thick,dashed] (3,3.2) .. controls (3.1,3.2) and (3.2,3) .. (3.2,2.85) .. controls (3.2,2.7) and (3.1,2.5) .. (3,2.5);
\end{scope}
\begin{scope}[xshift=-1.3cm]
\draw [->] (-.5,3.25) .. controls (.3,3.2) and (.6,3.2) .. (1.3,3.2) node[below]{\scriptsize $\gamma$};
\draw (1.3,3.2) .. controls (2,3.2) and (2,4.4) .. (2.85,4.4);
\end{scope}
\begin{scope}[xshift=-4.7cm]
\draw (1.3,3.4) -- (4.7,3.4);
\draw (3,2) circle (.5);
\draw [red,thick,<-] (3.5,2) arc (0:180:.5);
\draw [red] (3.7,1.9) node[left]{\scriptsize $\alpha_1$};
\draw [blue,thick,->] (3,3.4) .. controls (2.9,3.4) and (2.8,3.15) .. (2.8,2.95);
\draw [blue,thick] (3,2.5) .. controls (2.9,2.5) and (2.8,2.75) .. (2.8,2.95);
\draw [blue] (2.95,2.8) node[left]{\scriptsize $\beta_1$};
\draw [blue,thick,dashed] (3,3.4) .. controls (3.1,3.4) and (3.2,3.15) .. (3.2,2.95) .. controls (3.2,2.75) and (3.1,2.5) .. (3,2.5);
\end{scope}
\draw [very thick,->] (3.6,2) -- (4.2,2);
\begin{scope}[xshift=6.7cm,yshift=1.7cm]
\draw (1.3,1.7) .. controls (2,1.7) and (2,3.2) .. (3,3.2);
\draw (4.7,1.5) .. controls (4,1.5) and (4,3.2) .. (3,3.2);
\draw [red] (3.65,2.05) node{\scriptsize $\alpha_2$};
\draw (3,2.2) circle (.3);
\draw [red,thick,>-] (3.3,2.2) arc (0:270:.3);
\draw [blue] (3.05,2.7) node[left]{\scriptsize $\beta_2$};
\draw [blue,thick] (3,2.5) .. controls (2.9,2.5) and (2.8,2.7) .. (2.8,2.85);
\draw [blue,thick,-<] (3,3.2) .. controls (2.9,3.2) and (2.8,3) .. (2.8,2.85);
\draw [blue,thick,dashed] (3,3.2) .. controls (3.1,3.2) and (3.2,3) .. (3.2,2.85) .. controls (3.2,2.7) and (3.1,2.5) .. (3,2.5);
\end{scope}
\begin{scope}[xshift=3.3cm]
\draw (1.3,3.4) -- (4.7,3.4);
\draw (3,2) circle (.5);
\draw [red,thick,<-] (3.5,2) arc (0:180:.5);
\draw [red] (3.7,1.9) node[left]{\scriptsize $\alpha_1$};
\draw [blue,thick,-<] (3,2.5) .. controls (2.9,2.5) and (2.8,2.75) .. (2.8,2.95);
\draw [blue,thick] (2.8,2.95).. controls (2.8,3.05) and (4,3.1) .. (4.7,3.2) .. controls (5.4,3.3) and (6.3,3.6) .. (6.4,3.6);
\draw [blue] (3,2.8) node[left]{\scriptsize $\beta^{\prime}_1$};
\draw [blue,thick,dashed] (3,2.5) .. controls (3.1,2.5) and (3.2,2.75) .. (3.2,2.85) .. controls (3.2,2.95) and (4,3) .. (4.7,3.05) .. controls (5.4,3.1) and (6.6,3.6) .. (6.5,3.6);
\end{scope}
\begin{scope}[xshift=-8cm,yshift=-.5cm]
\draw (-1,-.8) -- (1,-.8);
\draw (1,-.8) .. controls (1.2,-.8) and (1,-.2) .. (0,-.2) .. controls (-1,-.2) and (-1.2,-.8) .. (-1,-.8);
\draw (-2,-.8) .. controls (-1.7,-.8) and (-1.5,.5) .. (0,.5) .. controls (1.5,.5) and (1.7,-.8) .. (2,-.8);
\draw (0,-2) circle (.3);
\draw [red,<-] (.3,-2) arc (0:180:.3);
\draw [red] (.2,-2) node[right]{\scriptsize $\alpha_1$};
\draw [blue,thick] (0,-.8) .. controls (-.1,-.8) and (-.2,-1) .. (-.2,-1.25);
\draw [blue,thick,-<] (0,-1.7) .. controls (-.1,-1.7) and (-.2,-1.35) .. (-.2,-1.25);
\draw [blue] (-.05,-1.3) node[left]{\scriptsize $\beta^{\prime}_1$};
\draw [blue,thick,dashed] (0,-.8) .. controls (.1,-.8) and (.2,-1.05) .. (.2,-1.25) .. controls (.2,-1.45) and (.1,-1.7) .. (0,-1.7);
\draw [blue,thick] (0,.5) .. controls (-.1,.5) and (-.1,.3) .. (-.2,.15);
\draw [blue,thick,->] (0,-.2) .. controls (-.1,-.2) and (-.2,0) .. (-.2,.15);
\draw [blue] (0,.1) node[left]{\scriptsize $\beta_2$};
\draw [blue,thick,dashed] (0,.5) .. controls (.1,.5) and (.2,.3) .. (.2,.15) .. controls (.2,0) and (.1,-.2) .. (0,-.2);
\end{scope}
\begin{scope}[yshift=-.5cm]
\draw (0,-.8) .. controls (1,-.8) and (1,-2.3) .. (1.65,-2.3) .. controls (2.1,-2.3) and (1,-.2) .. (0,-.2) .. controls (-1,-.2) and (-2.1,-2.3) .. (-1.65,-2.3) .. controls (-1,-2.3) and (-1,-.8) .. (0,-.8);
\draw (-2.9,-2.3) .. controls (-2.4,-2.3) and (-2,.5) .. (0,.5) .. controls (2,.5) and (2.4,-2.3) .. (2.9,-2.3);
\draw (0,-2) circle (.3);
\draw [red,<-] (.3,-2) arc (0:180:.3);
\draw [red] (.2,-2) node[right]{\scriptsize $\alpha_1$};
\draw [blue,thick] (0,-.8) .. controls (-.1,-.8) and (-.2,-1) .. (-.2,-1.25);
\draw [blue,thick,-<] (0,-1.7) .. controls (-.1,-1.7) and (-.2,-1.35) .. (-.2,-1.25);
\draw [blue] (-.05,-1.4) node[left]{\scriptsize $\beta^{\prime}_1$};
\draw [blue,thick,dashed] (0,-.8) .. controls (.1,-.8) and (.2,-1.05) .. (.2,-1.25) .. controls (.2,-1.45) and (.1,-1.7) .. (0,-1.7);
\draw [blue,thick] (0,.5) .. controls (-.1,.5) and (-.1,.3) .. (-.2,.15);
\draw [blue,thick,->] (0,-.2) .. controls (-.1,-.2) and (-.2,0) .. (-.2,.15);
\draw [blue] (0,.1) node[left]{\scriptsize $\beta_2$};
\draw [blue,thick,dashed] (0,.5) .. controls (.1,.5) and (.2,.3) .. (.2,.15) .. controls (.2,0) and (.1,-.2) .. (0,-.2);
\draw [very thick,->] (-11.6,-1.8)-- (-11,-1.8);
\draw [very thick,->] (-4.2,-1.8)-- (-3.6,-1.8);
\draw [very thick,->] (3.4,-1.8) -- (4,-1.8);
\end{scope}
\begin{scope}[xshift=6.7cm,yshift=-4.3cm]
\draw (1.3,1.5) .. controls (2,1.5) and (2,3.2) .. (3,3.2);
\draw (4.7,1.5) .. controls (4,1.5) and (4,3.2) .. (3,3.2);
\draw [red] (2.4,2) node[right]{\scriptsize $\alpha_2$};
\draw (3,2) circle (.5);
\draw [red,thick,->] (3.5,2) arc (0:180:.5);
\draw [blue] (2.8,3.44) node{\scriptsize $\beta_2$};
\draw [blue,thick] (3,2.5) .. controls (2.9,2.5) and (2.8,2.7) .. (2.8,2.85);
\draw [blue,thick,-<] (3,3.2) .. controls (2.9,3.2) and (2.8,3) .. (2.8,2.85);
\draw [blue,thick,dashed] (3,3.2) .. controls (3.1,3.2) and (3.2,3) .. (3.2,2.85) .. controls (3.2,2.7) and (3.1,2.5) .. (3,2.5);
\end{scope}
\begin{scope}[xshift=3.3cm,yshift=-4.3cm]
\draw (1.3,1.5) .. controls (2,1.5) and (2,3.2) .. (3,3.2);
\draw (4.7,1.5) .. controls (4,1.5) and (4,3.2) .. (3,3.2) ;
\draw (3,2) circle (.5);
\draw [red,thick,<-] (3.5,2) arc (0:180:.5);
\draw [red] (3.7,1.9) node[left]{\scriptsize $\alpha_1$};
\draw [blue,thick,->] (3,3.2) .. controls (2.9,3.2) and (2.8,3) .. (2.8,2.85);
\draw [blue,thick] (3,2.5) .. controls (2.9,2.5) and (2.8,2.7) .. (2.8,2.85);
\draw [blue] (3,2.65) node[left]{\scriptsize $\beta^{\prime}_1$};
\draw [blue,thick,dashed] (3,3.2) .. controls (3.1,3.2) and (3.2,3) .. (3.2,2.85) .. controls (3.2,2.7) and (3.1,2.5) .. (3,2.5);
\end{scope}
\end{tikzpicture}
\caption{Handle slide}
\label{fighansli}
\end{center}
\end{figure}

See $(-\ComX)$ in $\handB$ as the upward vertical field in the first picture of Figure~\ref{fighansli}. This field is an outward normal to $\handB$ except around $\pointw$ that is not shown in our figures, and around the crossings of $\matchingo$, more precisely on the gray disks $D_i$ shown in Figure~\ref{figgammaidisk}. Inside the disks $D_i$, $(-\ComX)$ is an inward normal to $\handB$. On the boundary of this disk, it is tangent to the surface.
Our homotopy will fix $(-\ComX)$ in the neighborhood of $\pointw$ where $(-\ComX)$ is not an outward normal to $\handB$, and the locus of $\partial \handB$ where $Y_t$ is a positive (resp. negative) normal to $\handB$ will not depend on $t$. Thus this homotopy can be canonically modified (without changing the locus where $Y_t$ is a positive (resp. negative) normal to $\handB$) so that $Y_t$ is fixed on $\partial \handB$.

Observe that there is no loss in assuming that the path $\gamma$ from $\beta_1$ to $\beta_2$ that parametrizes the handle slide is as in  the first picture of Figure~\ref{fighansli}. The next pictures describe various positions of $\handB$ under an ambient isotopy $(h_t)_{t\in [0,1]}$ of $\RR^3$ that first moves the handle of $\beta_2$ upward (second picture), slides it over the handle of $\beta_1$ (fourth picture), moves the handle of $\beta_1$ upward (fifth picture) and replaces the slid foot of $H_2$ in its original position by letting it slide away from the handles (last picture). The isotopy $(h_t)_{t\in [0,1]}$ starts with $h_0$ that is the Identity and finishes with a homeomorphism $h_1$ of $\RR^3$ that maps $\handB$ to itself.
Let $\vec{N}$ be the upward vector field of $\RR^3$. Then $(h_t)^{-1}_{\ast}(\vec{N}_{|h_t(\handB)})$ defines a homotopy of nowhere zero vector fields from $Y_0=(-\ComX)$ and $Y_1=(-\ComX^{\prime})$ on $\handB$ that behaves as wanted on the boundary.

\eop

Let us start with common preliminaries for the proofs of the remaining three lemmas.

Set $\CJ^{\prime}_{2i}=\CJ_{2i}-\CJ_{1i}$.
For any interval $I$ of an $\alpha_i$,
$\langle I, \beta^{\prime}_1\rangle=\langle I, \beta_1 + \beta_2^+\rangle$
and 
$$\langle I, \CJ_{1i}\beta^{\prime}_1 + \CJ^{\prime}_{2i}\beta_2\rangle=\langle I,\CJ_{1i}\beta_1 + \CJ_{2i}\beta_2 +\CJ_{1i}(\beta_2^+ -\beta_2)\rangle.$$
Set $\CJ^{\prime}_{ji}=\CJ_{ji}$ for any $(j,i)$ such that $j \neq 2$.
Every quantity associated with $\CD^{\prime}$ will have a prime superscript.
Our definitions of the $\CJ^{\prime}_{ji}$ ensure that
$$\langle \alpha_k,\sum_{j} \CJ^{\prime}_{ji}\beta^{\prime}_j\rangle  = \langle \alpha_k,\sum_{j}\CJ_{ji} \beta_j\rangle =\delta_{ik},$$ for any $i$ and $k$, as required.

Let $\CaC_2$ be the set of crossings of $\CD$ on $\beta_2$, and let $\CaC_2^+$ be the set of crossings of $\CD^{\prime}$ on $\beta_2^{+}$, $\CaC_2^+$ is in natural one-to-one correspondence with $\CaC_2$ and the crossing of $\CaC_2^+$ that corresponds to $c$ will be denoted by $c^+$.
$$\CaC^{\prime}=\CaC \cup \CaC_2^+.$$

\begin{figure}[h]
\begin{center}
\begin{tikzpicture}
\useasboundingbox (-2,-.2) rectangle (6.4,1.9);
\begin{scope}[xshift=-5cm,yshift=.5cm]
\draw [->] (.8,.1) -- node[very near end, below]{\scriptsize $\beta_2$} (0,.1);
\draw [->] (0,.7) -- node[very near end, above]{\scriptsize $\beta_1$} (.8,.7);
\draw (1.2,.4) node{$\rightarrow$};
\draw [->] (2.6,.1) -- node[very near end, below]{\scriptsize $\beta_2$} (1.6,.1);
\draw [->] (1.6,.7) .. controls (2,.7) and (2,.3) .. (1.6,.3);
\draw [->] (2.6,.3) .. controls (2.2,.3) and (2.2,.7) .. node[very near end, above]{\scriptsize $\beta^{\prime}_1$} (2.6,.7); (.8,.5);
\end{scope}
\begin{scope}[xshift=-1cm]
\draw [->] (1,.5) arc (-90:270:.4);
\draw (1,0) node{\scriptsize $\alpha^{\prime}_2$};
\draw [->] (3.4,.5) arc (270:-90:.4);
\draw (3.4,0) node{\scriptsize $\alpha^{\prime \prime}_2$};
\draw (1.55,.7) node{\scriptsize $\matchm_2$} (2.85,.7) node{\scriptsize $\matchm_2$} (2.55,1.35) node{\scriptsize $\beta^{\prime}_1$} (1.7,1.35) node{\scriptsize $\beta^{\prime}_1$} (2.2,1.8) node{\scriptsize $\sigma(\matchm_2)=1$};
\draw [->] (1.7,.9) -- (1.1,.9);
\draw [->] (3.15,.9) -- (2.55,.9);
\draw [->] (3.5,1).. controls (4.1,1) and (3.8,1.45) .. (3.4,1.45) .. controls (3.2,1.45) and (3,1.05) .. (2.55,1.05);
\draw [<-] (.9,1).. controls (.3,1) and (.6,1.45) .. (1,1.45) .. controls (1.2,1.45) and (1.4,1.05) .. (1.7,1.05);
\end{scope}
\begin{scope}[xshift=4cm]
\draw [->] (1,.5) arc (-90:270:.4);
\draw (1,0) node{\scriptsize $\alpha^{\prime}_2$};
\draw [->] (3.4,.5) arc (270:-90:.4);
\draw (3.4,0) node{\scriptsize $\alpha^{\prime \prime}_2$};
\draw (1.55,1.1) node{\scriptsize $\matchm_2$} (2.85,1.1) node{\scriptsize $\matchm_2$} (2.55,.45) node{\scriptsize $\beta^{\prime}_1$} (1.7,.45) node{\scriptsize $\beta^{\prime}_1$} (2.2,1.8) node{\scriptsize $\sigma(\matchm_2)=-1$};
\draw [<-] (1.7,.9) -- (1.1,.9);
\draw [<-] (3.15,.9) -- (2.55,.9);
\draw [<-] (3.5,.8).. controls (4.1,.8) and (3.8,.35) .. (3.4,.35) .. controls (3.2,.35) and (3,.75) .. (2.55,.75);
\draw [->] (.9,.8).. controls (.3,.8) and (.6,.35) .. (1,.35) .. controls (1.2,.35) and (1.4,.75) .. (1.7,.75);
\end{scope}

\end{tikzpicture}
\caption{Variation of $\tilth$}
\label{figvartheta}
\end{center}
\end{figure}
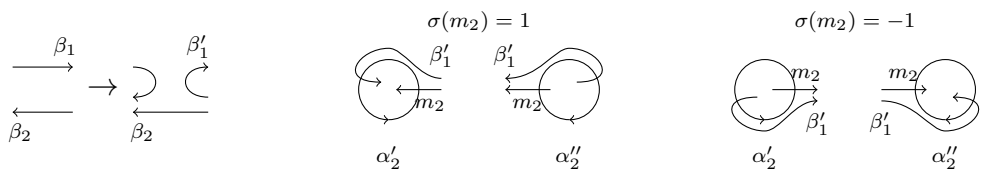

\noindent{\sc Proof of Lemma~\ref{lemhandleslidetheta}:}
Without loss, assume that $\beta_2$ goes from right to left at the place of the band sum as in Figure~\ref{figvartheta}. Then $\beta_1$ is above $\beta_2$ and it goes from left to right. Thus after the band sum, the degree of $\beta^{\prime}_1$ is increased by $(-1/2)$ before and after $\beta_2^+$ and by $(1/2)$ before and after $\matchm_2$.

Therefore $\tilth^{\prime}(\beta^{\prime}_1)=\tilth(\beta_1)+\tilth(\beta_2)$, and, for any $i$,
$$\sum_{j=1}^g\CJ^{\prime}_{ji} \tilth^{\prime}(\beta^{\prime}_j) = \sum_{j=1}^g\CJ_{ji} \tilth(\beta_j).$$
Then for any $c$ that is not in $[e^+,\matchm_1[_{\beta^{\prime}_1}$, $\tilth^{\prime}(c)=\tilth(c)$.
Since $\tilth(\beta_2)-\sum_{(r,s) \in \underline{g}^2}\CJ_{sr}\langle \alpha_r,\beta_2\rangle\tilth(\beta_{s})=0$, for any $c \in [e^+,\matchm_1[_{\beta^{\prime}_1} \setminus \beta_2^+$, 
$\tilth^{\prime}(c)=\tilth(c)$, too,
 so that
$$e(\CD^{\prime},\pointw,\matchingo)-e(\CD,\pointw,\matchingo)=\sum_{c \in \beta_2}\CJ_{1i(c)}\sigma(c)(\tilth^{\prime}(c^+)-\tilth(c)).$$

For $c\in \beta_2$, $$\tilth^{\prime}(c^+)-\tilth^{\prime}(e^+)
=\left\{\begin{array}{ll}\tilth(c)
-\tilth(e)\;&\mbox{if}\; c \in ]e,\matchm_2[\\
\tilth(c)
-\tilth(e)+1\;&\mbox{if}\; c \in ]\matchm_2,e[\\
\tilth(c)
-\tilth(e)+\frac12\;&\mbox{if}\; c =\matchm_2.\end{array}
\right.$$
\eop

For the remaining two lemmas, for any $2$--cycle $\twocycG=\sum_{(c,d) \in (\CaC^{\prime})^2}g_{c d} (\gamma(c) \times \gamma(d)_{\parallel})$ of $M^2$, we compute $\ell^{(2)}(\twocycG)$ with Proposition~\ref{propevalbis}
with
$$\ell(c,d)=\langle [\pointp(\alpha(c)),c\halfbra_{\alpha}, [\pointp(\beta(d)),d\halfbra_{\beta}\rangle -\sum_{(i,j) \in \underline{g}^2}\CJ^{\prime}_{ji} \langle [\pointp(\alpha(c)),c\halfbra_{\alpha},\beta^{\prime}_j \rangle \langle \alpha_i ,[\pointp(\beta(d)),d\halfbra_{\beta} \rangle$$
where $\pointp(\beta_2)=e$ and $\pointp(\beta^{\prime}_1)$ is the first crossing of $\beta^{\prime}_1$ after $\beta^{+}_2$ on $\beta_1$, the $\pointp(\alpha_i)$ are not on $\beta^{+}_2$, and, if $\pointp(\alpha_i) \in \beta_2$, then $\sigma(\pointp(\alpha_i))=1$ (up to changing the orientation of $\alpha_i$). This map $\ell^{(2)}$ may be used for any $2$--cycle $\twocycG=\sum_{(c,d) \in \CaC^2}g_{c d} (\gamma(c) \times \gamma(d)_{\parallel})$ of $M^2$, as well, and we use it.

\begin{lemma} \label{lemellprime}
Recall $\CaC^{\prime}=\CaC \cup \CaC_2^+.$
Let $(c,d) \in \CaC^2$.
If $c \in \beta_2$, then
 $$\ell(c^+,d)-\ell(c,d) + \frac12 \sum_{i=1}^g \CJ_{2i} \langle \alpha_i ,[\pointp(\beta(d)),d\halfbra_{\beta} \rangle =\left\{\begin{array}{ll} 0 &\mbox{if} \;d \notin \beta_2\\
0  &\mbox{if} \;d \in \beta_2 \;\mbox{and}\; c \notin [e,d]_{\beta}\\
\frac12  &\mbox{if} \;d \in \beta_2 \;\mbox{and}\; c \in [e,d[_{\beta}\\
\frac14  &\mbox{if} \;c=d.\\
\end{array}\right.  $$
If $d \in \beta_2$, then $$\ell(c,d^+)-\ell(c,d)=\left\{\begin{array}{ll}-\frac12\;&\mbox{if}\; c \in [e,d[_{\beta}\\
-\frac14\;&\mbox{if}\; c =d\\
0\;&\mbox{if}\; c \notin [e,d]_{\beta}.\end{array}\right.
$$
If $(c,d) \in \CaC_2^2$, then $$\ell(c^+,d^+)-\ell(c^+,d)=\ell(c,d^+)-\ell(c,d).
$$
\end{lemma}
\bp Let $(c,d) \in \CaC^2$.
Assume $c \in \beta_2$.
If $\sigma(c)=1$, then
$$\begin{array}{ll}\ell(c^+,d)-\ell(c,d)&=\langle \halfbra c,c^+\halfbra_{\alpha}, [\pointp(\beta(d)),d\halfbra_{\beta}\rangle -\sum_{(i,j) \in \underline{g}^2}\CJ^{\prime}_{ji} \langle \halfbra c,c^+\halfbra_{\alpha},\beta^{\prime}_j \rangle \langle \alpha_i ,[\pointp(\beta(d)),d\halfbra_{\beta} \rangle\\
   &=\langle \halfbra c,c^+\halfbra_{\alpha}, [\pointp(\beta(d)),d\halfbra_{\beta}\rangle- \frac12 \sum_{i=1}^g (\CJ^{\prime}_{2i}+\CJ^{\prime}_{1i}) \langle \alpha_i ,[\pointp(\beta(d)),d\halfbra_{\beta} \rangle.
  \end{array}
$$
If $\sigma(c)=-1$,
$$\ell(c^+,d)-\ell(c,d)=-\langle \halfbra c^+,c\halfbra_{\alpha}, [\pointp(\beta(d)),d\halfbra_{\beta}\rangle - \frac12 \sum_{i=1}^g \CJ_{2i} \langle \alpha_i ,[\pointp(\beta(d)),d\halfbra_{\beta} \rangle.$$
Let $d \in \beta_2$. For any interval $I$ of an $\alpha_i$, $\langle I,[\pointp(\beta^{\prime}_1),d^+\halfbra_{\beta}\rangle=\langle I,\beta_1+ [e^+,d^+\halfbra_{\beta}\rangle$.\\
If $c \notin \beta_2$, $$\begin{array}{ll}\ell(c,d^+)-\ell(c,d)&=\langle [\pointp(\alpha(c)),c\halfbra_{\alpha}, \beta_1 \rangle -\sum_{(i,j) \in \underline{g}^2}\CJ^{\prime}_{ji} \langle [\pointp(\alpha(c)),c\halfbra_{\alpha},\beta^{\prime}_j \rangle \langle \alpha_i , \beta^{\prime}_1 - \beta^{\prime}_2\rangle\\
&=\langle [\pointp(\alpha(c)),c\halfbra_{\alpha}, \beta_1 \rangle -\sum_{j \in  \underline{g}}(\delta_{j1}-\delta_{j2}) \langle [\pointp(\alpha(c)),c\halfbra_{\alpha},\beta^{\prime}_j \rangle\\ 
&=\langle [\pointp(\alpha(c)),c\halfbra_{\alpha}, \beta_1 \rangle - \langle [\pointp(\alpha(c)),c\halfbra_{\alpha},\beta^{\prime}_1 - \beta^{\prime}_2\rangle\\
&=0.
\end{array}$$
When $c \in \beta_2$, we similarly get
$$\ell(c,d^+)-\ell(c,d)=\langle [\pointp(\alpha(c)),c\halfbra_{\alpha},  [e^+,d^+\halfbra_{\beta} -[e,d\halfbra_{\beta} \rangle=\left\{\begin{array}{ll}-\frac12\;&\mbox{if}\; c \in [e,d[_{\beta}\\
-\frac14\;&\mbox{if}\; c =d\\
0\;&\mbox{if}\; c \notin [e,d]_{\beta}.\end{array}\right.
$$
Since $\langle [\pointp(\alpha(c)),c\halfbra_{\alpha},  [e^+,d^+\halfbra_{\beta} -[e,d\halfbra_{\beta} \rangle=\langle [\pointp(\alpha(c)),c^+\halfbra_{\alpha}, [e^+,d^+\halfbra_{\beta} -[e,d\halfbra_{\beta} \rangle$,
$\ell(c^+,d^+)-\ell(c^+,d)=\ell(c,d^+)-\ell(c,d).
$
\eop

\noindent{\sc Proof of Lemma~\ref{lemhandleslidelk}:}
Set $L=\Link(\CD,\matchingo)=\sum_{i=1}^g \gamma_i- \sum_{c\in \CaC} \CJ_{j(c)i(c)} \sigma(c) \gamma(c)$
and $L^{\prime}=\Link(\CD^{\prime},\matchingo)$.
Then $$L^{\prime}-L=\sum_{c\in \CaC_2}\CJ_{1i(c)}\sigma(c)(\gamma(c)-\gamma(c^+))$$
is a cycle and $$lk(\Link^{\prime}, \Link^{\prime}_{\parallel})-
lk(\Link, \Link_{\parallel})= \ell((L^{\prime}-L) \times (L^{\prime}-L)) + 2 \ell((L^{\prime}-L) \times L)$$
thanks to the symmetry of the linking number in Proposition~\ref{propevallk}.

The last assertion of Lemma~\ref{lemellprime} guarantees that
$$\ell((L^{\prime}-L)\times (L^{\prime}-L))=0.$$

Now, $\ell((L^{\prime}-L) \times L)
=\ell_1+\ell_2$
with
$$\ell_1=\sum_{c\in \CaC_2,i\in \underline{g}}\CJ_{1i(c)}\sigma(c)(\ell(c,\matchm_i)-\ell(c^+,\matchm_i))$$
where $\matchingo = \{\matchm_i\}_{i \in \underline{g}}$ and $\matchm_i \in \alpha_i \cap \beta_i$ and
$$\ell_2=\sum_{c\in \beta_2,d\in \CaC}\CJ_{1i(c)}\sigma(c)\CJ_{j(d)i(d)} \sigma(d)(\ell(c^+,d)-\ell(c,d)).$$

Since the part $\left(\frac12 \sum_{i=1}^g \CJ_{2i} \langle \alpha_i ,[\pointp(\beta(d)),d\halfbra_{\beta} \rangle\right)$ that occurs in the expressions of
$(\ell(c^+,d)-\ell(c,d))$ in Lemma~\ref{lemellprime} is independent of $c$, the factor
$\sum_{c\in \beta_2}\CJ_{1i(c)}\sigma(c)$ that vanishes makes it disappear so that
$$\ell((L^{\prime}-L) \times L)
=\tilde{\ell}_1+\tilde{\ell}_2$$
where
$$\tilde{\ell}_1=-\frac{1}{2}\left(\sum_{c \in [e,\matchm_2[_{\beta}}\sigma(c)\CJ_{1i(c)} + \frac{1}{2}\sigma(\matchm_2)\CJ_{12}\right) = -\frac{1}{2}\sum_{c \in [e,\matchm_2\halfbra_{\beta}}\sigma(c)\CJ_{1i(c)}
$$ and
$$\tilde{\ell}_2=\frac12\sum_{d \in \beta_2, c \in [e,d\halfbra_{\beta}}\sigma(c)\sigma(d)\CJ_{1i(c)}\CJ_{2i(d)}.$$
\eop

\noindent{\sc Proof of Lemma~\ref{lemhandleslideCh}:}
Recall $$\ell_2(\CD)=\sum_{(c,d) \in \CaC^2}\CJ_{j(c)i(d)}\CJ_{j(d)i(c)}\sigma(c)\sigma(d) \ell(c,d) -\sum_{c \in \CaC} \CJ_{j(c)i(c)}\sigma(c)\ell(c,c).$$
Define the projection $q \colon \CaC^{\prime} \to \CaC$ such that $q(c)=c$ if $c \in \CaC$ and
$q(c^+)=c$ if $c \in \beta_2$.
Since a crossing $c$ of $\beta_2$ gives rise to two crossings $c$ and $c^+$ of $\CaC^{\prime}$ whose coefficients $\CJ^{\prime}_{2r}$ and $\CJ^{\prime}_{1r}$ add up to $\CJ_{2r}$,
$$\ell_2(\CD)=\sum_{(c,d) \in (\CaC^{\prime})^2}\CJ^{\prime}_{j(c)i(d)}\CJ^{\prime}_{j(d)i(c)}\sigma(c)\sigma(d) \ell(q(c),q(d)) -\sum_{c \in \CaC^{\prime}} \CJ^{\prime}_{j(c)i(c)}\sigma(c)\ell(q(c),q(c))$$
so that
$$\begin{array}{ll}\ell_2(\CD^{\prime})-\ell_2(\CD)=&\sum_{(c,d) \in (\CaC^{\prime})^2}\CJ^{\prime}_{j(c)i(d)}\CJ^{\prime}_{j(d)i(c)}\sigma(c)\sigma(d) (\ell(c,d)-\ell(q(c),q(d)))\\& -\sum_{c \in \CaC_2} \CJ^{\prime}_{1i(c)}\sigma(c)(\ell(c^+,c^+)-\ell(c,c)).\end{array}$$
Write $\ell(c,d)-\ell(q(c),q(d))=\ell(c,d)-\ell(c,q(d))+\ell(c,q(d))-\ell(q(c),q(d))$.
$$\ell(c,d^+)-\ell(c,d)
=\ell(q(c),d^+)-\ell(q(c),d)
=\left\{\begin{array}{ll}-\frac12\;&\mbox{if}\; q(c) \in [e,d[_{\beta}\\
-\frac14\;&\mbox{if}\; q(c) =d\\
0\;&\mbox{if}\; q(c) \notin [e,d]_{\beta}\end{array}\right.
$$
for $d \in \CaC_2$ so that
$$\begin{array}{ll}\ell_2(\CD^{\prime})-\ell_2(\CD)=&-\frac12\sum_{d \in \CaC_2}\sum_{c \in [e,d\halfbra_{\beta}}(\CJ^{\prime}_{2i(d)}+\CJ^{\prime}_{1i(d)})\CJ^{\prime}_{1i(c)}\sigma(c)\sigma(d)+ A
\\
&+\frac14\sum_{c \in \CaC_2}\CJ^{\prime}_{1i(c)}\sigma(c)-\sum_{c \in \CaC_2} \CJ^{\prime}_{1i(c)}\sigma(c)(\ell(c^+,c)-\ell(c,c))\end{array}$$
where
$$\begin{array}{lll}A&=&
\sum_{(c,d) \in (\CaC^{\prime})^2}\CJ^{\prime}_{j(c)i(d)}\CJ^{\prime}_{j(d)i(c)}\sigma(c)\sigma(d) (\ell(c,q(d))-\ell(q(c),q(d)))\\
&=&
\sum_{(c,d) \in \CaC^{\prime} \times \CaC }\CJ^{\prime}_{j(c)i(d)}\CJ_{j(d)i(c)}\sigma(c)\sigma(d) (\ell(c,d)-\ell(q(c),d))\\
&=&\sum_{(c,d) \in \CaC_2\times \CaC}\CJ^{\prime}_{1i(d)}\CJ_{j(d)i(c)}\sigma(c)\sigma(d) (\ell(c^+,d)-\ell(c,d))\\
&=&-\frac12\sum_{(c,d) \in \CaC_2\times \CaC}\CJ^{\prime}_{1i(d)}\CJ_{j(d)i(c)}\sigma(c)\sigma(d) ( \sum_{i=1}^g \CJ_{2i} \langle \alpha_i ,[\pointp(\beta(d)),d\halfbra_{\beta} \rangle)\\
&&+\frac12\sum_{d \in \CaC_2,c \in [e,d\halfbra_{\beta}}\CJ^{\prime}_{1i(d)}\CJ_{2i(c)}\sigma(c)\sigma(d)\\
&=&-\frac12 \sum_{i=1}^g\sum_{d \in \CaC_2}\CJ^{\prime}_{1i(d)}\sigma(d)\CJ_{2i}\langle \alpha_i ,[\pointp(\beta(d)),d\halfbra_{\beta} \rangle\\&&
+\frac12\sum_{d \in \CaC_2,c \in [e,d\halfbra_{\beta}}\CJ^{\prime}_{1i(d)}\CJ_{2i(c)}\sigma(c)\sigma(d)\\
&=&0.\end{array}
$$
$$\begin{array}{lll}\sum_{c \in \CaC_2} \CJ_{1i(c)}\sigma(c)(\ell(c^+,c)-\ell(c,c))&=&-\frac12\sum_{c \in \CaC_2} \CJ_{1i(c)}\sigma(c)( \sum_{i=1}^g \CJ_{2i} \langle \alpha_i ,[\pointp(\beta(c)),c\halfbra_{\beta} \rangle)\\
   &&+\frac14\sum_{c \in \CaC_2} \CJ_{1i(c)}\sigma(c)\\
&=&-\frac12\sum_{(c,d) \in \CaC_2^2; d \in[e,c\halfbra_{\beta} } \CJ_{2i(d)}\CJ_{1i(c)}\sigma(c)  \sigma(d).
  \end{array}
$$
$$\begin{array}{lll}\ell_2(\CD^{\prime})-\ell_2(\CD)&=&-\frac12\sum_{d \in \CaC_2}\sum_{c \in [e,d\halfbra_{\beta}}\CJ_{2i(d)}\CJ_{1i(c)}\sigma(c)\sigma(d)
\\&&+\frac12\sum_{(c,d) \in \CaC_2^2; d \in[e,c\halfbra_{\beta} } \CJ_{2i(d)}\CJ_{1i(c)}\sigma(c)  \sigma(d).\end{array}$$
For $r,s\in\underline{g}$, set
$$V_{r,s}=\sum_{c \in \CaC_2}\sum_{d \in[e,c\halfbra_{\beta}}\CJ_{ri(d)}\CJ_{si(c)}\sigma(c)\sigma(d).$$
Note that $V_{r,s}+V_{s,r}=\delta_{r2}\delta_{s2}$ (recall the argument after the statement of Lemma~\ref{lemhandleslidetheta}).
Thus $\ell_2(\CD^{\prime})-\ell_2(\CD)=\frac12 (V_{2,1}- V_{1,2})=V_{2,1}$.
\eop

\subsection{Connected sums and stabilizations}

The previous subsections guarantee that $\tilde{\lambda}$ is an invariant of Heegaard decompositions.

\begin{figure}[h]
\begin{center}
\begin{tikzpicture}
\useasboundingbox (-7,-.9) rectangle (7,.7);
\draw (-2,-.8) rectangle (2,.6) (0,.2) node{\scriptsize $\beta_1$};
\draw [->] (-1,-.3) arc (-90:270:.3);
\draw (-1,-.55) node{\scriptsize $\alpha^{\prime}_1$};
\draw [->] (1,-.3) arc (270:-90:.3);
\draw (1,-.55) node{\scriptsize $\alpha^{\prime\prime}_1$};
\draw [->] (.7,0) -- (0,0);
\draw (0,0) -- (-.7,0);
\end{tikzpicture}
\caption{Genus one Heegaard diagram of $S^3$}
\label{figrectsimp}
\end{center}
\end{figure}
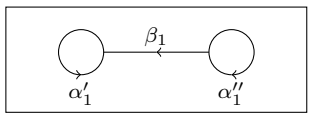

\begin{lemma} \label{lemlambst}
 Let $$S^3=T_A \cup_{\partial T_A\sim -\partial T_B} T_B$$ be the genus one decomposition of $S^3$ as a union of two solid tori $T_A$ and $T_B$ glued along their boundaries so that the meridian $\alpha_1$ of $T_A$ meets the meridian $\beta_1$ of $T_B$ once.
$$\tilde{\lambda}(T_A \cup_{\partial T_A\sim -\partial T_A} T_B)=0.$$
\end{lemma}
\bp
Orient $\alpha_1$ and $\beta_1$ so that $\langle \alpha_1, \beta_1 \rangle_{\partial T_A}=1$.
Then $\CJ_{11}=1$. Let $\matchingo=\{\alpha_1 \cap \beta_1\}$ be the unique matching. Let $\pointw$ be a point of the connected $\partial T_A \setminus (\alpha_1 \cup \beta_1)$. Then the intersection of $T_A$ with $B^3$ can be embedded in $B^3$ as in Figure~\ref{figgammai}, and $T_B$ is its complement in $B^3$. In particular, $\ComX(\pointw,\matchingo)$ is the vertical field of $\RR^3$ and
$p_1(\ComX(\pointw,\matchingo))=0$. The rectangular picture of the Heegaard diagram is simply Figure~\ref{figrectsimp} so that $e(\CD,\pointw,\matchingo)=0$.
Since $\twocycG(\CD)=\emptyset$ and $\Link(\CD,\matchingo)=\emptyset$, $\ell_2(\CD)=0$ and  $s_{\ell}(\CD,\matchingo)=0$.
\eop

The \emph{connected sum} $M \sharp M^{\prime}$ of two connected closed manifolds $M$ and $M^{\prime}$ of dimension $d$ is obtained by removing the interior of an open ball from $M$ and from $M^{\prime}$ and by gluing the obtained manifolds along their spherical boundaries $$M \sharp M^{\prime} = (M \setminus \mathring{B}^d) \cup_{S^{d-1}} (M^{\prime} \setminus \mathring{B}^{\prime d}).$$
When the manifolds are $3$--manifolds equipped with Heegaard decompositions $M = \handA \cup_{\partial \handA} \handB$ and $M^{\prime} = \handA^{\prime} \cup_{\partial \handA^{\prime}} \handB^{\prime}$, 
the \emph{connected sum} of the Heegaard decompositions is the Heegaard decomposition
$$M \sharp M^{\prime} = \handA \sharp_{\partial}  \handA^{\prime} \cup_{\partial \handA \sharp \partial \handA^{\prime}} \handB \sharp_{\partial}  \handB^{\prime}$$
where the open ball $\ballB$ (resp. $\ballB^{\prime}$) removed from $M$ (resp. from $M^{\prime}$) intersects the Heegaard surface $\partial \handA$ (resp. $\partial \handA^{\prime}$) as a properly embedded two dimensional disk that separates $\ballB$ into two half-balls $\mathring{\handA} \cap \ballB$ and $\mathring{\handB} \cap \ballB$ (resp. $\mathring{\handA}^{\prime} \cap \ballB^{\prime}$ and $\mathring{\handB}^{\prime} \cap \ballB^{\prime}$), the connected sum along the boundaries
$$\handA \sharp_{\partial}  \handA^{\prime}=\left(\handA \setminus (\handA \cap \mathring{\ballB})\right) \cup_{\handA \cap \partial \ballB \sim (-\handA^{\prime} \cap \partial \ballB^{\prime})}  \left(\handA^{\prime}\setminus (\handA^{\prime} \cap \mathring{\ballB}^{\prime})\right)$$ is homeomorphic to the manifold obtained by identifying $\handA$ and $ \handA^{\prime}$ along a two-dimensional disk of the boundary, and $\handB \sharp_{\partial}  \handB^{\prime}$ is defined similarly.

\begin{proposition} \label{propaddconnsum}
Under the hypotheses above, if $M$ and $M^{\prime}$ are rational homology spheres, then
$$ \tilde{\lambda}(\handA \sharp_{\partial}  \handA^{\prime} \cup_{\partial \handA \sharp \partial \handA^{\prime}} \handB \sharp_{\partial}  \handB^{\prime})= \tilde{\lambda}(\handA \cup_{\partial \handA} \handB )+\tilde{\lambda}( \handA^{\prime} \cup_{\partial \handA^{\prime}} \handB^{\prime})$$
\end{proposition}
\bp
When performing such a connected sum on manifolds equipped with Heegaard diagrams $\CD=(\partial \handA, (\alpha_i)_{i\in \underline{g}}, (\beta_j)_{j\in \underline{g}})$ and $\CD^{\prime}=(\partial \handA^{\prime}, (\alpha^{\prime}_i)_{i\in \underline{g^{\prime}}}, (\beta^{\prime}_j)_{j\in \underline{g^{\prime}}})$ and with exterior points $\pointw$ and $\pointw^{\prime}$ of $\CD$ and $\CD^{\prime}$, we
assume that the balls $D$ and $D^{\prime}$ meet the Heegaard surfaces inside the connected component of $\pointw$ or $\pointw^{\prime}$ outside the diagram curves, without loss, and we choose a basepoint $\pointw^{\prime \prime}$ in the corresponding region of $\partial \handA \sharp \partial \handA^{\prime}$.
Then we obtain the obvious Heegaard diagram $$\CD^{\prime \prime}=(\partial \handA \sharp \partial \handA^{\prime}, (\alpha^{\prime \prime}_i)_{i\in \underline{g^{\prime \prime}}}, (\beta^{\prime \prime}_j)_{j\in \underline{g^{\prime \prime}}})$$
where $g^{\prime \prime}=g+g^{\prime}$, $\alpha^{\prime \prime}_i=\alpha_i$ and $\beta^{\prime \prime}_i=\beta_i$ when $i \leq g$ and $\alpha^{\prime \prime}_i=\alpha^{\prime}_{i-g}$ and $\beta^{\prime \prime}_i=\beta^{\prime}_{i-g}$ when $i > g$,
with the associated intersection matrix and its inverse that are diagonal with respect to the two blocks corresponding to the former matrices associated with $\CD$ and $\CD^{\prime}$.

When $\CD$ and $\CD^{\prime}$ are furthermore equipped with matchings $\matchingo$ and $\matchingo^{\prime}$,
$\matchingo^{\prime \prime}=\matchingo \cup \matchingo^{\prime}$ is a matching for $\CD^{\prime \prime}$
and a rectangular figure for $(\CD^{\prime \prime},\pointw^{\prime \prime},\matchingo^{\prime \prime})$ similar to Figure~\ref{figplansplit} is obtained from the corresponding figures for $\CD$ and $\CD^{\prime}$ by juxtapositions of the two rectangles of $\CD$ and $\CD^{\prime}$.
In particular, 
$$e(\CD^{\prime \prime},\pointw^{\prime \prime},\matchingo^{\prime \prime})=e(\CD^{\prime},\pointw^{\prime},\matchingo^{\prime})+e(\CD,\pointw,\matchingo).$$
Furthermore, we can see $B_{M^{\prime \prime}}$ as the juxtaposition of two half-balls glued along a vertical disk equipped with the vertical field (over the intersection of the two rectangles above) such that the two half-balls are obtained from $B_{M}$ and $B_{M^{\prime}}$ by removing standard vertical half-balls equipped with the vertical field, so that the vector field $\ComX(\pointw^{\prime \prime},\matchingo^{\prime \prime})$ coincides with $\ComX(\pointw,\matchingo)$ on the remaining part of $B_M$ and with $\ComX(\pointw^{\prime},\matchingo^{\prime})$ on the remaining part of $B_{M^{\prime}}$. This makes clear that
$$p_1(\ComX(\pointw^{\prime \prime},\matchingo^{\prime \prime}))=p_1(\ComX(\pointw,\matchingo))+p_1(\ComX(\pointw^{\prime},\matchingo^{\prime})).$$
Now it is easy to observe that $\twocycG(\CD^{\prime \prime})=\twocycG(\CD)+\twocycG(\CD^{\prime})$, that $$\ell_2(\CD^{\prime \prime})=\ell_2(\CD)+\ell_2(\CD^{\prime}),$$ that
$\Link(\CD^{\prime \prime},\matchingo^{\prime \prime})=\Link(\CD,\matchingo)+\Link(\CD^{\prime},\matchingo^{\prime})$
and that
$$s_{\ell}(\CD^{\prime \prime},\matchingo^{\prime \prime})=s_{\ell}(\CD,\matchingo)+s_{\ell}(\CD^{\prime},\matchingo^{\prime}).$$
\eop

A connected sum of a Heegaard decomposition with the genus one decomposition of $S^3$ is called
a \emph{stabilization}.
A well-known Reidemeister-Singer theorem proved by Siebenmann in \cite{SiebenmannRS}, asserts that any two Heegaard decompositions of the same $3$-manifold become isomorphic after some stabilizations.
This Reidemeister-Singer theorem can also be proved using Cerf theory \cite{Cerfpseudo} as in \cite[Proposition~2.2]{OS1}.

Together with Proposition~\ref{propaddconnsum} and Lemma~\ref{lemlambst}, it implies that $\tilde{\lambda}$ does not depend on the Heegaard decomposition and allows us to prove the following theorem.

\begin{theorem}
 There exists a unique invariant $\tilde{\lambda}$ of $\QQ$-spheres
such that for any Heegaard diagram $\CD$ of a $\QQ$-sphere $M$, equipped with a matching $\matchingo$ and with an exterior point $\pointw$,
$$24\tilde{\lambda}(M)=4\ell_2(\CD)+4s_{\ell}(\CD,\matchingo)-4e(\CD,\pointw,\matchingo)-{p_1(\ComX(\pointw,\matchingo))}.$$
Furthermore, $\tilde{\lambda}$ satisfies the following properties.
\begin{itemize}
 \item For any two rational homology spheres $M_1$ and $M_2$,
$$\tilde{\lambda}(M_1 \sharp M_2)=\tilde{\lambda}(M_1) + \tilde{\lambda}(M_2).$$
\item For any rational homology sphere $M$, if $(-M)$ denotes the manifold $M$ equipped with the opposite orientation, then
$$\tilde{\lambda}(-M)=-\tilde{\lambda}(M).$$
\end{itemize}
\end{theorem}
\bp The invariance of $\tilde{\lambda}$ is already proved. 
Proposition~\ref{propaddconnsum} now implies that $\tilde{\lambda}$ is additive under connected sum. Reversing the orientation of $M$ reverses the orientation of the surface that contains a diagram $\CD$ of $M$. This changes the signs of the intersection points and reverses the sign of $\CJ$.
Thus $\Link(\CD,\matchingo)$, $\twocycG(\CD)$ and $\ComX(\pointw,\matchingo)$ are unchanged, while $\ell$ is changed to its opposite. Changing the orientation of the ambient manifold reverses the sign of $p_1$.
A rectangular diagram of $(-M)$ as in Figure~\ref{figplansplit} is obtained from the diagram of $M$ by a orthogonal symmetry that fixes a vertical line so that the $\tilth$ are changed to their opposites. Thus all the terms of the formula are multiplied by $(-1)$ when the orientation of $M$ is reversed.

\eop

\section{The Casson surgery formula for \texorpdfstring{$\tilde{\lambda}$}{lambda} }
\setcounter{equation}{0}
\label{secsurgcomb}

\subsection{The statement and its consequences}

In this section, we prove that $\tilde{\lambda}$ coincides with the Casson invariant for integral homology spheres by proving that it satisfies the same surgery formula.
More precisely, we prove the following theorem.

\begin{figure}[h]
\begin{center}
\begin{tikzpicture}
\useasboundingbox (-7,-2) rectangle (7,2.8);
\draw [thick,->] (-2.15,1) .. controls (-1.8,1) and (-1.7,0) .. (-2,0) .. controls (-2.3,0) .. (-2.3,1.3) .. controls  (-2.3,2.8) .. (-3.8,2.8) .. controls (-5,2.8) and (-5.3,2.4) .. (-5.3,0) .. controls (-5.3,-2) .. (0,-2);
\draw [thick,->] (-3.05,1.6) .. controls (-4.1,1.6) and (-4,0) .. (-4.4,0) .. controls (-4.7,0) .. (-4.7,1.3) .. controls (-4.7,1.7) and (-4.3,2.2) .. (-3.8,2.2) .. controls (-2.9,2.2) .. (-2.9,1.3) .. controls (-2.9,0) .. (-3.2,0) .. controls (-3.5,0) and (-3.5,1) .. (-3.05,1);
\draw [thick,->] (-.5,0) .. controls (-1.5,0) and (-1.1,1.6) .. (-2.15,1.6);
\draw [->] (-1.4,-.5) .. controls (-1.4,1) and (-1.7,1.3) .. (-2.15,1.3) (-3.05,1.3) .. controls (-3.6,1.3) and (-3.8,1) .. (-3.8,-.5) ..  controls (-3.8,-2) and (-1.4,-2) .. (-1.4,-.5);
\draw [->] (-2.6,-.3) .. controls (-2.6,2.5) and (-2.6,2.5) .. (-3.8,2.5) .. controls (-4.8,2.5) and (-5,2) .. (-5,-.3) ..  controls (-5,-1) and (-2.6,-1) .. (-2.6,-.3);
\draw (-2.6,-.5) node[right]{\scriptsize $z_2$} (-1.4,-.5) node[right]{\scriptsize $z_1$} (0,-2) node[above]{\scriptsize $K$} (0,0) node{\dots};
\begin{scope}[xshift=6.4cm]
\draw [thick,->] (-2.15,1) .. controls (-1.8,1) and (-1.7,0) .. (-2,0) .. controls (-2.3,0) .. (-2.3,1.3) .. controls  (-2.3,2.8) .. (-3.8,2.8) .. controls (-5,2.8) and (-5.3,2.4) .. (-5.3,1) .. controls (-5.3,0) .. (-5.9,0);
\draw [thick,->] (-3.05,1.6) .. controls (-4.1,1.6) and (-4,0) .. (-4.4,0) .. controls (-4.7,0) .. (-4.7,1.3) .. controls (-4.7,1.7) and (-4.3,2.2) .. (-3.8,2.2) .. controls (-2.9,2.2) .. (-2.9,1.3) .. controls (-2.9,0) .. (-3.2,0) .. controls (-3.5,0) and (-3.5,1) .. (-3.05,1);
\draw [thick,->] (-6.4,-2) .. controls (-1.1,-2) .. (-1.1,0) .. controls (-1.1,1.2) and (-1.6,1.6) .. (-2.15,1.6);
\draw [->] (-1.4,-.5) .. controls (-1.4,1) and (-1.7,1.3) .. (-2.15,1.3) (-3.05,1.3) .. controls (-3.6,1.3) and (-3.8,1) .. (-3.8,-.5) ..  controls (-3.8,-2) and (-1.4,-2) .. (-1.4,-.5);
\draw [->] (-2.6,-.3) .. controls (-2.6,2.5) and (-2.6,2.5) .. (-3.8,2.5) .. controls (-4.8,2.5) and (-5,2) .. (-5,-.3) ..  controls (-5,-1) and (-2.6,-1) .. (-2.6,-.3);
\draw (-2.5,-.3) node[left]{\scriptsize $z_{2g(\Sigma)}$} (-1.3,-.65) node[left]{\scriptsize $z_{2g(\Sigma)-1}$};
\end{scope}
\end{tikzpicture}
\caption{Curves on the surface $\Sigma$}
\label{figsigmasur}
\end{center}
\end{figure}
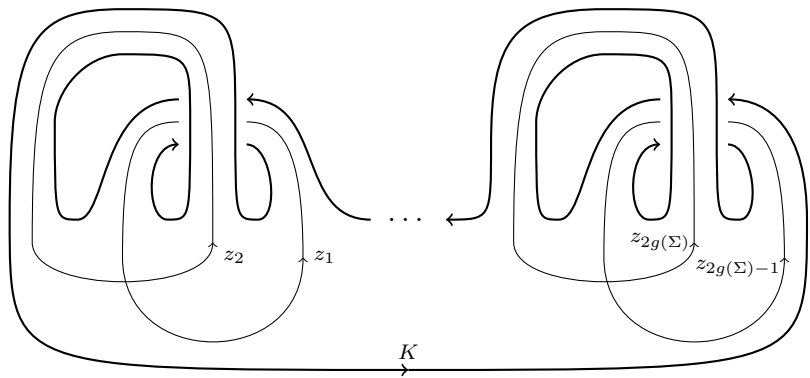

\begin{theorem}
\label{thmcassur}
Let $K$ be a null-homologous knot in a $\QQ$-sphere $M$.

Let $\Sigma$ be an oriented connected surface of genus $g(\Sigma)$ in $M$ bounded by $K$
such that the closure of the complement of a collar $$\handA =\Sigma \times [-1,1]$$ 
of $\Sigma=\Sigma \times \{0\}$ in $M$ is homeomorphic to a handlebody $\handB$.
This gives rise to the Heegaard decomposition $$M=\handA \cup_{\Psi_M} \handB $$ where $\Psi_M$ is an orientation reversing diffeomorphism from $\partial \handB$ to $\partial \handA$.
Let $M(K)$ be the manifold obtained from $M$ by surgery of coefficient $1$ along $K$ that can be defined by its Heegaard decomposition
$$M(K)=\handA \cup_{\Psi_M \circ \twist_K} \handB$$
where $\twist_K$ is the right-handed Dehn twist of $(-\partial \handB)$ about $K$.
Let $g=2g(\Sigma)$.
Let $(z_i)_{i \in \underline{g}}$ be closed curves of $\Sigma=\Sigma \times \{0\}$ that form a geometric symplectic basis of $H_1(\Sigma)$ as in Figure~\ref{figsigmasur}, and let $z_i^+=z_i \times \{1\}$. For any $i\in \underline{g(\Sigma)}$, $\langle z_{2i-1},z_{2i} \rangle=1$. 
Then 
$$\tilde{\lambda}(M(K))-\tilde{\lambda}(M)=\sum_{(i,r) \in \underline{g(\Sigma)}^2}\left(lk(z_{2i}^+,z_{2r})lk(z_{2i-1}^+,z_{2r-1})-lk(z_{2i}^+,z_{2r-1})lk(z_{2i-1}^+,z_{2r})\right).$$
\end{theorem}

We will prove the theorem exactly as it is stated.
A \emph{Seifert surface} $\Sigma$ of $K$ as in the statement is said to be \emph{unknotted}. It is well-known that any Seifert surface can be transformed to an unknotted one by adding some tubes (to remove unwanted $2$-handles from its exterior). (See \cite[Lemme 5.1]{mar}, \cite[p.84]{akmc} or \cite[Lemme 4.1]{gm} in the original surveys \cite{mar,akmc,gm} of the Casson invariant, for example.)
Thus any null-homologous knot bounds an unknotted surface as in the statement.
The manifold $M(K;\frac{p}{q})$ obtained from $M$ by \emph{Dehn surgery} with coefficient $p/q$ along $K$, for two coprime integers $p$ and $q$, is usually defined as
$$M(K;\frac{p}{q})=\left(M \setminus \mathring{N}(K)\right) \cup_{\partial N(K) \sim \partial D^2 \times S^1} \left(D^2 \times S^1\right)$$ 
where $N(K)$ is a tubular neighborhood of $K$, and the gluing homeomorphism from $\partial D^2 \times S^1$ to $\partial N(K)$ identifies the meridian $\partial D^2 \times \{x\}$ of $D^2 \times S^1$ with a curve homologous to $pm(K) +q\ell(K)$ where $m(K)$ is the meridian of $K$ such that $lk(m(K),K)=1$ and $\ell(K)$ is the curve parallel to $K$ such that $lk(\ell(K),K)=0$.
In our case, for $n \in \ZZ \setminus \{0\}$, the manifold $M(K;\frac1{n})$ obtained from $M$ by surgery of coefficient $\frac1{n}$ along $K$ can also be defined by its Heegaard decomposition
$$M(K;\frac1{n})=\handA \cup_{\Psi_M \circ \twist_K^n} \handB,$$
and it is easy to observe that the variation $(\tilde{\lambda}(M(K;\frac1{n}))-\tilde{\lambda}(M))$ can be deduced from the general knowledge of $(\tilde{\lambda}(M(K))-\tilde{\lambda}(M))$.
In our case, Theorem~\ref{thmcassur} implies that 
$$\tilde{\lambda}(M(K;\frac1{n}))-\tilde{\lambda}(M)=n(\tilde{\lambda}(M(K))-\tilde{\lambda}(M)).$$
In our proof, we will obtain the variation $\tilde{\lambda}(M(K))-\tilde{\lambda}(M)$ as it is stated, directly, so that our proof shows that \index{lambdaprime@$\lambda^{\prime}$}
$$\lambda^{\prime}=\sum_{(i,r) \in \underline{g(\Sigma)}^2}\left(lk(z_{2i}^+,z_{2r})lk(z_{2i-1}^+,z_{2r-1})-lk(z_{2i}^+,z_{2r-1})lk(z_{2i-1}^+,z_{2r})\right) $$ is a knot invariant.
In Lemma~\ref{lemalexand}, we will identify $\lambda^{\prime}$ with $\frac12\Delta_K^{\prime \prime}(1)$ where $\Delta_K$ denotes the Alexander polynomial of $K$ so that the surgery formula of Theorem~\ref{thmcassur} coincides with the Casson surgery formula of \cite[Thm. 1.1 (v)]{mar}, \cite[p. (xii)]{akmc} or \cite[Thm. 1.5]{gm}. Since any integral homology sphere can be obtained from $S^3$ by a finite sequence of surgeries with coefficients $\pm 1$ (see \cite[Section 4]{mar} or \cite[Lemme 2.1]{gm}, for example), it follows that $\tilde{\lambda}$ coincides with the Casson invariant for integral homology spheres.

Our proof will also yield the following theorem. Recall that the \emph{Euler class} of a nowhere zero vector field of a $3$-manifold $M$ is the Euler class of its orthogonal plane bundle in $M$.

\begin{theorem} \label{thmponesurg}
Let $F$ be a genus $g(F)$ oriented compact surface with connected boundary embedded in an oriented compact $3$-manifold $M$ whose boundary $\partial M$ is either empty or identified with $\partial B(1)$.
Let $[-2,2] \times F$ be a neighborhood of $\mathring{F} = \{0\} \times \mathring{F}$ in $M$,
and let $\ComX$ be a nowhere zero vector field of $M$ whose Euler class is a torsion element of $H^2(M;\ZZ)$, that is tangent to $[-2,2] \times \{x\}$ at any
point $(u,x)$ of $[-2,2] \times F$, and that is constant on $\partial B(1)$ when $\partial M =\partial B(1)$.
Let $K$ be a parallel of $\partial F$ inside $F$,
and let  $([-2,2] \times F)(K)$ be obtained from $[-2,2] \times F$ by $+1$--Dehn surgery along $K$. Let $\twist_K$ denote the right-handed Dehn twist about $K$.
Then $$([-2,2] \times F)(K)=[-2,0] \times F \cup_{\{0\}\times F \stackrel{\twist_K}{\leftarrow} \{0\}\times F^+} [0,2] \times F^+ $$
where $F^+$ is a copy of $F$ and $(0,x) \in \{0\}\times F^+$ is identified with $(0,\twist_K(x)) \in \{0\}\times F$.
Define the diffeomorphism
$$\begin{array}{lllll} \psi_F \colon & ([-2,2] \times F)(K) & \to &[-2,2] \times F &\\
   & (t,x) &\mapsto &(t,x)& \mbox{if}\; (t,x) \in [-2,0] \times F\\
& (t,x) &\mapsto &(t,\twist_K(x))& \mbox{if}\; (t,x) \in [0,2] \times F^+
  \end{array}$$
and let $\ComY$ be a nowhere zero vector field of $M(K)$ that coincides with $\ComX$ outside $]-1,1[ \times \mathring{F}$ and that is normal to $\psi_F^{-1}(\{t\} \times F)$ on $\psi_F^{-1}(\{t\} \times F)$ for any
$t \in [-2,2]$.
Then $$p_1(\ComY) - p_1(\ComX)=(4g(F)-1)g(F).$$
\end{theorem}

\subsection{A preliminary lemma on Pontrjagin numbers}

\begin{lemma} \label{lemponesurg}
Under the assumptions of Theorem~\ref{thmponesurg}, the variation $\left(p_1(\ComY) - p_1(\ComX)\right)$ does not depend on $M$, $K$ and $F$, it only depends
on $g(F)$. It will be denoted by $p_1(g(F))$.
\end{lemma}
\bp
Let $\tau_{F} \colon F \times \RR^2 \to TF$ be a parallelization of $F$ such that the parallelization $\ComX \oplus \tau_{F}$  of $[-2,2] \times F$ extends to a trivialization $\tau$ of $M$ --that is standard on $\partial M$ if $\partial M=S^2$--.
(Since $M$ is parallelizable and since $\pi_1(SO(3))$ is generated by a loop of rotations with arbitrary fixed axis, there exists a parallelization of $M$ that has this prescribed form on $[-2,2] \times F$.)
Observe that the degree of the tangent map to $K$ is $(1-2g)$ with respect to $\tau_{F}$.
(This degree does not depend on $\tau_{F}$ and can be computed in Figure~\ref{figsigmasur}.)
Let $K \times [-1,1]$ be a tubular neighborhood of $K$ in $F$ such that $K \times \{-1\} = \partial F$.
Then $[-1,1] \times K \times [-1,1]$ is a neighborhood $N_{\Box}(K)$ of $K$ in $M$ that has a standard parallelization $\tau_\nornu=(\ComX, TK, \nornu)$ where $TK$ stands for the unit tangent vector to $K$ and $\nornu$ is tangent to $\{(h,x)\} \times [-1,1]$.
 We assume that $$\tau_\nornu^{-1}\tau((t,k=\exp(2i\pi\theta),u),v \in \RR^3)
=((t,k,u),\rho_{(2g-1)\theta}(v))$$
where $\rho_{(2g-1)\theta}$ is the rotation whose axis is directed by the first basis vector $e_1$ of $\RR^3$ with angle $(2g-1)\theta$.

Let $\hat{K}$ be the image of $K$ (that is fixed by $\twist_K$) in $M(K)$. 
The neighborhood $N_{\Box}(\hat{K})=\psi_F^{-1}(N_{\Box}(K))$ of $\hat{K}$ in $([-2,2] \times F)(K)$ is also equipped with a standard parallelization $\hat{\tau}_\nornu=(\ComY, T\hat{K}, \hat{\nornu})=\psi_{F\ast}^{-1} \circ \tau_\nornu$. 

Define the parallelization $\tau^{\prime}$ of $M(K)$ that coincides with $\tau$ outside $]-1,1[ \times \mathring{F}$ and that is the following stabilization of the positive normal $\ComY$ to $F$ on $[-1,1] \times F$.
Let $\check{F} = F \setminus (K \times [-1,1[)$. 
On $[-1,1] \times \check{F}$, $$\tau^{\prime}(t,x,v \in \RR^3)=\tau(t,x,\rho_{(1-2g)\pi(h+1)}(v)).$$

This parallelization extends to $N_{\Box}(\hat{K})$ as a stabilization of $\ComY$ because
it extends to a square bounded by the following square meridian $\mu_K$ of $\hat{K}$
$$\mu_K=\{-1\} \times \{k\}\times [-1,1] + ([-1,1] \times (k,1)) - \{1\} \times \twist_K^{-1}(\{k\}\times [-1,1]) - ([-1,1] \times (k,-1))$$
written with respect to coordinates of $\partial N_{\Box}(K)$.

Write a (round) tubular neighborhood $N(K)$ in $N_{\Box}(K)$ as $S^1 \times D^2 =\partial D^2 \times D^2$ so that $\mu_K$ induces the same parallelization of $K$ as the longitude $(\{x\} \times D^2)$.
Let $$W_F=\left(\left([0,1] \times [-1,1] \times F\right) \cup_{\{1\}\times N(K) \sim \partial D^2 \times D^2} D^2 \times D^2 \right) \sharp (-\CC P^2)$$ be a cobordism 
from $[-1,1] \times F$ to $([-1,1] \times F)(K)$ obtained from $[0,1] \times [-1,1] \times F$ by gluing a $2$-handle $D^2 \times D^2$ along $N(K)$ using the identification of $N(K)$ with $\partial D^2 \times D^2$ above, by smoothing in a standard way, and by next performing a connected sum with a copy of $(-\CC P^2)$ in the interior of the $2$-handle.
We compute $(p_1(\tau^{\prime})-p_1(\tau))$ by using the cobordism $W_F$
completed to a signature $0$ cobordism by the product $[0,1] \times \left(M \setminus \mbox{Int}([-1,1] \times F)\right)$
where $T{[0,1]} \oplus \tau$ extends both $\tau$ and $\tau^{\prime}$. 
Since $\pi_1(SU(2))$ is trivial, the induced complex parallelization over $\partial ([0,1] \times [-1,1]) \times \check{F}$
extends as a stabilization of $T{[0,1]} \oplus \ComX$ whose restriction to $[0,1] \times [-1,1] \times \partial{\check{F}}$ only depends on the genus of $F$.
Thus $(p_1(\tau^{\prime})-p_1(\tau))$ is the obstruction to extending this extension to $\left([0,1] \times N_{\Box}(K)\cup_{\{1\}\times N(K) \sim \partial D^2 \times D^2} D^2 \times D^2 \right) \sharp (-\CC P^2)$ and it only depends on $g(F)$. Call it $p_1(g(F))$.

Now compose $\tau$ and $\tau^{\prime}$ by a small rotation whose axis is the second basis vector $e_2$ of $\RR^3$ around $[-2,2] \times F$, so that  $\ComX\neq \pm\tau(e_1)$ on $[-1,1] \times F$, and $\ComX$ and $\tau(e_1)$ are transverse.
Then $L_{\ComX=\tau(e_1)}=L_{\ComY=\tau^{\prime}(e_1)}$, $L_{\ComX=-\tau(e_1)}=L_{\ComY=-\tau^{\prime}(e_1)}$.
Furthermore, since $L_{\ComX=\tau(e_1)}$ does not meet $[-1,1] \times F$, and since it is rationally null-homologous --because the Euler class of $\ComX$ is a torsion element of $H^2(M;\ZZ)$ (see \cite[Theorem 1.1]{lescomb} for details)--, $L_{\ComX=\tau(e_1)}$ bounds a Seifert surface disjoint from $N_{\Box}(K)$ and $L_{\ComY=\tau^{\prime}(e_1)}$ bounds the same Seifert surface in $M(K) \setminus N_{\Box}(\hat{K})$
so that $$lk(L_{\ComX=\tau(e_1)},L_{\ComX=-\tau(e_1)})=lk(L_{\ComY=\tau^{\prime}(e_1)},L_{\ComY=-\tau^{\prime}(e_1)})$$
and
$$p_1(\ComX)-p_1(\tau)=p_1(\ComY)-p_1(\tau^{\prime})$$
according to Theorem~\ref{thmlescomb}, if $H_1(M;\QQ)=0$, and according to \cite[Theorem 1.2]{lescomb}, more generally.
\eop

\subsection{Introduction to the proof of the surgery formula}

Let us now begin our proof of Theorem~\ref{thmcassur} by fixing the Heegaard diagrams that we are going to use.

\begin{figure}[h]
\begin{center}
\begin{tikzpicture} 
\useasboundingbox (-7,-2) rectangle (7,2.8);
\draw [gray!50,->] (-1.4,-.5) .. controls (-1.4,1) and (-1.7,1.3) .. (-2.15,1.3) (-3.05,1.3) .. controls (-3.6,1.3) and (-3.8,1) .. (-3.8,-.5) ..  controls (-3.8,-2) and (-1.4,-2) .. (-1.4,-.5);
\draw [gray!50,->] (-2.6,-.3) .. controls (-2.6,2.5) and (-2.6,2.5) .. (-3.8,2.5) .. controls (-4.8,2.5) and (-5,2) .. (-5,-.3) ..  controls (-5,-1) and (-2.6,-1) .. (-2.6,-.3);
\draw [gray!50] (-2.6,-.5) node[right]{\scriptsize $z_2$} (-1.4,-.5) node[right]{\scriptsize $z_1$};
\draw [thick,->] (-2.15,1) .. controls (-1.8,1) and (-1.7,0) .. (-2,0) .. controls (-2.3,0) .. (-2.3,1.3) .. controls  (-2.3,2.8) .. (-3.8,2.8) .. controls (-5,2.8) and (-5.3,2.4) .. (-5.3,0) .. controls (-5.3,-2) .. (0,-2);
\draw [thick,->] (-3.05,1.6) .. controls (-4.1,1.6) and (-4,0) .. (-4.4,0) .. controls (-4.7,0) .. (-4.7,1.3) .. controls (-4.7,1.7) and (-4.3,2.2) .. (-3.8,2.2) .. controls (-2.9,2.2) .. (-2.9,1.3) .. controls (-2.9,0) .. (-3.2,0) .. controls (-3.5,0) and (-3.5,1) .. (-3.05,1);
\draw [thick,->] (-.5,0) .. controls (-1.5,0) and (-1.1,1.6) .. (-2.15,1.6);
\draw [->] (-1.4,-2) -- (-1.4,-1);
\draw (-1.4,-1) .. controls (-1.4,1) and (-1.7,1.3) .. (-2.15,1.3) (-3.05,1.3) .. controls (-3.6,1.3) and (-3.8,1) .. (-3.8,-.5) -- (-3.8,-1.9);
\draw [->] (-2.6,-2) -- (-2.6,-1);
\draw (-2.6,-1) .. controls (-2.6,2.5) and (-2.6,2.5) .. (-3.8,2.5) .. controls (-4.8,2.5) and (-5,2) .. (-5,-1.5);
\draw (-2.6,-1) node[right]{\scriptsize $u_2$} (-1.4,-1) node[right]{\scriptsize $u_1$} (0,-2) node[above]{\scriptsize $K$} (0,0) node{\dots};
\begin{scope}[xshift=6.4cm]
\draw [gray!50,->] (-1.4,-.5) .. controls (-1.4,1) and (-1.7,1.3) .. (-2.15,1.3) (-3.05,1.3) .. controls (-3.6,1.3) and (-3.8,1) .. (-3.8,-.5) ..  controls (-3.8,-2) and (-1.4,-2) .. (-1.4,-.5);
\draw [gray!50,->] (-2.6,-.3) .. controls (-2.6,2.5) and (-2.6,2.5) .. (-3.8,2.5) .. controls (-4.8,2.5) and (-5,2) .. (-5,-.3) ..  controls (-5,-1) and (-2.6,-1) .. (-2.6,-.3);
\draw [gray!50] (-2.5,-.3) node[left]{\scriptsize $z_{2g(\Sigma)}$} (-1.3,-.65) node[left]{\scriptsize $z_{2g(\Sigma)-1}$};
\draw [->] (-1.4,-1.5) -- (-1.4,-1);
\draw [thick,->] (-2.15,1) .. controls (-1.8,1) and (-1.7,0) .. (-2,0) .. controls (-2.3,0) .. (-2.3,1.3) .. controls  (-2.3,2.8) .. (-3.8,2.8) .. controls (-5,2.8) and (-5.3,2.4) .. (-5.3,1) .. controls (-5.3,0) .. (-5.9,0);
\draw [thick,->] (-3.05,1.6) .. controls (-4.1,1.6) and (-4,0) .. (-4.4,0) .. controls (-4.7,0) .. (-4.7,1.3) .. controls (-4.7,1.7) and (-4.3,2.2) .. (-3.8,2.2) .. controls (-2.9,2.2) .. (-2.9,1.3) .. controls (-2.9,0) .. (-3.2,0) .. controls (-3.5,0) and (-3.5,1) .. (-3.05,1);
\draw [thick,->] (-6.4,-2) .. controls (-1.1,-2) .. (-1.1,0) .. controls (-1.1,1.2) and (-1.6,1.6) .. (-2.15,1.6);
\draw (-1.4,-1) .. controls (-1.4,1) and (-1.7,1.3) .. (-2.15,1.3) (-3.05,1.3) .. controls (-3.6,1.3) and (-3.8,1) .. (-3.8,-.5) -- (-3.8,-2);
\draw [->] (-2.6,-1.9) -- (-2.6,-1);
\draw (-2.6,-1) .. controls (-2.6,2.5) and (-2.6,2.5) .. (-3.8,2.5) .. controls (-4.8,2.5) and (-5,2) .. (-5,-2);
\draw (-2.5,-1) node[left]{\scriptsize $u_{2g(\Sigma)}$} (-1.25,-.9) node[left]{\scriptsize $u_{2g(\Sigma)-1}$};
\end{scope}
\end{tikzpicture}
\caption{The curves $u_i$ on the surface $\Sigma$}
\label{figsigmasurui}
\end{center}
\end{figure}
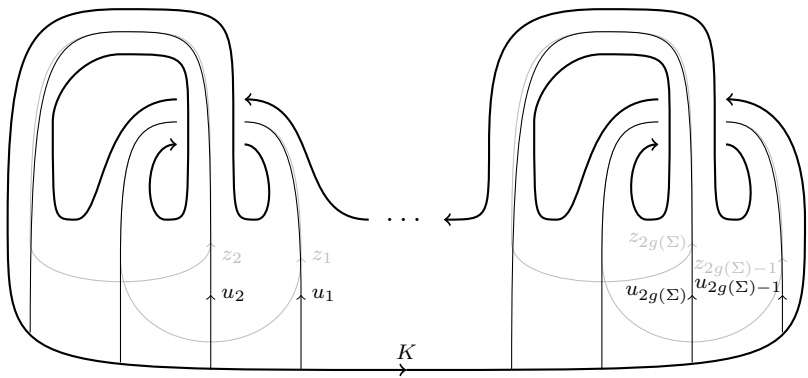

Let $u_i$ be non-intersecting curves of $\Sigma$ as in Figure~\ref{figsigmasurui} with boundaries in $\partial \Sigma$ such that $u_i$ is homologous to $z_i$ in $H_1(\Sigma,\partial \Sigma)$.
Then the $u_i \times [-1,1]$ form a system of (topological) meridian disks for the handlebody
$\handA$. Set $\alpha_i=-\partial(u_i \times [-1,1])$.
Fix a system of meridians $(\beta_j)_{j \in \underline{g}}$ that meet the $\alpha$ curves transversally and that meet $K\times [-1,1]$ as a product by $[-1,1]$.
Set $\Sigma^+=\Sigma \times \{1\}$ and $\Sigma^-=\Sigma \times \{-1\}$.
Assume that the Heegaard diagram $\CD=((\alpha_i)_{i \in \underline{g}},(\beta_j)_{j \in \underline{g}})$ has a matching $\matchingo=\{\matchm_i\}_{i \in \underline{g}}$ where $\matchm_i \in \alpha_i \cap \beta_i$ and $\matchm_i \in \Sigma^-$ (up to isotopies of the curves $\beta$). The invariant $\tilde{\lambda}(M)$ will be computed with 
the diagram $\CD$,
and the invariant $\tilde{\lambda}(M(K))$ will be computed with 
the diagram $$\CD^{\prime}=((\alpha_i)_{i \in \underline{g}},(\beta^{\prime}_j=\twist_K(\beta_j))_{j \in \underline{g}}).$$
We fix a common exterior point $\pointw$ for $\CD$ and $\CD^{\prime}$ in $\Sigma^-$.

\begin{lemma}
\label{lemvarsurpone}
 The variation $\left(p_1(\ComX(\CD^{\prime},\pointw,\matchingo)) - p_1(\ComX(\CD,\pointw,\matchingo))\right)$ is equal to the number $p_1(g(\Sigma))$ defined in Lemma~\ref{lemponesurg}.
\end{lemma}
\bp Apply Lemma~\ref{lemponesurg} to $$F=\Sigma^+ \cup_{K \times \{1\}} (K \times [-1,1]) \subset \partial \handA, $$ 
$\ComX=\ComX(\CD,\pointw,\matchingo)$ and $\ComY=\ComX(\CD^{\prime},\pointw,\matchingo)$.
\eop

Let $u_i$ also denote $u_i \times \{1\}=\alpha_i \cap \Sigma^+$.

\begin{figure}[h]
\begin{center}
\begin{tikzpicture} 
\useasboundingbox (-1.5,-1) rectangle (10.4,3);
\draw (.9,0) -- (8,0) arc (-90:90:1.5) (8,3) -- (.9,3) arc (90:270:1.5) (.9,1) -- (8,1) arc (-90:90:.5) (8,2) -- (.9,2) arc (90:270:.5);
\draw [very thick,->] (.9,2.5) arc (90:270:1) (.9,.5) -- (4,.5);
\draw [very thick,->] (4,.5) -- (8,.5) arc (-90:90:1) (8,2.5) -- (.9,2.5);
\draw (.9,2.5) node[above]{\scriptsize $K$} (6.4,.5) node[above]{\scriptsize \dots};
\draw [->] (7.7,.75) -- (7.7,1) (7.7,0) -- (7.7,.75) node[right]{\scriptsize $\alpha_{g-1}$};
\draw [->] (7.5,.75) -- (7.5,1) (7.5,0) -- (7.5,.75);
\draw [->] (7.3,.75) -- (7.3,0) (7.3,1) -- (7.3,.75);
\draw [->] (7.1,.75) -- (7.1,0) (7.1,1) -- (7.1,.75) node[left]{\scriptsize $\alpha_g$};
\draw [dashed,->] (7.7,1.4) arc (0:180:.2) (7.5,1.4) arc (0:180:.2) (7.7,1.3) -- (7.7,1.4) (7.7,1) -- (7.7,1.3);
\draw [dashed,->] (7.5,1) -- (7.5,1.4) (7.3,1) -- (7.3,1.4) (7.1,1.2) -- (7.1,1) (7.1,1.4) -- (7.1,1.2);
\draw (7.1,1.3) node[left]{\scriptsize $u_g$} (7.65,1.3) node[right]{\scriptsize $u_{g-1}$};
\begin{scope}[xshift=-2.1cm]
\draw [->] (7.7,.75) -- (7.7,1) (7.7,0) -- (7.7,.75) node[right]{\scriptsize $\alpha_1$};
\draw [->] (7.5,.75) -- (7.5,1) (7.5,0) -- (7.5,.75);
\draw [->] (7.3,.75) -- (7.3,0) (7.3,1) -- (7.3,.75);
\draw [->] (7.1,.75) -- (7.1,0) (7.1,1) -- (7.1,.75) node[left]{\scriptsize $\alpha_2$};
\draw [dashed,->] (7.7,1.4) arc (0:180:.2) (7.5,1.4) arc (0:180:.2) (7.7,1.3) -- (7.7,1.4) (7.7,1) -- (7.7,1.3);
\draw [dashed,->] (7.5,1) -- (7.5,1.4) (7.3,1) -- (7.3,1.4) (7.1,1.2) -- (7.1,1) (7.1,1.4) -- (7.1,1.2);
\draw (7.1,1.3) node[left]{\scriptsize $u_2$} (7.7,1.3) node[right]{\scriptsize $u_1$};
\end{scope}
\begin{scope}[xshift=.9cm]
\draw [->] (.3,.75) -- (.3,0) (.3,1) -- (.3,.75);
\draw [->] (.5,.75) -- (.5,1) (.5,0) -- (.5,.75);
\draw [->] (.9,.5) node[above]{\scriptsize \dots} (1.1,.75) -- (1.1,1) (1.1,0) -- (1.1,.75);
\draw [->] (1.7,.5) node[above]{\scriptsize \dots} (2.1,.75) -- (2.1,0) (2.1,1) -- (2.1,.75);
\draw [dashed]  (.5,1) arc (180:0:.8) (1.1,1) arc (0:180:.4);
\draw (.3,1.2) node[above]{\scriptsize $q_+$} (2.1,1.4) node[above]{\scriptsize $q_{-}$};
\draw (.3,-.15) -- (.4,-.25) -- (2,-.25) -- (2.1,-.15) (1.2,-.2) node[below]{\scriptsize Intersection of $K$ with the $\beta$ curves};
\end{scope}
\end{tikzpicture}
\caption{The intersections of $K$ with the curves of $\CD$}
\label{figintalongk}
\end{center}
\end{figure}
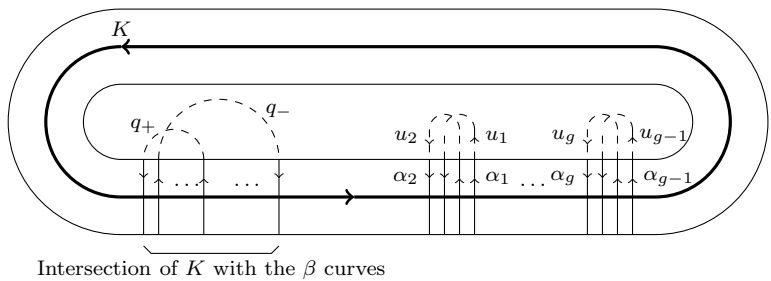

Assume that along $K$, from some basepoint of $K$, we first meet all the intersection points of $K$ with the $\beta_j$ and next the intersection points of $K$ with the $\alpha_i$, that correspond to the endpoints of the $u_i$, as in Figure~\ref{figintalongk}.

Recall $\lambda^{\prime}= \sum_{(i,r) \in \underline{g(\Sigma)}^2}\left(lk(z_{2i}^+,z_{2r})lk(z_{2i-1}^+,z_{2r-1})-lk(z_{2i}^+,z_{2r-1})lk(z_{2i-1}^+,z_{2r})\right).$

We are going to prove the following lemmas.

\begin{lemma}
\label{lemvarsurlamb}
$$\ell_2(\CD^{\prime})-\ell_2(\CD)=8\lambda^{\prime}.$$
\end{lemma}

\begin{lemma}
\label{lemvarsurlk}
$$s_{\ell}(\CD^{\prime},\matchingo)-s_{\ell}(\CD,\matchingo)=-g(\Sigma)^2 - 2 \lambda^{\prime}.$$
\end{lemma}

\begin{lemma}
\label{lemvarsureuler}
$$e(\CD^{\prime},\pointw,\matchingo)-e(\CD,\pointw,\matchingo)=(1-g)g(\Sigma).$$
\end{lemma}

It follows from these lemmas that 
$$24\tilde{\lambda}(M(K))-24\tilde{\lambda}(M)=24\lambda^{\prime} + 4 g(\Sigma)(g(\Sigma)-1) - p_1(g(\Sigma))$$

Applying this formula to a trivial knot $U$ seen as the boundary of a genus $g(\Sigma)$ surface $\Sigma_U$ for which $\lambda^{\prime}=0$ shows that 
$$p_1(g(\Sigma))=4 g(\Sigma)(g(\Sigma)-1)$$
since $M(U)$ is diffeomorphic to $M$.

Thus Lemmas~\ref{lemvarsurlamb}, \ref{lemvarsurlk} and \ref{lemvarsureuler} imply Theorems~\ref{thmcassur} and \ref{thmponesurg}, and we are left with their proofs that occupy most of the end of this section.

\subsection{Preliminaries for the proofs of the remaining three lemmas}

Set $\overline{2r}=2r-1$ and $\overline{2r-1}=2r$.

\begin{lemma}
\label{lemlkfundasur}
For any $(i,r) \in \underline{g}^2$,
 $$\sum_{j=1}^g\CJ_{ji} \langle u_r, \beta_j \rangle=\langle z_i,z_{\overline{i}} \rangle lk(z^+_r,z_{\overline{i}}).$$
\end{lemma}
\bp
Think of $\handA$ as a thickening of a wedge of the $z_i$. Let $m(z_i)$ denote a meridian of $z_i$ on $\partial \handA$.
Then $z_r^+=\sum_{k=1}^g lk(z^+_r,z_k)m(z_k)$ in $H_1(\handB;\QQ)$.
Since $m(z_k)$ is homologous to $\langle z_{\overline{k}},z_k \rangle (z^+_{\overline{k}}  - z^-_{\overline{k}})$ in $\partial \handA$,
$$\begin{array}{ll}
   \langle m(z_k),\beta_j \rangle&=\langle z_{\overline{k}},z_k \rangle \langle u_{\overline{k}}  - u^-_{\overline{k}},\beta_j \rangle
\\&=\langle z_{\overline{k}},z_k \rangle \langle \alpha_{\overline{k}},\beta_j \rangle
  \end{array}$$
$\langle u_r,\beta_j \rangle=\langle z_r^+,\beta_j \rangle =\sum_{k=1}^g lk(z^+_r,z_k)\langle z_{\overline{k}},z_k \rangle \langle \alpha_{\overline{k}},\beta_j \rangle,$
and 
$\sum_{j=1}^g\CJ_{ji}\langle u_r,\beta_j \rangle=\langle z_i,z_{\overline{i}} \rangle lk(z^+_r,z_{\overline{i}})$.
\eop

\begin{lemma}
\label{lemlkfundasurgg}
 $$\sum_{(i,j) \in \underline{g}^2}\CJ_{ji} \langle u_i, \beta_j \rangle=g(\Sigma).$$
\end{lemma}
\bp $\sum_{i=1}^g\langle z_i,z_{\overline{i}} \rangle lk(z^+_i,z_{\overline{i}})=\sum_{r=1}^{g(\Sigma)}(lk(z^+_{2r-1},z_{2r})-lk(z^+_{2r},z_{2r-1}))$.
\eop 

\begin{figure}[h]
\begin{center}
\begin{tikzpicture} 
\useasboundingbox (-1.5,0) rectangle (10.4,3);
\draw (.9,0) -- (8,0) arc (-90:90:1.5) (8,3) -- (.9,3) arc (90:270:1.5) (.9,1) -- (8,1) arc (-90:90:.5) (8,2) -- (.9,2) arc (90:270:.5);
\draw (6.4,.5) node[above]{\scriptsize \dots};
\draw [->] (7.7,.75) -- (7.7,1) (7.7,0) -- (7.7,.75) node[right]{\scriptsize $\alpha_{g-1}$};
\draw [->] (7.5,.75) -- (7.5,1) (7.5,0) -- (7.5,.75);
\draw [->] (7.3,.75) -- (7.3,0) (7.3,1) -- (7.3,.75);
\draw [->] (7.1,.75) -- (7.1,0) (7.1,1) -- (7.1,.75) node[left]{\scriptsize $\alpha_g$};
\draw [dashed,->] (7.7,1.4) arc (0:180:.2) (7.5,1.4) arc (0:180:.2) (7.7,1.3) -- (7.7,1.4) (7.7,1) -- (7.7,1.3);
\draw [dashed,->] (7.5,1) -- (7.5,1.4) (7.3,1) -- (7.3,1.4) (7.1,1.2) -- (7.1,1) (7.1,1.4) -- (7.1,1.2);
\draw (7.1,1.3) node[left]{\scriptsize $u_g$} (7.65,1.3) node[right]{\scriptsize $u_{g-1}$};
\begin{scope}[xshift=-2.1cm]
\draw [->] (7.7,.75) -- (7.7,1) (7.7,0) -- (7.7,.75) node[right]{\scriptsize $\alpha_1$};
\draw [->] (7.5,.75) -- (7.5,1) (7.5,0) -- (7.5,.75);
\draw [->] (7.3,.75) -- (7.3,0) (7.3,1) -- (7.3,.75);
\draw [->] (7.1,.75) -- (7.1,0) (7.1,1) -- (7.1,.75) node[left]{\scriptsize $\alpha_2$};
\draw [dashed,->] (7.7,1.4) arc (0:180:.2) (7.5,1.4) arc (0:180:.2) (7.7,1.3) -- (7.7,1.4) (7.7,1) -- (7.7,1.3);
\draw [dashed,->] (7.5,1) -- (7.5,1.4) (7.3,1) -- (7.3,1.4) (7.1,1.2) -- (7.1,1) (7.1,1.4) -- (7.1,1.2);
\draw (7.1,1.3) node[left]{\scriptsize $u_2$} (7.7,1.3) node[right]{\scriptsize $u_1$};
\end{scope}
\begin{scope}[xshift=.9cm]
\draw [dashed]  (.5,1) arc (180:0:.8) (1.1,1) arc (0:180:.4);
\draw (.1,1.2) node[above]{\scriptsize $\twist_K(q_+)$} (2.3,1.4) node[above]{\scriptsize $\twist_K(q_{-})$};
\end{scope}
\draw [->] (3,1) .. controls (3,.7) and (2.4,.52) .. (2,.52) -- (.9,.52) arc (270:90:1.1) -- (.9,2.72) -- (8,2.72) arc (90:-90:1.3) -- (3.4,.12) ..  controls (3,.12) .. (3,0);
\draw [->] (2,0) .. controls (2,.24) and (3,.24) .. (3.4,.24) -- (8,.24) arc (-90:90:1.16) (8,2.58) -- (.9,2.58) arc (90:270:.97) (.9,.64) -- (1.4,.64) .. controls (1.7,.64) and (2,.7) .. (2,1);
\draw [->] (1.4,0) .. controls (1.4,.28) and (3,.36) .. (3.4,.36) -- (8,.36) arc (-90:90:1.04) (8,2.44) -- (.9,2.44) arc (90:270:.84) (.9,.76) -- (1.2,.76) .. controls (1.4,.76) .. (1.4,1);
\draw [->] (1.2,1) .. controls (1.2,.88) .. (1.05,.88) -- (.9,.88) arc (270:90:.72) -- (.9,2.32) -- (8,2.32) arc (90:-90:.92) -- (3.4,.48) ..  controls (3,.48) and (1.2,.3) .. (1.2,0);
\end{tikzpicture}
\caption{The diagram $\CD^{\prime}$ in a neighborhood of $K$ on $\partial \handA$}
\label{figintalongktwist}
\end{center}
\end{figure}
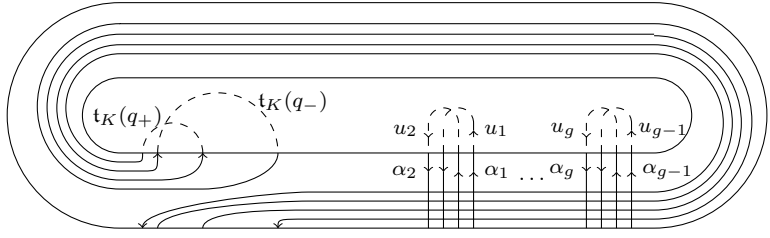

For $j \in \underline{g}$, let $Q_j$ denote the set of connected components of $\beta^{\prime}_j \cap (\Sigma^+ \cup (\partial \Sigma \times [-1,1])$. Let $Q=\cup_{j=1}^g Q_j$. For an arc $q$ of $Q_j$, set $j(q)=j$. The intersection of an arc $q$ of $Q$ with $\Sigma^+ \times \{1\}$ will be denoted by $q^+$.
Let $\CaC$ and $\CaC^{\prime}$ denote the set of crossings of $\CD$ and $\CD^{\prime}$, respectively.

For each $(q,i) \in Q \times \underline{g}$, there is a set $\CaC(q,i)=\alpha_i \cap q$ of $4$ crossings. 
Then $$\CaC^{\prime} = \CaC \coprod \coprod_{(q,i) \in Q \times \underline{g}}\CaC(q,i).$$

\begin{figure}[h]
\begin{center}
\begin{tikzpicture} 
\useasboundingbox (-5,-1) rectangle (5,1.9);
\begin{scope}[xshift=.7cm]
\draw [->] (3,.9) -- (3,1.4) (3,0) -- (3,.9); 
\draw [dashed,->] (2,1.4) arc (180:135:.5) (3,1.4) arc (0:135:.5);
\draw [->] (2,.9) -- (2,0) (2,1.4) -- (2,.9);
\draw (1.95,1.85) node{\scriptsize $u_i$};
\draw [->] (2.5,1.2) -- (4,1.2) (0,0) .. controls (0,1.2) .. (2,1.2) -- (2.5,1.2);
\draw [dashed,->] (4,.6) arc (-90:0:.3) (4,1.2) arc (90:0:.3);
\draw (4.3,.9) node[right]{\scriptsize $q$} (2.5,.7) node[below]{\scriptsize $h(q,i)$} (2.5,1.1) node[above]{\scriptsize $t(q,i)$} (3.5,.4) node{\scriptsize $d_1(q,i)$} (3.5,1.4) node{\scriptsize $d_2(q,i)$} (1.8,.4) node{\scriptsize $d_4$} (1.8,1.4) node{\scriptsize $d_3$} (2.15,-.5) node{\scriptsize $\CaC(q,i)$ when $\sigma(q,i)=-1$};
\draw [->] (4,.6) -- (2.5,.6);
\draw [->] (2.5,.6) -- (2,.6) .. controls (1,.6) .. (1,0);
\fill (2,.6) circle (0.05) (3,.6) circle (0.05) (2,1.2) circle (0.05) (3,1.2) circle (0.05);
\end{scope}
\begin{scope}[xshift=-5cm]
\useasboundingbox (-5,-1) rectangle (5,1.9);
\draw [->] (3,.9) -- (3,1.4) (3,0) -- (3,.9); 
\draw [dashed,->] (2,1.4) arc (180:135:.5) (3,1.4) arc (0:135:.5);
\draw [->] (2,.9) -- (2,0) (2,1.4) -- (2,.9);
\draw (1.95,1.85) node{\scriptsize $u_i$};
\draw [->] (4,1.2) -- (2.5,1.2);
\draw [<-] (0,0) .. controls (0,1.2) .. (2,1.2) -- (2.5,1.2);
\draw [dashed,->] (4,1.2) arc (90:0:.3) (4,.6) arc (-90:0:.3);
\draw (4.3,.9) node[right]{\scriptsize $q$} (2.5,.7) node[below]{\scriptsize $t(q,i)$} (2.5,1.1) node[above]{\scriptsize $h(q,i)$} (3.5,.4) node{\scriptsize $d_1(q,i)$} (3.5,1.4) node{\scriptsize $d_2(q,i)$} (1.8,.4) node{\scriptsize $d_4$} (1.8,1.4) node{\scriptsize $d_3$} (2.15,-.5) node{\scriptsize $\CaC(q,i)$ when $\sigma(q,i)=1$};
\draw (4,.6) -- (2.5,.6);
\draw [<-] (2.5,.6) -- (2,.6) .. controls (1,.6) .. (1,0);
\fill (2,.6) circle (0.05) (3,.6) circle (0.05) (2,1.2) circle (0.05) (3,1.2) circle (0.05);
\end{scope}
\end{tikzpicture}
\caption{The new crossings of $\CD^{\prime}$}
\label{fignewcrosssur}
\end{center}
\end{figure}

Denote $\CaC(q,i)=\{d_1(q,i),d_2(q,i),d_3(q,i),d_4(q,i)\}$ where following $\alpha_i$ from $\matchm_i$, $d_1(q,i)$, $d_2(q,i)$, $d_3(q,i)$ and $d_4(q,i)$ are met in this order. Set $\sigma(q,i)=\sigma(d_2(q,i))$. Then
$$\sigma(q,i)=\sigma(d_2(q,i))=\sigma(d_4(q,i))=-\sigma(d_1(q,i))=-\sigma(d_3(q,i)).$$
Let $t(q,i)$ denote the (tail) arc of $q$ before $q^+$ with its ends in $\CaC(q,i)$ and
let $h(q,i)$ denote the (head) arc of $q$ after $q^+$ with its ends in $\CaC(q,i)$.

If $\sigma(q,i)=-1$, then $q$ goes from left to right as $q_-$ in Figures~\ref{figintalongk} and \ref{figintalongktwist},
following $q$ we meet $\CaC(q,i)$ in the order $d_3$, $d_2$, $d_1$, $d_4$,
$t(q,i)=\halfbra d_3,d_2 \halfbra_{\beta}$ and $h(q,i)=\halfbra d_1,d_4 \halfbra_{\beta}$.

If $\sigma(q,i)=1$, then $q$ goes from right to left as $q_+$ in Figures~\ref{figintalongk} and \ref{figintalongktwist},
following $q$ we meet $\CaC(q,i)$ in the reversed order $d_4$, $d_1$, $d_2$, $d_3$,
$t(q,i)=\halfbra d_4,d_1 \halfbra_{\beta}$ and $h(q,i)=\halfbra d_2,d_3 \halfbra_{\beta}$.
Thus $t(q,i)$ begins at $d_{b(q,i)}(q,i)$ where $b(q,i)=3$ if $\sigma(q,i)=-1$ and $b(q,i)=4$ if $\sigma(q,i)=1$.

Note that for any $(i,j) \in \underline{g}^2$, $\langle \alpha_i, \beta_j \rangle=\langle \alpha_i, \beta^{\prime}_j \rangle$ so that the coefficients $\CJ_{ji}$ are the same for $\CD$ and $\CD^{\prime}$.

The set of crossings of $\CD$ on $\Sigma^+$ (resp. on $\Sigma^-$)  will be denoted by $\CaC^+$ (resp. by $\CaC^-$). 

\noindent{\sc Proof of Lemma~\ref{lemvarsureuler}:} On the rectangle $R_{\CD}$ of Figure~\ref{figplansplit} for $(\CD,\pointw,\matchingo)$, let $p^{\prime}_i$ (resp. $p^{\prime \prime}_i$) denote the other end of the diameter of $\alpha^{\prime}_i$ (resp. $\alpha^{\prime \prime}_i$) that contains the crossing $\matchm_i$ of $\matchingo$.
Draw the knot $K$ on a picture of the Heegaard diagram as in Figure~\ref{figplansplit} so that $K$ meets the curves
$\alpha^{\prime}$ and $\alpha^{\prime \prime}$ as the $\beta_j$ do, away from the points of $\matchingo$,
with horizontal tangent vectors near the $p^{\prime}_i$ and the $p^{\prime \prime}_i$.
Let $N(\matchingo)$ denote an open tubular neighborhood of $\matchingo$ in $\partial \handA$ made of $2g(\Sigma)$ open disks.
See $\partial \handA \setminus \mathring{N}(\matchingo)$ as obtained from the rectangle $R_{\CD}$ with holes bounded by the $\alpha^{\prime}_i$  and the $\alpha^{\prime \prime}_i$, by gluing horizontal thin
rectangles $D_i$ along their two vertical small sides that are neighborhoods of $p^{\prime}_i$ or $p^{\prime \prime}_i$ in $\alpha^{\prime}_i$ or $\alpha^{\prime \prime}_i$.
The standard parallelization of this picture equips $\partial \handA \setminus \mathring{N}(\matchingo)$ with a parallelization so that the degree $\tilth(K)$ of the tangent to $K$   is $1-2g(\Sigma)$ in this figure.
A similar picture for $(\CD^{\prime},\pointw,\matchingo)$ is obtained by performing the Dehn twist about $K$ on the $\beta$-curves in this figure. Since these curves do not intersect $K$ algebraically,
the $\tilth(\beta_j)$ are unchanged by this operation. Similarly, for any crossing $c$ of $\CaC^-$,
$\tilth(\halfbra \matchm_{j(c)},c\halfbra_{\beta})$ is unchanged and so is $\tilth(c)$.
For any crossing $c$ of $\CaC^+$, we have $\tilth^{\prime}(c)=\tilth(c)+1-2g(\Sigma)$ since 
$\langle K, \halfbra \matchm_{j(c)},c\halfbra_{\beta}\rangle =1$.
Now, let 
$(q,i) \in Q \times \underline{g}$.
The contribution of $\CaC(q,i)$ to $\left(e(\CD^{\prime},\pointw,\matchingo)-e(\CD,\pointw,\matchingo)\right)$ is 
$$\pm\CJ_{j(q)i}\left(\tilth(h(q,i))+\tilth(t(q,i)) - \sum_{(r,s) \in \underline{g}^2}\CJ_{sr}\langle \alpha_r,h(q,i)+t(q,i)\rangle\tilth(\beta_{s})\right)$$
that is zero.
Finally,
$$\begin{array}{ll}e(\CD^{\prime},\pointw,\matchingo)-e(\CD,\pointw,\matchingo)&=(1-g)\sum_{c \in \CaC^+}\CJ_{j(c)i(c)} \sigma(c)\\
&=(1-g)\sum_{(i,j) \in \underline{g}^2}\CJ_{ji} \langle u_i,\beta_j\rangle\\
&=(1-g)g(\Sigma)\end{array}$$
according to Lemma~\ref{lemlkfundasurgg}.
\eop

\subsection{Study of \texorpdfstring{$\ell$}{ell} }

Let $\ell$ and $\ell^{\prime}$ be the maps of Proposition~\ref{propevalbis} associated with $\CD$ and $\CD^{\prime}$, respectively, with respect to the basepoints $\matchm_i$ of $\Sigma^-$.

$$\ell(c,d)=\langle [\matchm_{i(c)},c\halfbra_{\alpha}, [\matchm_{j(d)},d\halfbra_{\beta}\rangle -\sum_{(i,j) \in \underline{g}^2}\CJ_{ji} \langle [\matchm_{i(c)},c\halfbra_{\alpha},\beta_j \rangle \langle \alpha_i ,[\matchm_{j(d)},d\halfbra_{\beta} \rangle$$

\begin{lemma}
\label{lemcpcp}
Let $(c,d) \in \CaC^2$. If $(c,d) \in (\CaC^+)^2$, then
$$\ell^{\prime}(c,d)=\ell(c,d)-1$$
Otherwise,
$$\ell^{\prime}(c,d)=\ell(c,d)$$
\end{lemma}
\bp Recall that the $\matchm_i$ are in $\Sigma^-$.
Note that $\twist_K([\matchm_{j(d)},d|_{\beta})$ is obtained from $[\matchm_{j(d)},d|_{\beta}$ by adding some multiple of $K$ located in $K \times [-1,1]$, algebraically, so that 
$$\langle \alpha_i ,\twist_K([\matchm_{j(d)},d|_{\beta}) \rangle=\langle \alpha_i ,[\matchm_{j(d)},d\halfbra_{\beta} \rangle$$ for any $i \in \underline{g}$.
Since $\twist_K(\beta_j)$ differs from $\beta_j$ by an algebraically null sum of copies of $K$ in $K \times [-1,1]$,
$$\langle [\matchm_{i(c)},c\halfbra_{\alpha},\beta^{\prime}_j\rangle =\langle [\matchm_{i(c)},c\halfbra_{\alpha},\beta_j\rangle $$ for any $j \in \underline{g}$.
Thus in any case, 
$$\ell^{\prime}(c,d)-\ell(c,d)= \langle [\matchm_{i(c)},c|_{\alpha},\twist_K([\matchm_{j(d)},d|_{\beta}) -[\matchm_{j(d)},d|_{\beta} \rangle. $$
If $d \in \Sigma^-$, $\twist_K([\matchm_{j(d)},d|_{\beta})$ differs from $[\matchm_{j(d)},d|_{\beta}$ by an algebraically null sum of copies of $K$ in $K \times [-1,1]$ so that $\ell^{\prime}(c,d)=\ell(c,d)$.
If $d \in \Sigma^+$,
$\ell^{\prime}(c,d)-\ell(c,d)= \langle [\matchm_{i(c)},c|_{\alpha},K \rangle$.
If $c \in \Sigma^-$, the arc $[\matchm_{i(c)},c|_{\alpha}$ meets $K \times [-1,1]$ as the empty set or
as two parallel arcs with opposite direction and $\ell^{\prime}(c,d)=\ell(c,d)$.
If $c \in \Sigma^+$, then the arc $[\matchm_{i(c)},c|_{\alpha}$ meets $K \times [-1,1]$ as an arc that crosses $K$ once with a negative sign.
\eop

\begin{lemma}
\label{lemcpqi}
Let $c \in \CaC$ and let $(q,i) \in Q \times \underline{g}$.
$$\sum_{d \in \CaC(q,i)}\sigma(d)\ell^{\prime}(c,d)=\sum_{d \in \CaC(q,i)}\sigma(d)\ell^{\prime}(d,c)=0.$$
\end{lemma}
\bp
For any interval $I$ of a $\beta^{\prime}$-curve,
$$\sum_{d \in \CaC(q,i)}\sigma(d)\langle [\matchm_i,d\halfbra_{\alpha}, I\rangle = \sigma(q,i)\langle \halfbra d_1(q,i),d_2(q,i)\halfbra_{\alpha} + \halfbra d_3(q,i),d_4(q,i)\halfbra_{\alpha}, I\rangle$$ that is zero if $I$ has no end points in $K \times [-1,1]$.
This shows that $\sum_{d \in \CaC(q,i)}\sigma(d)\ell^{\prime}(d,c)=0$.
For any interval $I$ of an $\alpha$-curve,
$$\sum_{d \in \CaC(q,i)}\sigma(d)\langle I, [\matchm_{j(q)},d\halfbra_{\beta}, \rangle =-\langle I, t(q,i)+h(q,i)\rangle.$$
Again, this is zero if $I$ has no end points in $K \times [-1,1]$.
\eop

\begin{lemma}
\label{lemqiqi}
Let $(q,i)$ and $(q^{\prime},r)$ belong to $Q \times \underline{g}$.
If $q \neq q^{\prime}$ and $i \neq r$, then
$$\sum_{(c,d) \in \CaC(q,i) \times \CaC(q^{\prime},r)}\sigma(c)\sigma(d)\ell^{\prime}(c,d)=- lk_{K \times \{1\}}(\partial q^+,\partial q^{\prime+})lk_{K \times \{1\}}(\partial u_i,\partial u_{r}) .$$
If $q = q^{\prime}$ or $i = r$, then $\sum_{(c,d) \in \CaC(q,i) \times \CaC(q^{\prime},r)}\sigma(c)\sigma(d)\ell^{\prime}(c,d)=0$.
\end{lemma}
\bp Set $A=\sum_{(c,d) \in \CaC(q,i) \times \CaC(q^{\prime},r)}\sigma(c)\sigma(d)\ell^{\prime}(c,d)$. As in the proof of Lemma~\ref{lemcpqi},
$$\begin{array}{ll}A&=-\sum_{c \in \CaC(q,i)}\sigma(c)\langle [\matchm_i,c\halfbra_{\alpha}, t(q^{\prime},r)+h(q^{\prime},r)\rangle\\&=
-\sigma(q,i)\langle \halfbra d_1(q,i),d_2(q,i)\halfbra_{\alpha} + \halfbra d_3(q,i),d_4(q,i)\halfbra_{\alpha},t(q^{\prime},r)+h(q^{\prime},r)\rangle\end{array}$$
This is zero unless $q \neq q^{\prime}$, $i \neq r$ and  
 $lk(\partial q,\partial q^{\prime})lk(\partial u_i,\partial u_{r}) \neq 0$.
When the sign of $q^{\prime}$ changes, so does the result.
Furthermore, the result is symmetric when $(q,i)$ and $(q^{\prime},r)$ are exchanged, thanks to the symmetry of the linking number (see Proposition~\ref{propevallk}).

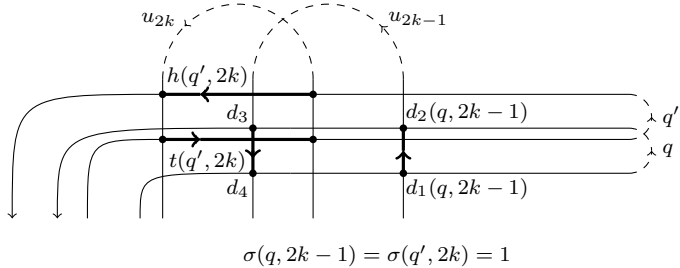
\begin{figure}[h]
\begin{center}
\begin{tikzpicture} 
\useasboundingbox (-5,-1) rectangle (5,2.4);
\begin{scope}[xshift=-5cm]
\useasboundingbox (-5,-1) rectangle (5,1.9);
\draw [-] (4,.9) -- (4,1.4) (4,-.45) -- (4,.9); 
\draw [dashed,->] (2,1.4) arc (180:135:1) (4,1.4) arc (0:135:1);
\draw [-] (2,.9) -- (2,-.45) (2,1.4) -- (2,.9);
\draw (1.95,2.2) node{\scriptsize $u_{2k}$};
\draw [->] (8.2,1.2) -- (2.5,1.2);
\draw [very thick,->] (2.5,1.2) -- (2,1.2) (2.5,.6) -- (4,.6) (2,.6) -- (2.5,.6);
\draw [very thick,<-] (2.5,1.2) -- (4,1.2);
\draw [<-] (0,-.45) .. controls (0,1.2) .. (2,1.2) -- (2.5,1.2);
\draw [dashed,->] (8.2,1.2) arc (90:0:.3) (8.2,.6) arc (-90:0:.3);
\draw (8.5,.9) node[right]{\scriptsize $q^{\prime}$} (2.6,.6) node[below]{\scriptsize $t(q^{\prime},2k)$} (2.6,1.15) node[above]{\scriptsize $h(q^{\prime},2k)$};
\draw (8.2,.6) -- (2.5,.6);
\draw [<-] (2.5,.6) -- (2,.6) .. controls (1,.6) .. (1,-.45);
\fill (2,.6) circle (0.05) (4,.6) circle (0.05) (2,1.2) circle (0.05) (4,1.2) circle (0.05);
\end{scope}
\begin{scope}[xshift=-3.8cm, yshift=-.45cm]
\useasboundingbox (-5,-1) rectangle (5,1.9);
\draw [-] (4,.9) -- (4,1.85) (4,0) -- (4,.9); 
\draw [dashed,->] (2,1.85) arc (180:45:1) (4,1.85) arc (0:45:1);
\draw [-] (2,.9) -- (2,0) (2,1.85) -- (2,.9);
\draw (4.2,2.65) node{\scriptsize $u_{2k-1}$};
\draw [-] (7,1.2) -- (2.5,1.2);
\draw [very thick,->] (2,.6) -- (2,.9) (4,.9) -- (4,1.2) (4,.6) -- (4,.9);
\draw [very thick,->] (2,1.2) -- (2,.8);
\draw [<-] (-.6,0) .. controls (-.6,1.2) .. (2,1.2) -- (2.5,1.2);
\draw [dashed,->] (7,1.2) arc (90:0:.3) (7,.6) arc (-90:0:.3);
\draw (7.3,.9) node[right]{\scriptsize $q$} (2.5,.7)
(4.85,.4) node{\scriptsize $d_1(q,2k-1)$} 
(4.85,1.4) node{\scriptsize $d_2(q,2k-1)$} (1.8,.4) node{\scriptsize $d_4$} (1.8,1.4) node{\scriptsize $d_3$} (3.65,-.5) node{\scriptsize $\sigma(q,2k-1)= \sigma(q^{\prime},2k)=1$};
\draw (7,.6) -- (2.5,.6);
\draw [-] (2.5,.6) -- (2,.6) .. controls (.5,.6) .. (.5,0);
\fill (2,.6) circle (0.05) (4,.6) circle (0.05) (2,1.2) circle (0.05) (4,1.2) circle (0.05);
\end{scope}
\end{tikzpicture}
\caption{Computation of $lk(\partial q,\partial q^{\prime})lk(\partial u_{2k-1},\partial u_{2k})$}
\label{figcomplk}
\end{center}
\end{figure}

Therefore, it suffices to prove the lemma when $\sigma(q,i)=\sigma(q^{\prime},r)=1$ and $(i,r)=(2k-1,2k)$.
When we have the order $h(q^{\prime},r)h(q,i)t(q^{\prime},r)t(q,i)$ on $K$ that coincides with $h(u_{r})h(u_i)t(u_{r})t(u_i)$, we get $A=-1$ as in Figure~\ref{figcomplk}.
For the order $h(q,i)h(q^{\prime},r)t(q,i)t(q^{\prime},r)$, we get $A=1$.
\eop

\begin{lemma} \label{lemlksur}
When $i\neq r$,
 $lk_{K \times \{1\}}(\partial u_i,\partial u_{r})=-\langle z_i,z_{r} \rangle.$\\
When $q\neq q^{\prime}$,
$$lk_{K \times \{1\}}(\partial q^+,\partial q^{\prime+})=-\sum_{k=1}^{g}\langle  z_k,z_{\overline{k}} \rangle \langle  u_k, q \rangle\langle  u_{\overline{k}}, q^{\prime}\rangle$$
and, for any $q \in Q$, $\sum_{k=1}^{g}\langle  z_k,z_{\overline{k}} \rangle \langle  u_k, q \rangle\langle  u_{\overline{k}}, q^{\prime}\rangle=0$.
\end{lemma}
\bp
Let $\gamma(q^{\prime})$ be a curve on $K \times \{1\}$ that does not meet the $\alpha$-curves, such that $\partial \gamma(q^{\prime})=\partial q^{\prime +}$.
Then $lk_{K \times \{1\}}(\partial q^+,\partial q^{\prime+})=\langle \partial q^+, \gamma(q^{\prime}) \rangle_K$.
Since $q$ and $q^{\prime}$ do not intersect, this also reads
$lk_{K \times \{1\}}(\partial q^+,\partial q^{\prime+})=-\langle q, q^{\prime +}-\gamma(q^{\prime}) \rangle_{\partial \handA}$
where  $(q^{\prime +}-\gamma(q^{\prime}))$ is a closed curve of $\Sigma^+$ whose homology class reads
$$(q^{\prime +}-\gamma(q^{\prime}))=\sum_{k=1}^{g} \langle z_k,  z_{\overline{k}}\rangle \langle q^{\prime+}-\gamma(q^{\prime}), z_{\overline{k}} \rangle_{\Sigma^+} z_k= \sum_{k=1}^{g}  \langle z_k,  z_{\overline{k}}\rangle \langle q^{\prime}, u_{\overline{k}} \rangle_{\partial \handA} z_k.$$
\eop

\begin{lemma}
\label{lembilqiqi}
$$\sum_{(q,q^{\prime}) \in Q_j \times Q_s}\sum_{(c,d) \in \CaC(q,i) \times \CaC(q^{\prime},r)}\sigma(c)\sigma(d)\ell^{\prime}(c,d)
=-\langle z_i,z_r \rangle\sum_{k=1}^{g}\langle  z_k,z_{\overline{k}} \rangle\langle  u_k,\beta_j \rangle\langle u_{\overline{k}},\beta_s \rangle.$$
\end{lemma}
\bp
According to Lemmas~\ref{lemqiqi} and \ref{lemlksur}, 
$$\sum_{(c,d) \in \CaC(q,i) \times \CaC(q^{\prime},r)}\sigma(c)\sigma(d)\ell^{\prime}(c,d)=-\langle z_i,z_r \rangle \sum_{k=1}^{g}\langle  z_k,z_{\overline{k}} \rangle \langle  u_k, q \rangle\langle  u_{\overline{k}}, q^{\prime} \rangle.$$
\eop

\begin{lemma}
\label{lemcacprime}
$$\begin{array}{lll}\sum_{c \in \CaC^{\prime} \setminus \CaC}\CJ_{j(c)i(c)}\sigma(c)\ell^{\prime}(c,c)&=&\sum_{(i,j,r,s) \in \underline{g}^4}\CJ_{ji}\CJ_{sr}\langle u_i ,\beta_s\rangle \langle u_r,\beta_j \rangle-g(\Sigma)\\&&
-\sum_{(i,j,k,s) \in \underline{g}^4}\CJ_{ji}\CJ_{s\overline{i}}\langle z_i,z_{\overline{i}} \rangle\langle  z_k,z_{\overline{k}} \rangle \langle  u_k,\beta_j \rangle\langle u_{\overline{k}}, \beta_s \rangle \end{array}$$
\end{lemma}
\bp 
Let us fix $(q,i) \in Q \times \underline{g}$ and compute $\sum_{c \in \CaC(q,i)}\sigma(c)\ell^{\prime}(c,c)$.

Since the arc $[\matchm_i,d_1(q,i)[_{\alpha}$ does not intersect the arcs $[d_{b(q,i)}(q,i),c]_{\beta}$,
$$\begin{array}{lll}\sum_{c \in \CaC(q,i)}\sigma(c)\langle [\matchm_i,c\halfbra_{\alpha}, [\matchm_{j(q)},c\halfbra_{\beta}\rangle
&=& \sum_{c \in \CaC(q,i)}\sigma(c)\langle [d_1(q,i),c\halfbra_{\alpha}, [d_{b(q,i)}(q,i),c\halfbra_{\beta}\rangle
\\
&&+\sum_{c \in \CaC(q,i)}\sigma(c)\langle [d_1(q,i),c\halfbra_{\alpha}, [\matchm_{j(q)},d_{b(q,i)}(q,i)[_{\beta}\rangle\end{array}$$
where $\sum_{c \in \CaC(q,i)}\sigma(c)\langle [d_1(q,i),c\halfbra_{\alpha}, [\matchm_{j(q)},d_{b(q,i)}(q,i)[_{\beta}\rangle$ equals $$ \sigma(q,i) \langle \halfbra d_1(q,i),d_2(q,i) \halfbra_{\alpha} + \halfbra d_3(q,i),d_4(q,i) \halfbra_{\alpha}, [\matchm_{j(q)},d_{b(q,i)}(q,i)[_{\beta}\rangle=0$$
since the arc $[\matchm_{j(q)},d_{b(q,i)}(q,i)[_{\beta}$ intersects the arcs $[d_1(q,i),c]_{\alpha}$
as whole arcs $q^{\prime}$ of $Q_{j(c)}$.

Thus
$$ \sum_{c \in \CaC(q,i)}\sigma(c)\langle [\matchm_i,c\halfbra_{\alpha}, [\matchm_{j(q)},c\halfbra_{\beta}\rangle
= 1+\sum_{c \in \CaC(q,i)}\sigma(c)\langle [d_1(q,i),c[_{\alpha}, [d_{b(q,i)}(q,i),c[_{\beta}\rangle$$
Note that neither $d_1(q,i)$ nor $d_{b(q,i)}$ contributes to the new sum.

$$\langle [d_1(q,i),d_2(q,i)[_{\alpha}, [d_{b(q,i)}(q,i),d_2(q,i)[_{\beta}\rangle =\left\{\begin{array}{ll} \sigma(d_1(q,i))=-1\;&\mbox{if}\;\sigma(q,i)=1 \\ 0 &\mbox{if}\;\sigma(q,i)=-1.\end{array} \right.$$
If $\sigma(q,i)=1$, we are left with the computation of
$$\langle [d_1(q,i),d_3(q,i)[_{\alpha},  [d_{b(q,i)}(q,i),d_3(q,i)[_{\beta}\rangle =\langle u_i,q\rangle.$$
If $\sigma(q,i)=-1$, we are left with the computation of
$$\langle [d_1(q,i),d_4(q,i)[_{\alpha},  [d_{b(q,i)}(q,i),d_4(q,i)[_{\beta}\rangle =\langle u_i,q \rangle + \sigma(d_3(q,i)).$$
In any case, $$\sum_{c \in \CaC(q,i)}\sigma(c)\langle [\matchm_i,c\halfbra_{\alpha},  [\matchm_{j(q)},c\halfbra_{\beta}\rangle=-\langle u_i,q \rangle.$$

Let us fix $(r,s) \in \underline{g}^2$ and compute
$A=\sum_{c \in \CaC(q,i)}\sigma(c)\langle[\matchm_i,c\halfbra_{\alpha} ,\beta^{\prime}_s\rangle \langle \alpha_r,[\matchm_{j(q)},c\halfbra_{\beta} \rangle$.
Observe 
$$\langle[\matchm_i,d_4(q,i)\halfbra_{\alpha} ,\beta^{\prime}_s\rangle=\langle[\matchm_i,d_1(q,i)\halfbra_{\alpha} ,\beta^{\prime}_s\rangle + \langle u_i ,\beta_s\rangle$$
and
$$\langle[\matchm_i,d_3(q,i)\halfbra_{\alpha} ,\beta^{\prime}_s\rangle=\langle[\matchm_i,d_2(q,i)\halfbra_{\alpha} ,\beta^{\prime}_s\rangle + \langle u_i ,\beta_s\rangle.$$

Let $B=\sigma(q,i)\langle u_i ,\beta_s\rangle \left(\langle \alpha_r,[\matchm_{j(q)},d_4(q,i)\halfbra_{\beta} \rangle - \langle \alpha_r,[\matchm_{j(q)},d_3(q,i)\halfbra_{\beta} \rangle \right)
=-\langle u_i ,\beta_s\rangle \langle \alpha_r,q^+ \rangle$.
$$\begin{array}{lll}A-B&=&\sigma(q,i)\langle[\matchm_i,d_1(q,i)\halfbra_{\alpha} ,\beta^{\prime}_s\rangle \left(\langle \alpha_r,[\matchm_{j(q)},d_4(q,i)\halfbra_{\beta} \rangle - \langle \alpha_r,[\matchm_{j(q)},d_1(q,i)\halfbra_{\beta} \rangle \right)\\
&&+\sigma(q,i)\langle[\matchm_i,d_2(q,i)\halfbra_{\alpha},\beta^{\prime}_s\rangle \left(\langle \alpha_r,[\matchm_{j(q)},d_2(q,i)\halfbra_{\beta} \rangle - \langle \alpha_r,[\matchm_{j(q)},d_3(q,i)\halfbra_{\beta} \rangle \right)\\
&=&\sigma(q,i)\left(\langle[\matchm_i,d_1(q,i)\halfbra_{\alpha} ,\beta^{\prime}_s\rangle - \langle[\matchm_i,d_2(q,i)\halfbra_{\alpha},\beta^{\prime}_s\rangle\right)\langle \alpha_r,h(q,i) \rangle\\
&=&-\sigma(q,i)\langle\halfbra d_1(q,i),d_2(q,i)\halfbra_{\alpha},\beta^{\prime}_s\rangle \langle \alpha_r,h(q,i) \rangle
\end{array}$$ where $\langle \alpha_r,h(q,i) \rangle =lk_{K \times \{1\}} (\partial u_r,\partial u_i)$ when $r \neq i$, so that $\langle \alpha_r,h(q,i) \rangle=\langle z_i,z_r \rangle$ in any case.
Summarizing, we get
$$\begin{array}{ll}\sum_{c \in \CaC(q,i)}\sigma(c)\ell^{\prime}(c,c)=&-\langle u_i,q \rangle
\\&+\sum_{(r,s) \in \underline{g}^2}\CJ_{sr}\left( \langle u_i ,\beta_s\rangle \langle u_r,q \rangle +\sigma(q,i)\langle\halfbra d_1(q,i),d_2(q,i)\halfbra_{\alpha},\beta_s\rangle \langle z_i,z_r \rangle\right).\end{array}$$
where $$\begin{array}{ll}\sigma(q,i)\langle\halfbra d_1(q,i),d_2(q,i)\halfbra_{\alpha},\beta_s\rangle&=-\sum_{q^{\prime} \in Q_s; q^{\prime} \neq q}lk_{K \times \{1\}} (\partial q^{\prime},\partial q)\\
&=\sum_{q^{\prime} \in Q_s; q^{\prime} \neq q}lk_{K \times \{1\}} (\partial q,\partial q^{\prime})
\\&=-\sum_{k=1}^{g}\langle  z_k,z_{\overline{k}} \rangle \langle u_k, q \rangle\langle u_{\overline{k}}, \beta_s \rangle\end{array}$$
according to Lemma~\ref{lemlksur}.

Now, let us fix $j \in \underline{g}$ and compute
$$\begin{array}{lll}\sum_{q \in Q_j}\sum_{c \in \CaC(q,i)}\sigma(c)\ell^{\prime}(c,c)&=&-\langle u_i,\beta_j\rangle+\sum_{(r,s) \in \underline{g}^2}\CJ_{sr}\langle u_i ,\beta_s\rangle \langle u_r,\beta_j \rangle\\&&
-\sum_{(r,s) \in \underline{g}^2}\CJ_{sr}\langle z_i,z_r \rangle\sum_{k=1}^{g}\langle  z_k,z_{\overline{k}} \rangle \langle  u_k,\beta_j \rangle\langle u_{\overline{k}}, \beta_s \rangle
\\&=&-\langle u_i,\beta_j\rangle+\sum_{(r,s) \in \underline{g}^2}\CJ_{sr}\langle u_i ,\beta_s\rangle \langle u_r,\beta_j \rangle\\&&
-\sum_{(k,s) \in \underline{g}^2}\CJ_{s\overline{i}}\langle z_i,z_{\overline{i}} \rangle\langle  z_k,z_{\overline{k}} \rangle \langle  u_k,\beta_j \rangle\langle u_{\overline{k}}, \beta_s \rangle
\end{array}$$
$$\begin{array}{lll}\sum_{c \in \CaC^{\prime} \setminus \CaC}\CJ_{j(c)i(c)}\sigma(c)\ell^{\prime}(c,c)&=&\sum_{(i,j) \in \underline{g}^2}\CJ_{ji}\left(\sum_{(r,s) \in \underline{g}^2}\CJ_{sr}\langle u_i ,\beta_s\rangle \langle u_r,\beta_j \rangle-\langle u_i,\beta_j\rangle\right)\\&&
-\sum_{(i,j,k,s) \in \underline{g}^4}\CJ_{ji}\CJ_{s\overline{i}}\langle z_i,z_{\overline{i}} \rangle\langle  z_k,z_{\overline{k}} \rangle \langle  u_k,\beta_j \rangle\langle u_{\overline{k}}, \beta_s \rangle.
 \end{array}$$
Conclude with Lemma~\ref{lemlkfundasurgg}.
\eop

\subsection{Proofs of the remaining two lemmas}

\begin{lemma}
\label{lemexplamp}
$$\begin{array}{ll}2\lambda^{\prime}&=\sum_{(i,r) \in \underline{g}^2}lk(z_r^+,z_i)
lk(z_{\overline{r}}^+,z_{\overline{i}})\langle z_i, z_{\overline{i}}\rangle\langle z_r , z_{\overline{r}}\rangle\\&=\sum_{(i,j,k,s) \in \underline{g}^4}\CJ_{ji}\CJ_{s\overline{i}}\langle z_i,z_{\overline{i}} \rangle\langle  z_k,z_{\overline{k}} \rangle \langle  u_k,\beta_j \rangle\langle u_{\overline{k}}, \beta_s \rangle \end{array}$$
\end{lemma}
\bp Let $C$ be the expression of the second line. Computing $C$ with Lemma~\ref{lemlkfundasur} yields
$$C=\sum_{(i,k) \in \underline{g}^2}lk(z_k^+,z_{\overline{i}})\langle z_i, z_{\overline{i}}\rangle^2
lk(z_{\overline{k}}^+,z_i)\langle z_{\overline{i}}, z_i\rangle \langle  z_k,z_{\overline{k}} \rangle.$$
\eop

\noindent {\sc Proof of Lemma~\ref{lemvarsurlk}:}
Use Lemmas~\ref{lemcpcp}, \ref{lemcpqi}, \ref{lembilqiqi}, \ref{lemlkfundasurgg}
and \ref{lemexplamp} above to compute
$$\begin{array}{lll}s_{\ell}(\CD^{\prime},\matchingo)-s_{\ell}(\CD,\matchingo)&=
&-\sum_{(c,d) \in (\CaC^+)^2}\CJ_{j(c)i(c)}\sigma(c)\CJ_{j(d)i(d)}\sigma(d)\\
&&
- \sum_{(i,j,k,s) \in \underline{g}^4}\CJ_{ji}\CJ_{s\overline{i}}\langle z_i,z_{\overline{i}} \rangle \langle  z_k,z_{\overline{k}} \rangle\langle  u_k,\beta_j \rangle\langle u_{\overline{k}},\beta_s \rangle\\
&=& -\sum_{(i,j,r,s) \in \underline{g}^4}\CJ_{ji} \langle  u_i,\beta_j \rangle \CJ_{sr}\langle  u_r,\beta_s \rangle -2\lambda^{\prime}\\
&=& -g(\Sigma)^2 -2\lambda^{\prime}.
\end{array}$$
\eop

\begin{lemma}
\label{lemexplampp}
Set $$\lambda^{\prime}_+(\Sigma)=\sum_{(i,r) \in \underline{g}^2}lk(z_r^+,z_i)
lk(z_{\overline{i}}^+,z_{\overline{r}})\langle z_i, z_{\overline{i}}\rangle\langle z_r , z_{\overline{r}}\rangle.$$
Then
$$\lambda^{\prime}_+(\Sigma)=-\sum_{(i,j,r,s) \in \underline{g}^4}\CJ_{ji}\CJ_{sr} \langle  u_i,\beta_s \rangle\langle u_r, \beta_j \rangle 
=2\lambda^{\prime} -g(\Sigma).$$
\end{lemma}
\bp
Using Lemma~\ref{lemlkfundasur}, we get
$$\begin{array}{lll}\sum_{(i,j,r,s) \in \underline{g}^4}\left(-\CJ_{ji}\CJ_{sr} \langle  u_i,\beta_s \rangle\langle u_r, \beta_j \rangle \right)&=&-\sum_{(i,r) \in \underline{g}^2}lk(z_r^+,z_{\overline{i}})\langle z_i, z_{\overline{i}}\rangle
lk(z_i^+,z_{\overline{r}})\langle z_r , z_{\overline{r}}\rangle\\
&=&\sum_{(i,r) \in \underline{g}^2}lk(z_r^+,z_i)
lk(z_{\overline{i}}^+,z_{\overline{r}})\langle z_i, z_{\overline{i}}\rangle\langle z_r , z_{\overline{r}}\rangle\\
&=&\sum_{(i,r) \in \underline{g}^2}lk(z_i^+,z_r)
lk(z_{\overline{i}}^+,z_{\overline{r}})\langle z_i, z_{\overline{i}}\rangle\langle z_r , z_{\overline{r}}\rangle \\
&&-\sum_{(i,r) \in \underline{g}^2}\langle z_i, z_r\rangle
lk(z_{\overline{i}}^+,z_{\overline{r}})\langle z_i, z_{\overline{i}}\rangle\langle z_r , z_{\overline{r}}\rangle\\
&=&2\lambda^{\prime}+\sum_{i \in \underline{g}}
lk(z_{\overline{i}}^+,z_{i})\langle z_i, z_{\overline{i}}\rangle=2\lambda^{\prime}-g(\Sigma).
\end{array}$$
\eop

\noindent {\sc Proof of Lemma~\ref{lemvarsurlamb}:}
According to Lemmas~\ref{lemcacprime}, \ref{lemexplamp} and \ref{lemexplampp},
$$\sum_{c \in \CaC^{\prime} \setminus \CaC}\CJ_{j(c)i(c)}\sigma(c)\ell^{\prime}(c,c)=-\lambda^{\prime}_+(\Sigma)-g(\Sigma)-2\lambda^{\prime}=-4\lambda^{\prime}.$$
Therefore
$$\sum_{c \in \CaC}\CJ_{j(c)i(c)}\sigma(c)\ell(c,c)-\sum_{c \in \CaC^{\prime}}\CJ_{j(c)i(c)}\sigma(c)\ell^{\prime}(c,c)=4\lambda^{\prime}+ \sum_{(i,j) \in \underline{g}^2} \CJ_{ji}\langle u_i,\beta_j\rangle=4\lambda^{\prime}+g(\Sigma)$$
according to Lemmas \ref{lemcpcp} and \ref{lemlkfundasurgg}.
Using Lemmas~\ref{lemcpcp}, \ref{lemcpqi} and \ref{lembilqiqi} again, we get
$$\begin{array}{lll}\ell_2(\CD^{\prime})-\ell_2(\CD)&=&-\sum_{(i,j,k,s) \in \underline{g}^4}\CJ_{j\overline{i}} \CJ_{si}\langle z_i,z_{\overline{i}} \rangle \langle  z_k,z_{\overline{k}} \rangle\langle  u_k,\beta_j \rangle\langle u_{\overline{k}},\beta_s \rangle\\
   &&-\sum_{(i,j,k,s) \in \underline{g}^4}\CJ_{jr} \CJ_{si}\langle  u_i,\beta_j \rangle \langle  u_r,\beta_s \rangle
\\&&
+4\lambda^{\prime}+g(\Sigma)
\\&=& 2\lambda^{\prime} + \lambda^{\prime}_+(\Sigma) +4\lambda^{\prime}+g(\Sigma)=8\lambda^{\prime}
\end{array}$$ thanks to Lemma~\ref{lemexplampp}.
\eop

Finally, we identify $\lambda^{\prime}$ to $\frac12 \Delta_K^{\prime \prime}(1)$ where 
$$\Delta_K(t)=t^{-g(\Sigma)}\mbox{det}\left(\left[t lk(z_r^+,z_s) -lk(z_s^+,z_r)\right]_{(r,s) \in \underline{g}^2}\right)$$ denotes the Alexander polynomial of $K$.

\begin{lemma}
\label{lemalexand}
$$\frac12 \Delta_K^{\prime \prime}(1)=\lambda^{\prime}.$$
\end{lemma}
\bp 
Note $t lk(z_r^+,z_s) -lk(z_s^+,z_r) = (t-1)lk(z_r^+,z_s) + \langle z_r, z_s \rangle$.
$$\Delta_K(t) = t^{-g(\Sigma)} + t^{-g(\Sigma)} (t-1)\sum_{i \in \underline{g}}lk(z_i^+,z_{\overline{i}})\langle z_i, z_{\overline{i}} \rangle + t^{-g(\Sigma)} (t-1)^2 A+ O(t-1)^3$$
where $\sum_{i \in \underline{g}}lk(z_i^+,z_{\overline{i}})\langle z_i, z_{\overline{i}} \rangle=g(\Sigma)$ (see Lemma~\ref{lemlkfundasurgg})
and 
$$\begin{array}{ll}A&=\sum_{\{i,r\} \subset \underline{g}} \langle z_i, z_{\overline{i}}\rangle \langle z_r, z_{\overline{r}}\rangle \left( lk(z_i^+, z_{\overline{i}})lk(z_r^+, z_{\overline{r}}) - lk(z_i^+, z_{\overline{r}})lk(z_r^+, z_{\overline{i}})\right)\\
&= \frac12\sum_{(i,r) \in \underline{g}^2} \langle z_i, z_{\overline{i}}\rangle \langle z_r, z_{\overline{r}}\rangle \left( lk(z_i^+, z_{\overline{i}})lk(z_r^+, z_{\overline{r}}) - lk(z_i^+, z_{\overline{r}})lk(z_r^+, z_{\overline{i}})\right)\\
&= \frac{g(\Sigma)^2}2 + \frac12\sum_{(i,r) \in \underline{g}^2}lk(z_r^+,z_i)
lk(z_{\overline{i}}^+,z_{\overline{r}})\langle z_i, z_{\overline{i}}\rangle\langle z_r , z_{\overline{r}}\rangle
\\&=\frac12 \left(g(\Sigma)^2 + \lambda^{\prime}_+(\Sigma) \right)
 \end{array}$$ thanks to Lemma~\ref{lemexplampp}.\\
$\Delta_K^{\prime}(t) = -g(\Sigma)t^{-g(\Sigma)-1} + g(\Sigma)(t^{-g(\Sigma)}-g(\Sigma)t^{-g(\Sigma)-1}(t-1))+ 2t^{-g(\Sigma)} (t-1) A+ O(t-1)^2$.
$$\Delta_K^{\prime \prime}(1)=g(\Sigma)(g(\Sigma)+1 - 2 g(\Sigma)) +2 A=g(\Sigma)+\lambda^{\prime}_+(\Sigma)=2\lambda^{\prime}$$
according to Lemma~\ref{lemexplampp}.
\eop

\def\cprime{$'$}
\providecommand{\bysame}{\leavevmode ---\ }
\providecommand{\og}{``}
\providecommand{\fg}{''}
\providecommand{\smfandname}{et}
\providecommand{\smfedsname}{\'eds.}
\providecommand{\smfedname}{\'ed.}
\providecommand{\smfmastersthesisname}{M\'emoire}
\providecommand{\smfphdthesisname}{Th\`ese}

\newpage
\addcontentsline{toc}{section}{Index of notations}
\printindex


\begin{thebibliography}{LMO98}

\bibitem[AM90]{akmc}
{\scshape S.~Akbulut {\normalfont \smfandname} J.~D. McCarthy} --
  \emph{Casson's invariant for oriented homology {$3$}-spheres}, Mathematical
  Notes, vol.~36, Princeton University Press, Princeton, NJ, 1990, An
  exposition.

\bibitem[Cer70]{Cerfpseudo}
{\scshape J.~Cerf} -- {\og La stratification naturelle des espaces de fonctions
  diff\'erentiables r\'eelles et le th\'eor\`eme de la pseudo-isotopie\fg},
  \emph{Inst. Hautes \'Etudes Sci. Publ. Math.} (1970), no.~39, p.~5--173.

\bibitem[GH11]{GH}
{\scshape V.~Gripp {\normalfont \smfandname} Y.~Huang} -- {\og An absolute
  grading on {H}eegaard {F}loer homology by homotopy classes of oriented
  2-plane fields\fg}, arXiv:1112.0290v2, 2011.

\bibitem[GM92]{gm}
{\scshape L.~Guillou {\normalfont \smfandname} A.~Marin} -- {\og Notes sur
  l'invariant de {C}asson des sph\`eres d'homologie de dimension trois\fg},
  \emph{Enseign. Math. (2)} \textbf{38} (1992), no.~3-4, p.~233--290, With an
  appendix by Christine Lescop.

\bibitem[Gom98]{Go}
{\scshape R.~E. Gompf} -- {\og Handlebody construction of {S}tein surfaces\fg},
  \emph{Ann. of Math. (2)} \textbf{148} (1998), no.~2, p.~619--693.

\bibitem[Hir73]{hirzebruchEM}
{\scshape F.~E.~P. Hirzebruch} -- {\og Hilbert modular surfaces\fg},
  \emph{Enseignement Math. (2)} \textbf{19} (1973), p.~183--281.

\bibitem[KM99]{km}
{\scshape R.~Kirby {\normalfont \smfandname} P.~Melvin} -- {\og Canonical
  framings for {$3$}-manifolds\fg}, \emph{Proceedings of 6th {G}\"okova
  {G}eometry-{T}opology {C}onference}, vol.~23, Turkish J. Math., no.~1, 1999,
  p.~89--115.

\bibitem[Kon94]{ko}
{\scshape M.~Kontsevich} -- {\og Feynman diagrams and low-dimensional
  topology\fg}, First {E}uropean {C}ongress of {M}athematics, {V}ol.\ {II}
  ({P}aris, 1992), Progr. Math., vol. 120, Birkh\"auser, Basel, 1994,
  p.~97--121.

\bibitem[KT99]{kt}
{\scshape G.~Kuperberg {\normalfont \smfandname} D.~Thurston} -- {\og
  Perturbative 3--manifold invariants by cut-and-paste topology\fg},
  math.GT/9912167, 1999.

\bibitem[Les04a]{lesconst}
{\scshape C.~Lescop} -- {\og On the {K}ontsevich-{K}uperberg-{T}hurston
  construction of a configuration-space invariant for rational homology
  3-spheres\fg}, math.GT/0411088, 2004.

\bibitem[Les04b]{lessumgen}
\bysame , {\og Splitting formulae for the {K}ontsevich-{K}uperberg-{T}hurston
  invariant of rational homology 3-spheres\fg}, math.GT/0411431, 2004.

\bibitem[Les12a]{lesHC}
\bysame , {\og {A formula for the $\Theta$-invariant from {H}eegaard
  diagrams}\fg}, arXiv:1209.3219v2, 2012.

\bibitem[Les12b]{lesmek}
\bysame , {\og Introduction to finite type invariants of knots and
  $3$--manifolds\fg}, {M}ekn\`es research school 2012
  http://www-fourier.ujf-grenoble.fr/$\tilde{\;}$ lescop/preprints/meknes.pdf,
  2012.

\bibitem[Les13]{lescomb}
\bysame , {\og On homotopy invariants of combings of $3$--manifolds\fg},
  arXiv:1209.2785v2, 2013.

\bibitem[LMO98]{lmo}
{\scshape T.~T.~Q. Le, J.~Murakami {\normalfont \smfandname} T.~Ohtsuki} --
  {\og On a universal perturbative invariant of {$3$}-manifolds\fg},
  \emph{Topology} \textbf{37} (1998), no.~3, p.~539--574.

\bibitem[Mar88]{mar}
{\scshape A.~Marin} -- {\og Un nouvel invariant pour les sph\`eres d'homologie
  de dimension trois (d'apr\`es {C}asson)\fg}, \emph{Ast\'erisque} (1988),
  no.~161-162, p.~Exp.\ No.\ 693, 4, 151--164 (1989), S{\'e}minaire Bourbaki,
  Vol. 1987/88.

\bibitem[OS04]{OS1}
{\scshape P.~Ozsv{\'a}th {\normalfont \smfandname} Z.~Szab{\'o}} -- {\og
  Holomorphic disks and topological invariants for closed three-manifolds\fg},
  \emph{Ann. of Math. (2)} \textbf{159} (2004), no.~3, p.~1027--1158.

\bibitem[Sie80]{SiebenmannRS}
{\scshape L.~C. Siebenmann} -- {\og Les bisections expliquent le th\'eor\`eme
  de {R}eidemeister-{S}inger\fg}, {P}r\'epublications math\'ematiques d'{O}rsay
  Volume 16, Volume 80 de Pr\'epublications / Universit\'e de Paris-Sud,
  D\'epartement de Math\'ematiques, 1980.

\bibitem[Suz77]{Suzuki}
{\scshape S.~Suzuki} -- {\og On homeomorphisms of a 3-dimensional
  handlebody\fg}, \emph{Canad. J. Math.} \textbf{29} (1977), no.~1,
  p.~111--124.

\bibitem[Wal92]{wal}
{\scshape K.~Walker} -- \emph{An extension of {C}asson's invariant}, Annals of
  Mathematics Studies, vol. 126, Princeton University Press, Princeton, NJ,
  1992.

\bibitem[Wit89]{witten}
{\scshape E.~Witten} -- {\og Quantum field theory and the {J}ones
  polynomial\fg}, \emph{Comm. Math. Phys.} \textbf{121} (1989), no.~3,
  p.~351--399.

\end{thebibliography}
\end{document}